\documentclass[a4paper,12pt]{amsart}
\usepackage[frenchb]{babel}
\usepackage{pdfsync}
\usepackage[T1]{fontenc}
\usepackage{dsfont}
\usepackage{amsmath}
\usepackage{amssymb}
\usepackage{amsxtra}
\usepackage{amscd}
\usepackage{amsthm}
\usepackage{amsfonts}
\usepackage{eucal}
\usepackage[all]{xy}
\usepackage{xypic}
\usepackage{graphicx}
\usepackage{comment}
\usepackage{epsfig}
\usepackage{psfrag}
\usepackage{mathrsfs}
\usepackage{amscd}
\usepackage{rotating}
\usepackage{lscape}
\usepackage{amsbsy}
\usepackage{verbatim}
\usepackage{moreverb}
\usepackage{url}

\newcommand{\nc}{\newcommand}
\nc{\renc}{\renewcommand}

\nc\restr[2]{{
  \left.\kern-\nulldelimiterspace 
  #1 
  \vphantom{\big|} 
  \right|_{#2} 
  }}

\newtheorem{thm}{Théorème}[section]
\newtheorem{prop}[thm]{Proposition}
\newtheorem{lem}[thm]{Lemme}

\theoremstyle{definition}
\newtheorem{defi}[thm]{Définition}
\newtheorem{def-prop}[thm]{Définition-Proposition}
\newtheorem{rem}[thm]{Remarque}

\newtheorem{rappel}[thm]{Rappel}

\newtheorem{construction}[thm]{Construction}
\newtheorem{notation}[thm]{Notation}

\numberwithin{equation}{section}

\renc{\sec}{\section}
\nc{\ssec}{\subsection}
\nc{\sssec}{\subsubsection}

\nc{\thmref}[1]{théorème~\ref{#1}}
\nc{\secref}[1]{paragraphe~\ref{#1}}
\nc{\lemref}[1]{lemme~\ref{#1}}
\nc{\defiref}[1]{définition~\ref{#1}}
\nc{\propref}[1]{proposition~\ref{#1}}
\nc{\corref}[1]{corollaire~\ref{#1}}
\nc{\constructionref}[1]{construction~\ref{#1}}
\nc{\conjref}[1]{conjecture~\ref{#1}}
\nc{\remref}[1]{remarque~\ref{#1}}
\nc{\rappelref}[1]{rappel~\ref{#1}}
\nc{\questref}[1]{question~\ref{#1}}
\nc\Omegasour{\hbox{$\buildrel\smile\over{\vrule height 6pt depth 0pt width 0pt \smash \Omega}$}}

\nc{\on}{\operatorname}

\nc\wt{\widetilde}
\nc\wh{\widehat}
\nc\ol{\ov}
\nc{\oc}[1]{{\overset{\circ}{#1}}}
\nc{\ov}[1]{{\overline{#1}}}
\nc{\isor}[1]{{\xrightarrow[\raisebox{0.25 em}{\smash{\ensuremath{\sim}}}]{#1}}}
\nc{\isol}[1]{{\xleftarrow[\raisebox{0.25 em}{\smash{\ensuremath{\sim}}}]{#1}}}
\nc{\modmod}{/ \! \! /}

\nc{\mc}{\mathcal}
\nc{\mf}{\mathfrak}
\nc{\mr}{\mathrm}
\nc{\mb}{\mathbb}
\nc{\mbf}{\mathbf}

\nc{\R}{{\mathbb R}}
\nc{\Z}{{\mathbb Z}}
\nc{\N}{{\mathbb N}}
\nc{\C}{{\mathbb C}}
\nc{\Q}{{\mathbb Q}}

\nc{\Fq}{{\mathbb F}_q}
\nc{\Fl}{{\mathbb F}_\ell}
\nc{\Fqbar}{\ol{{\mathbb F}_q}}
\nc{\Flbar}{\ol{{\mathbb F}_\ell}}
\nc{\Zl}{{\mathbb Z}_\ell}
\nc{\Zlbar}{\ol{{\mathbb Z}_\ell}}
\nc{\Ql}{{\mathbb Q}_\ell}
\nc{\Qlbar}{\ol{{\mathbb Q}_\ell}}
\nc{\hl}{\overset{\leftarrow}h{}}
\nc{\hr}{\overset{\rightarrow}h{}}
\nc{\Gr}{{\on{Gr}}}
\nc{\Hecke}{\on{Hecke}}
 \nc{\Hom}{\on{Hom}}
 \nc{\Coker}{\on{Coker}}
 \nc{\Ker}{\on{Ker}}
 \nc{\Lie}{\on{Lie}}
\nc{\Loc}{\on{Loc}}
\nc{\Pic}{\on{Pic}}
\nc{\Bun}{\on{Bun}}
\nc{\IC}{\on{IC}}
\nc{\Aut}{\on{Aut}}
\nc{\Perv}{\on{Perv}}
\nc{\pos}{{\on{pos}}}
\nc{\Sym}{\on{Sym}}

\nc{\ta} {{}^\tau}
\nc {\tu}[1]{{}^{\tau^{#1}}\!}
\nc{\tav} {{}^\sigma}
\nc {\tuv}[1]{{}^{\sigma^{#1}}\!}

\nc{\Chr}{\mr{Cht}\mc R}
\nc{\Id}{\on{Id}}
\nc{\Fil}{\on{Fil}}
\nc{\pr}{\on{pr}}
\nc{\Res}{\on{Res}}
\nc{\cusp}{\on{cusp}}
\nc{\Frob}{\on{Frob}}
\nc{\diag}{\Delta}
\nc{\gr}{\on{gr}}
\nc{\Inj}{\on{Inj}}
\nc{\Bl}{\on{Bl}}
\nc{\dem}{\noindent {\bf Démonstration. }}
\nc{\cqfd}{{\ }\hfill $\square$ \vskip 1mm}
\nc{\s}[1]{\langle #1 \rangle}
\nc{\Cht}{\on{Cht}}
\nc{\isom}{\overset {\thicksim}{\to}}
\nc{\sm}{\smallsetminus}

\nc\Spf{\mathop{\mathrm {Spf }}}

\nc \lo{\mathrm{loc}}

\begin{document}

\title[ $G$-chtoucas restreints.  ]{Chtoucas restreints pour les groupes réductifs et paramétrisation de  Langlands locale.}
 \author{Alain Genestier et Vincent Lafforgue}

\dedicatory{A la m\'emoire de  G\'erard Laumon.} 

\address{Alain Genestier : 
Institut Elie Cartan, 
Université de Lorraine, 
B.P. 70239, 54506 Vandoeuvre-lès-Nancy,
France} 

\address{Vincent Lafforgue: CNRS et Institut de Math\'ematiques de Jussieu, Universit\'e Paris  Cit\'e.}
\date{\today}
\maketitle

 Soit $K$ un corps local d'égales caractéristiques, $\mc O_{K}$ son anneau d'entiers, $\pi_{K}$ une uniformisante et $\mf k$ son corps résiduel. Soit $\ell$ un nombre premier ne divisant pas $\sharp \mf k$.  Soit   $E$  une extension finie de 
    $\Ql$ contenant une racine carrée de $\sharp \mf k$. On note $\mc O_{E}$  son anneau d'entiers et $\lambda_{E}$ une uniformisante.

Soit $G$ un groupe réductif sur $K$. 
Sauf dans le paragraphe \ref{cas-non-deploye},  on le suppose déployé car cela simplifie les énoncés et les démonstrations. Le cas non déployé, qui  ne nécessite pas d'idée nouvelle mais entraîne des notations plus compliquées, fait l'objet du paragraphe \ref{cas-non-deploye}.  

  Un sous-groupe compact ouvert $U$ de $G(K)$ est dit d'ordre premier à  $\ell$ si c'est le cas de tous ses quotients finis. Une structure entière sur une représentation $E$-linéaire lisse admissible $H$ de $G(K)$ est la donnée  pour tout 
sous-groupe compact ouvert $U$  d'ordre premier à  $\ell$ d'une $\mc O_{E}$-structure $H^{U}_{\mc O_{E}}$ sur  le  $E$-espace vectoriel de dimension finie 
$H^{U}$,  de sorte que 
\begin{itemize}
\item pour tout $U$, $H^{U}_{\mc O_{E}}$ est stable par 
$C_{c}(U\backslash G(K)/U,\mc O_{E})$ (dont l'action est  normalisée pour que 
$\mathds 1_{U}$ agisse par l'identité), 
\item 
pour tous $U'\subset U$, on a $H^{U}_{\mc O_{E}}=H^{U'}_{\mc O_{E}}\cap H^{U}$. 
\end{itemize}
On dira qu'une représentation lisse  admissible est   entière si elle admet une structure entière.

  On note $\wh G$ le groupe dual de Langlands de $G$, considéré comme  un groupe réductif connexe déployé défini sur $\mc O_{E}$ (et donc a fortiori sur $E$).

Pour énoncer la compatibilité local-global, on rappelle le cadre de \cite{coh}. 
 Soit  $\Fq$ un corps fini,  $X$ une courbe projective lisse géométriquement irréductible  sur $\Fq$, $F_X$ son corps des fonctions, $\mb A$ ses adèles, 
 $N$ un sous-schéma fini de $X$ et $\Xi$ un sous-groupe discret cocompact 
  de $Z_{G}(F_X)\backslash Z_{G}(\mb A)$ (où $Z_{G}$ désigne le centre de $G$).  
 On rappelle que   \cite{coh}  fournit 
  une décomposition canonique de  
   $C_{c}(K_{N}\backslash G(\mb A)/K_{N},\Qlbar)$-modules 
 \begin{gather}\label{intro1-dec-canonique}
 C_{c}^{\mr{cusp}}(\Bun_{G,N}(\Fq)/\Xi,\Qlbar)=\bigoplus_{\sigma}
 \mf H_{\sigma},\end{gather}
 où la somme directe dans le membre de droite est indexée par des paramètres de Langlands globaux, c'est-à-dire des classes de  $\wh G(\Qlbar)$-conjugaison de  morphismes 
       $\sigma:\on{Gal}(\ov{ F_X}/F_X)\to \wh G(\Qlbar)$ 
       définis sur une extension finie de  $\Ql$, continus,    semi-simples et non ramifiés sur   $X\sm N$.

On montre dans  cet article   la paramétrisation  de Langlands locale à semi-simplification près, et la compatibilité local-global. Plus précisément on montre le théorème suivant. 

\begin{thm} \label{thm-intro}
Il existe une application 
\begin{gather}\label{param-local-sigma-pi}\pi\mapsto \sigma_{\pi}\end{gather}
\begin{itemize}
\item []
de l'ensemble des classes d'isomorphismes de représentations lisses admissibles et irréductibles de $G(K)$ définies sur une extension finie 
de $\Ql$ et entières, 
\item [] vers l'ensemble des classes d'isomorphismes de paramètres de Langlands locaux semi-simples, c'est-à-dire 
  des classes de  $\wh G(\Qlbar)$-conjugaison de  morphismes 
       $\sigma:\on{Gal}(\ov K/K)\to \wh G(\Qlbar)$ 
       définis sur une extension finie de  $\Ql$, continus et   semi-simples
\end{itemize}
 qui est déterminée de manière unique par les deux propriétés suivantes: 
\begin{itemize}
\item [] a)   $\sigma_{\pi}$ ne dépend que du caractère par lequel le centre de Bernstein agit sur $\pi$, et en dépend ``algébriquement'' au sens du \thmref{thm-mfz} ci-dessous, 
\item [] b) cette application 
 est compatible avec la paramétrisation globale construite dans \cite{coh}, au sens suivant. 
 \end{itemize}

  Soient  $\Fq,  X,  N,  \Xi$  comme ci-dessus. 
  Alors pour toute représentation $\pi=\bigotimes \pi_{v}$ de $G(\mb A)$ telle que $\pi^{K_{N}}$ apparaisse dans $\mf H_{\sigma}$ (dans \eqref{intro1-dec-canonique} ci-dessus), pour toute place $v$ de $X$ on a égalité entre 
  \begin{itemize}
  \item   le paramètre local $\sigma_{\pi_{v}}$ obtenu en appliquant 
  \eqref{param-local-sigma-pi}  avec   $K$ égal au complété 
 $F_{v}$ de $F_X$ en $v$,   
  \item   le semi-simplifié de la restriction  de $\sigma$ à $\on{Gal}(\ov{F_{v}}/F_{v})$.  
 \end{itemize}
 
 De plus cette application $\pi\mapsto \sigma_{\pi}$ 
 s'étend de fa\c con unique en une application 
 \begin{itemize}
\item []
de l'ensemble des classes d'isomorphismes de représentations lisses admissibles et irréductibles de $G(K)$ définies sur une extension finie 
de $\Ql$ (pas nécessairement entières)  
\item [] vers l'ensemble des classes d'isomorphismes de paramètres de Weil locaux semi-simples, c'est-à-dire 
  des classes de  $\wh G(\Qlbar)$-conjugaison de  morphismes 
       $\sigma:\on{Weil}(\ov K/K)\to \wh G(\Qlbar)$ 
       définis sur une extension finie de  $\Ql$, continus et   semi-simples
\end{itemize}
vérifiant a) ci-dessus. 

Cette application est compatible 
  à l'induction parabolique au sens suivant. Si $P$ est un sous-groupe parabolique de $G$, de quotient de Levi $M$, et si $\tau $ est une représentation lisse, admissible et  irréductible de $M(K)$ définie sur une extension finie 
de $\Ql$, et $\pi$ est un sous-quotient irréductible de la représentation induite compacte $\on{Ind}_{P(K)}^{G(K)}\tau$ (avec la normalisation unitaire), 
alors $\sigma_{\pi}$ est conjugué à la composée 
$ \on{Weil}(\ov K/K)\xrightarrow{\sigma_{\tau}}  \wh M(\Qlbar) \to  \wh G(\Qlbar)$. 

Enfin elle est compatible aux cas triviaux de  fonctorialité. 
Soit $G'$ un  groupe réductif déployé sur $K$ et $\Upsilon:G\to G'$ un morphisme de groupes sur $K$ dont  l'image est un sous-groupe distingué de $G'$. On en déduit  ${}^{L}\Upsilon: \wh  {G'}\to \wh  G$. Alors pour toute représentation lisse admissible irréductible $\pi$ de $G'(K)$, 
le centre de Bernstein de $G(K)$ agit sur $\pi$ par un caractère dont le paramètre local associé par \eqref{param-local-sigma-pi} et a) ci-dessus  est  
${}^{L}\Upsilon\circ \sigma_{\pi}$. 
\end{thm}

On verra dans le paragraphe \ref{facteurs-gamma}
 qu'on peut en déduire une théorie complète des facteurs $\gamma$.  

Dans \cite{wenwei}, Wen-Wei Li a montré que  
le paramètre local d'une représentation et celui de sa contragrédiente 
se déduisent l'un de l'autre par une involution de Cartan du groupe dual.

Dans le cas où $G=GL_{r}$, ce théorème était déjà connu par Laumon-Rapoport-Stuhler
\cite{laumon-rapoport-stuhler} et le corollaire VII.5 de \cite{laurent-inventiones} pour la compatibilité local-global (qui dans le cas de $GL_{r}$ est vraie sans semi-simplification).
En particulier pour $GL_{1}$, $\pi\mapsto \sigma_{\pi}$ est la théorie du corps de classe local. Par ailleurs 
 la construction des facteurs $\gamma$ avait déjà été faite dans de nombreux cas par des méthodes de Langlands-Shahidi \cite{hl11, hl13b, hl13a, gana-lomeli, lomeli, lomeli-AIF, 
lomeli15, gan-lomeli, hl17}. 

Dans les  travaux    \cite{berkeley}, \cite{fargues-ICM}, \cite{fargues-scholze},  Fargues et  Scholze
obtiennent la correspondance locale pour les tous les corps locaux, par des m\'ethodes  différentes. 

Dans  \cite{li-huerta}, Li-Huerta montre (par une compatibilit\'e commune avec le global) que la construction ci-dessus co\"incide avec celle de Fargues et Scholze \cite{fargues-scholze} dans le cas d'\'egales caract\'eristiques. 

\begin{rem}
En appliquant la dernière propriété du théorème à $G'=G\times \mb G_m$ et $\Upsilon=(\Id, \chi)$ où $\chi:G\to  \mb G_m$ est un caractère, on obtient la compatibilité 
de $\pi\mapsto \sigma_{\pi}$ 
au produit tensoriel  de $\pi$ par un caractère. 
\end{rem}

\begin{rem} Dans le théorème précédent (ou au moins dans sa généralisation aux groupes non déployés qui sera expliquée dans le dernier paragraphe)  on ne peut pas espérer mieux que la compatibilité local-global à semi-simplification près. En effet 
cela se voit déjà dans le cas des algèbres à division, comme nous l'ont fait remarquer Guy Henniart et Ioan Badulescu: si $D$ est une algèbre à division de rang $r$, anisotrope en $v$, et $\pi=\bigotimes_{w}\pi_{w}$ est une forme automorphe pour $D^{\times}$ telle que $\pi_{v}$ soit la représentation triviale, et $\sigma$ est le  paramètre de Langlands  associé à $\pi$ par \cite{coh}, alors $\restr{\sigma}{\on{Gal}(\ov {F_{v}}/F_{v})}$  
\begin{itemize}
\item est semi-simple si $\pi$ est la représentation triviale
\item n'est pas semi-simple si $\pi$ correspond par Jacquet-Langlands
à une forme automorphe cuspidale pour $GL_{r}$ dont la composante en $v$ est 
la Steinberg. 
\end{itemize}
\end{rem}

 La preuve utilise les champs de chtoucas restreints, qui sont des  analogues en égales caractéristiques des groupes de Barsotti-Tate tronqués, à ceci près que les chtoucas restreints  peuvent avoir plusieurs pattes. On possède aussi en égales caractéristiques une notion de chtoucas locaux (analogues aux groupes $p$-divisibles, à ceci près qu'ils peuvent avoir plusieurs pattes). Nous avons utilisé dans cet article les chtoucas restreints et non les chtoucas locaux car nous devons appliquer la théorie des  cycles proches sur une base générale, due à Deligne, Laumon, Gabber, Orgogozo, et elle n'a été écrite que dans un cadre algébrique.


Nous remercions Ioan Badulescu, Jean-Benoît Bost, Jean-Fran\c cois Dat, Arnaud Eteve, 
Laurent Fargues, 
Dennis Gaitsgory, Michael Harris, Guy Henniart, Laurent Lafforgue,  Wen-Wei Li,  Luis Lomel\'{\i}, Fabrice Orgogozo, Peter Scholze, Jared Weinstein, Cong Xue, Weizhe Zheng  et  Xinwen Zhu pour des discussions ou des remarques liées à  cet article. 

Nous sommes tous deux immensément redevables \`a G\'erard Laumon. Au fil des ann\'ees, il nous a consacr\'e son temps sans compter, \`a travers de nombreuses discussions (y compris au sujet de cet article). Il a laiss\'e une empreinte consid\'erable sur notre manière de concevoir les Math\'ematiques.

\section{Enoncé avec le centre de Bernstein et premières réductions}
\label{para-enonce}

Soit $U$ un sous-groupe ouvert   compact de $G(K)$  d'ordre premier à  $\ell$. 
On normalise la structure d'algèbre sur $C_{c}(U\backslash G(K)/U,E)$ pour que $\mathds 1_{U}$ soit une unité. Alors $C_{c}(U\backslash G(K)/U,\mc O_{E})$ est une $\mc  O_{E}$-algèbre. 

On note $\mf Z_{U}$ le centre de $C_{c}(U\backslash G(K)/U,E)$, $\mf Z_{U,\mc O_{E}}$ le centre de $C_{c}(U\backslash G(K)/U,\mc O_{E})$ (qui est donc une sous-alg\`ebre de  $\mf Z_{U}$) et $\wh{\mf Z_{U,\mc O_{E}}}$ son complété $\ell$-adique, muni de la topologie $\ell$-adique. 
Pour $U'\subset U$ on a des morphismes naturels  de $\mc  O_{E}$-algèbres 
\begin{gather}\label{inclusions-Z}\mf Z_{U'}\to \mf Z_{U},
 \mf Z_{U',\mc O_{E}}\to \mf Z_{U,\mc O_{E}}\text{ et }
 \wh{\mf Z_{U',\mc O_{E}}}\to \wh{\mf Z_{U,\mc O_{E}}} \end{gather}
donnés par $f\mapsto f\star \mbf 1_{U}$, et on note 
$$\mf Z=\varprojlim_{U} \mf Z_{U}, 
\mf Z_{\mc O_{E}}=\varprojlim_{U}\mf Z_{U,\mc O_{E}}\text{ et } 
 \wh{\mf Z_{\mc O_{E}}}= \varprojlim_{U}\wh{\mf Z_{U,\mc O_{E}}}.$$
 
 On munit $ \wh{\mf Z_{\mc O_{E}}}$ de la topologie limite projective des topologies $\ell$-adiques sur $\wh{\mf Z_{U,\mc O_{E}}}$. 
Toute représentation $E$-linéaire lisse  admissible irréductible $\rho$ de $G(K)$ fournit des caractères 
$\chi_{U,\rho}:\mf Z_{U}\to E$ pour $U$ assez petit, et ceux-ci sont compatibles avec les morphismes \eqref{inclusions-Z},  donc ils déterminent un caractère   
$\mf Z \to E$. Si de plus $\rho$  est  entière, elle  fournit des caractères 
$\chi_{U,\rho}:\wh{\mf Z_{U,\mc O_{E}}}\to \mc O_{E}$  pour $U$ assez petit, 
et ceux-ci sont compatibles avec les morphismes \eqref{inclusions-Z},  donc ils déterminent un caractère 
$\wh{\mf Z_{\mc O_{E}}}\to \mc O_{E}$. On notera $\mf Z^{G}$ et    
 $\wh{\mf Z_{\mc O_{E}}^{G}}$ quand on voudra spécifier $G$. 

Le \thmref{thm-intro} va résulter du théorème suivant. Soit $I$ un ensemble fini. On considère 
le quotient grossier $\wh G \backslash \wh G^{I}/\wh G$ (où les actions de $\wh G$ à gauche et à droite se font  par multiplication diagonale) comme un schéma  défini  sur $\mc O_{E}$. 
Soient   $f$ 
une fonction sur le quotient grossier $\wh G \backslash \wh G ^{I}/\wh G$ définie sur $\mc O_{E}$ et $(\gamma_{i})_{i\in I}$ 
 un $I$-uplet  d'éléments de  $\on{Weil}(\ov K/K)$. 
On rappelle que dans le paragraphe 13 de  \cite{coh}  on a défini l'opérateur d'excursion
     $S_{I,f,(\gamma_{i})_{i\in I}}\in \on{End}(C_{c}^{\mr{cusp}}(\Bun_{G,N}(\Fq)/\Xi,
     \mc O_{E}))$. 
     Cong Xue a construit une extension de cet opérateur à 
     $C_{c}(\Bun_{G,N}(\Fq)/\Xi,
     \mc O_{E})$,   dans \cite{cong-coeff-O} (pour le cas déployé, et  \cite{cong-zorro} pour tous les groupes par une approche un peu différente).
     Plus précisément pour tout entier $s$, et pour toute fonction $f$ sur $\wh G \backslash \wh G ^{I}/\wh G$ définie sur $\mc O_{E}/\lambda_{E}^{s}\mc O_{E}$,  Cong Xue définit 
       \begin{gather}       \label{excursion-etendus}
     S_{I,f,(\gamma_{i})_{i\in I}}:C_{c}(\Bun_{G,N}(\Fq)/\Xi,
     \mc O_{E}/\lambda_{E}^{s}\mc O_{E}))\to C_{c}(\Bun_{G,N}(\Fq)/\Xi,
     \mc O_{E}/\lambda_{E}^{s}\mc O_{E}). \end{gather}
 On note   $\ov{\eta_X}=\on{Spec}(\ov F)$ et comme dans   \cite{coh}, on choisit un point géométrique $\ov{\eta_X^{I}}$ au-dessus du point générique $\eta_X^{I}$ de $X^{I}$ et une flèche de spécialisation 
     $\mf{sp}: \ov{\eta_X^{I}}\to \Delta(\ov{\eta_X})$.  
     On fixe un plongement $\ov F\subset \ov K$, ce qui identifie $\ov K$ \`a une cl\^oture alg\'ebrique $\ov F_v$. 
     On choisit  
    un relèvement de $(\gamma_{i})_{i\in I}$ en un élément 
     $\delta\in \on{FWeil}(\eta^{I},\ov{\eta^{I}})$
     (voir le paragraphe 8 de  \cite{coh} pour la définition de ce groupe). On choisit 
    une représentation  de $\wh G ^{I}$ sur un 
    $ \mc O_{E}/\lambda_{E}^{s}\mc O_{E}$-module libre de type fini 
     $W$ 
et      $x\in W$ et  $\xi\in W^{*}$  invariants par l'action diagonale de $\wh G$ tels que 
      $ f$ soit égale au coefficient de matrice $(g_{i})_{i\in I}\mapsto \s{\xi, (g_{i})_{i\in I}\cdot x}$. 
     On rappelle alors que  \eqref{excursion-etendus} 
     est défini grâce à   \cite{cong-coeff-O} (pour le cas déployé, et  \cite{cong-zorro} pour tous les groupes par une approche un peu différente),
    comme la composée
        \begin{gather*}\label{diag-S}
    C_{c}(\Bun_{G,N}(\Fq)/\Xi, \mc O_{E}/\lambda_{E}^{s}\mc O_{E})
     \xrightarrow{ \mc C_{  x}^{\sharp}} 
   \varinjlim _{\mu} \restr{\mc H _{ N, I, W}^{0,\leq\mu, \mc O_{E}/\lambda_{E}^{s}\mc O_{E}}}{\Delta(\ov\eta)} \\ 
   \nonumber    \isor{\mf{sp}^{*}} 
    \varinjlim _{\mu}  \restr{\mc H _{ N, I, W}^{0,\leq\mu, \mc O_{E}/\lambda_{E}^{s}\mc O_{E}}}{ \ov{  \eta^{I}}}     \xrightarrow{\delta}  
     \varinjlim _{\mu}  \restr{\mc H _{ N, I, W}^{0,\leq\mu, \mc O_{E}/\lambda_{E}^{s}\mc O_{E}}}{ \ov{  \eta^{I}}} \\ 
      \nonumber    \isor{(\mf{sp}^{*})^{-1}}
       \varinjlim _{\mu} \restr{\mc H _{ N, I, W}^{0,\leq\mu, \mc O_{E}/\lambda_{E}^{s}\mc O_{E}}}{\Delta(\ov\eta)} 
    \xrightarrow{ \mc C_{  
\xi}^{\flat }}
    C_{c} (\Bun_{G,N}(\Fq)/\Xi, \mc O_{E}/\lambda_{E}^{s}\mc O_{E})
      \end{gather*}
    et est indépendant des choix de $\mf{sp}$ et $\delta$.

\begin{thm}\label{thm-mfz}
On possède 
\begin{itemize}
\item pour tout groupe $G$ comme ci-dessus, 
\item pour tout ensemble fini $I$, 
\item pour toute fonction $f$ sur le quotient grossier $\wh G \backslash \wh G ^{I}/\wh G$ définie sur $\mc O_{E}$, 
\item et pour tout $I$-uplet $(\gamma_{i})_{i\in I}$ d'éléments de    $\on{Weil}(\ov K/K)$, 
\end{itemize}
un élément $\mf z_{I,f,(\gamma_{i})_{i\in I}}\in  \wh{\mf Z_{\mc O_{E}}}$
(aussi noté $\mf z^{G}_{I,f,(\gamma_{i})_{i\in I}}\in  \wh{\mf Z^{G}_{\mc O_{E}}}$ quand on veut préciser $G$), qui vérifie  les propriétés suivantes:  
 \begin{itemize}
  \item [] (i) pour tout   $I$ et 
 $(\gamma_{i})_{i\in I}\in  \on{Weil}(\ov K/K)^{I}$, 
  $$f\mapsto 
  \mf z_{I,f,(\gamma_{i})_{i\in I}}$$ est un  morphisme 
  de $\mc O_{E}$-algèbres  commutatives  $\mc O(\wh G\backslash \wh G ^{I}/\wh G)\to  \wh{\mf Z_{\mc O_{E}}}$, 
  \item [] (ii) pour toute application 
  $\zeta:I\to J$, toute fonction  $f\in \mc O(\wh G\backslash \wh G ^{I}/\wh G)$ et  tout 
  $(\gamma_{j})_{j\in J}\in  \on{Weil}(\ov K/K)^{J}$, on a 
  $$\mf z_{J,f^{\zeta},(\gamma_{j})_{j\in J}}=\mf z_{I,f,(\gamma_{\zeta(i)})_{i\in I}}$$
   où $f^{\zeta}\in \mc O(\wh G\backslash \wh G ^{J}/\wh G)$ est définie par    $$f^{\zeta}((g_{j})_{j\in J})=f((g_{\zeta(i)})_{i\in I}),$$
   \item [] (iii) 
  pour tout   $f\in \mc O(\wh G\backslash \wh G ^{I}/\wh G)$
  et  $(\gamma_{i})_{i\in I},(\gamma'_{i})_{i\in I},(\gamma''_{i})_{i\in I}\in
    \on{Weil}(\ov K/K)^{I}$ on a     $$\mf z_{I\cup I\cup I,\wt f,(\gamma_{i})_{i\in I}\times (\gamma'_{i})_{i\in I}\times (\gamma''_{i})_{i\in I}}=
  \mf z_{I,f,(\gamma_{i}(\gamma'_{i})^{-1}\gamma''_{i})_{i\in I}}$$
   où  $I\cup I\cup I$ est une réunion disjointe et 
   $\wt f\in \mc O(\wh G\backslash \wh G ^{I\cup I\cup I}/\wh G)$ est définie  par  
   $$\wt f((g_{i})_{i\in I}\times (g'_{i})_{i\in I}\times (g''_{i})_{i\in I})=f((g_{i}(g'_{i})^{-1}g''_{i})_{i\in I}).$$
      \item [] (iv) pour tout $I$ et tout $f$, le morphisme 
      \begin{gather}\label{mor-excursion-f} \on{Weil}(\ov K/K)^{I}\to  \wh{\mf Z_{\mc O_{E}}}, \ \ (\gamma_{i})_{i\in I}\mapsto \mf z_{I,f,(\gamma_{i})_{i\in I}}\end{gather} est continu, c'est-à-dire que pour tout $U$ l'image de 
   $\mf z_{I,f,(\gamma_{i})_{i\in I}}$ dans     $ \wh{\mf Z_{U, \mc O_{E}}}$ (qui est muni de la topologie $\mc O_{E}$-adique) dépend continûment de $(\gamma_{i})_{i\in I}$, 
      \item []  (v) pour $\Fq,X, N, \Xi$ comme ci-dessus, pour toute place $v$ de $X$,  
     en prenant $K=F_{v}$, et pour tout plongement $\ov{ F_X}\subset \ov K$
     (d'où $\on{Weil}(\ov K/K)\subset \on{Gal}(\ov{F_X}/F_X)$)
     et pour $I,f,(\gamma_{i})_{i\in I}$ comme ci-dessus, 
     l'opérateur d'excursion $S_{I,f,(\gamma_{i})_{i\in I}}$   de \cite{coh,cong-finite}
     agit sur $  C_{c}(\Bun_{G,N}(\Fq)/\Xi,\mc O_{E})$ par multiplication par $\mf z_{I,f,(\gamma_{i})_{i\in I}}$, 
       \item [] (vi) si $M$ est un sous-groupe de Levi de $G$, $f$ et $(\gamma_{i})_{i\in I}$ sont comme ci-dessus, et $f^{M}$ désigne la fonction sur 
       $\wh M\backslash \wh M ^{I}/\wh M$ composée de $f$ avec le morphisme 
    $ \wh M\backslash \wh M ^{I}/\wh M\to  \wh G\backslash \wh G ^{I}/\wh G$, 
    alors pour toute représentation lisse admissible irréductible entière $\tau$ de $M(K)$, $\on{Ind}_{P(K)}^{G(K)}\tau$ 
    est une représentation lisse admissible entière de $G(K)$ et 
       $\mf z^{G}_{I,f,(\gamma_{i})_{i\in I}}$ agit sur  $\on{Ind}_{P(K)}^{G(K)}\tau$ par un scalaire égal au scalaire par lequel l'élément $\mf z^{M}_{I,f^{M},(\gamma_{i})_{i\in I}}\in \wh{\mf Z^{M}_{\mc O_{E}}}$ agit sur $\tau$. 
       
                \end{itemize}
      \end{thm}

\begin{rem}\label{li-huerta-s} En fait $\mf z_{I,f,(\gamma_{i})_{i\in I}}$ appartient au centre de Bernstein non compl\'et\'e. 
Cela r\'esulte  de \cite{li-huerta} et  \cite{fargues-scholze}, ou bien de la \remref{rem-supports1}. 
\end{rem}

\noindent{\bf Démonstration du \thmref{thm-intro}  en admettant le  \thmref{thm-mfz}.} On utilise les résultats de Richardson \cite{richardson}, plus précisément on applique la proposition 11.7 de \cite{coh} avec $\Gamma= \on{Weil}(\ov K/K)$ pour tout caractère entier du centre de Bernstein $ \wh{\mf Z_{\mc O_{E}}}$. 
En effet un tel caractère se factorise par $ \wh{\mf Z_{U, \mc O_{E}}}$  pour un certain $U$ et est automatiquement continu pour la topologie $\ell$-adique. 
 Tout $\sigma$ ainsi construit est défini sur $\mc O_{E}$ (quitte à élargir l'extension finie $E$ de $\Ql$), 
c'est-à-dire que pour tout $n$-uplet $(\gamma_{1},...,\gamma_{n})$ d'éléments de $\on{Weil}(\ov K/K)$, l'image de $(\sigma(\gamma_{1}),...,\sigma(\gamma_{n}))$ dans le quotient grossier  $\wh G ^{n} \modmod \wh G$ est définie sur $\mc O_{E}$. On  déduit de la construction de $\sigma$ fournie dans la démonstration de la proposition 11.7 de \cite{coh} que, quitte à augmenter $E$, on peut conjuguer $\sigma$ pour que son  image soit incluse dans 
$\wh G(\mc O_{E})$. Comme $\wh  G(\mc O_{E})$ est profini,   $\sigma$ 
  se factorise   par le morphisme $\on{Weil}(\ov K/K)\to \on{Gal}(\ov K/K)$. 
  On a donc obtenu l'application $\pi\mapsto \sigma_{\pi}$ dans le cas des représentations entières, et vérifié a) et b). La compatibilité avec l'induction parabolique dans le cas des représentations entières résulte du (vi) du \thmref{thm-mfz}. On en déduit que $\pi\mapsto \sigma_{\pi}$ s'étend à toutes les représentations (pas nécessairement entières). 
  En effet  toute supercuspidale ayant un caractère central entier admet une structure entière (et donc toute supercuspidale peut s'écrire comme le produit d'un caractère et d'une supercuspidale  admettant une structure entière, voir \cite{wenwei} pour plus de détails). Donc connaissant $\pi\mapsto \sigma_{\pi}$
  pour les supercuspidales entières on l'étend à toutes les supercuspidales, 
  puis à toutes les représentations en for\c cant la compatibilité à l'induction
  et comme cette compatibilité était connue pour les représentations entières, 
  l'application $\pi\mapsto \sigma_{\pi}$
 ainsi obtenue étend celle déjà construite dans le cas entier. 
  L'unicité de $\pi\mapsto \sigma_{\pi}$ vient de  la compatibilité local-global (grâce au \lemref{lem-comme-poincare} ci-dessous).  L'application étendue vérifie a) et est compatible avec l'induction parabolique. 
  Enfin la compatibilité aux   cas triviaux de  fonctorialité
  vient de l'énoncé analogue dans le cas global (proposition 12.5 de \cite{coh}). \cqfd

Le \thmref{thm-mfz} va lui-même résulter de la proposition suivante.

   Soit  $s$ un entier et $N$ un niveau.

     Pour tout entier $m$ on note 
 \begin{gather}\label{defi-U-m}
 U_{m}:=\Ker(G(\mc O_{K})\to G(\mc O_{K}/\pi_{K}^{m}\mc O_{K})).
 \end{gather}

\begin{prop}\label{prop-mfz}
Pour toute donnée  
\begin{itemize}
\item d'un ensemble fini $I$, 
\item d'un entier $s$, 
\item d'une fonction $f$ sur le quotient grossier $\wh G \backslash \wh G ^{I}/\wh G$ définie sur $\mc O_{E}/\lambda_{E}^{s}\mc O_{E}$, 
\item d'un $I$-uplet $(\gamma_{i})_{i\in I}$ d'éléments de  $\on{Weil}(\ov K/K)$, 
\item d'un entier $m$ assez grand pour que 
$ U_{m} $
soit d'ordre premier à $\ell$, 
\end{itemize}
il existe un élément $$\mf z_{m,s,I,f,(\gamma_{i})_{i\in I}}\in  
C_{c}(U_{m}\backslash G(K)/U_{m},\mc O_{E}/\lambda_{E}^{s}\mc O_{E})$$ 
dépendant continûment de $(\gamma_{i})_{i\in I}$ 
tel que 
 pour $\Fq,X,  \Xi$ comme ci-dessus, pour toute place $v$ de $X$,  
     en prenant $K=F_{v}$, et pour tout plongement $\ov{ F_X}\subset \ov K$
     (d'où $\on{Weil}(\ov K/K)\subset \on{Gal}(\ov{ F_X}/F_X)$) et pour tout sous-schéma  fini  $N^{v}$ de $X\sm v$,  
         l'opérateur d'excursion $S_{I,f,(\gamma_{i})_{i\in I}}$   de \cite{coh, cong-finite}
     agit sur $  C_{c}(\Bun_{G,N^{v}+mv}(\Fq)/\Xi,
     \mc O_{E}/\lambda_{E}^{s}\mc O_{E})$ par convolution à droite par $\mf z_{m,s,I,f,(\gamma_{i})_{i\in I}}$. 
                \end{prop}

La preuve de cette proposition occupera l'essentiel de l'article 
(du paragraphe \ref{para-cht-res}  au paragraphe \ref{rem-preuve-corresp}). 

\begin{lem}\label{lem-comme-poincare}
Soit $\Fq,X, \Xi$ comme ci-dessus, et  $v$ une place  de $X$. 
   On note  $K=F_{v}$. Soit $m$ un entier et $U_{m}$ le sous-groupe ouvert compact de $G(K)$ défini dans \eqref{defi-U-m}. 
Soit $h\in C_{c}(U_{m}\backslash G(K)/U_{m}, \mc O_{E}/\lambda_{E}^{s}\mc O_{E})$  non nul. Alors il existe un niveau $N^{v}$ en dehors de $v$ et une extension finie $E'$  de $ E$ tels que  l'action de $h$ sur
$  C_{c}(\Bun_{G,N^{v}+mv}(\Fq)/\Xi,
      \mc O_{E'}/\lambda_{E}^{s}\mc O_{E'})$  soit non nulle. 
\end{lem}
\dem On fixe une autre place $v'$ (distincte de $v$). 
On note $G(F_{v'})^{0}$ le sous-groupe distingué de $G(F_{v'})$ engendré par ses sous-groupes ouverts compacts (c'est aussi le noyau du morphisme universel de $G(F_{v'})$ vers un $\Z$-module libre de type fini). 
On fixe 
  une représentation $\pi_{v'}$ supercuspidale de $G(F_{v'})^{0}$, définie sur $E'$ extension finie de $ E$. 
On fixe $m'$ tel que cette représentation ait  au moins un  vecteur non nul invariant par $U'_{m'}$ (le sous-groupe ouvert compact de $G(F_{v'})$ défini comme dans  \eqref{defi-U-m}). On note $h'\in C_{c}(U'_{m'}\backslash G(F_{v'})^{0}/U'_{m'},E')$ un coefficient de matrice non nul de cette représentation, qui est à support compact. On normalise $h'$ pour qu'il prenne ses valeurs dans $\mc O_{E'}$ mais 
que sa réduction modulo l'uniformisante $\lambda_{E'}$ de $\mc O_{E'}$ soit non nulle. Enfin on choisit un niveau $N^{v,v'}$ en dehors de $v$ et $v'$ tel que, en notant $K^{v,v'}\subset \prod_{w\neq v,v'}G(F_{w})$ 
le sous-groupe ouvert compact correspondant à $N^{v,v'}$, les translatés de 
$K^{v,v'}\on{supp}(h) \on{supp}(h')$ par $G(F)\Xi$ soient deux à deux disjoints. 
La série de Poincaré correspondante $g\mapsto \sum_{\gamma\in G(F)\Xi}
(\mathds 1_{K^{v,v'}}hh')(\gamma g)$ est alors non nulle  
dans $C_{c}(\Bun_{G,N^{v,v'}+mv+m'v'}(\Fq)/\Xi,
       \mc O_{E'}/\lambda_{E}^{s}\mc O_{E'})$.  Or c'est l'image par multiplication par $h$ de la série de Poincaré
     $g\mapsto \sum_{\gamma\in G(F)\Xi}
(\mathds 1_{K^{v,v'}}\mathds 1_{U_{m}}h')(\gamma g )$
dans $C_{c}(\Bun_{G,N^{v,v'}+mv+m'v'}(\Fq)/\Xi,  \mc O_{E'}/\lambda_{E}^{s}\mc O_{E'})$. 
\cqfd

\noindent{\bf Démonstration du  \thmref{thm-mfz}  (sauf l'assertion (vi)) en admettant la \propref{prop-mfz}.} D'abord $\mf z_{m,s,I,f,(\gamma_{i})_{i\in I}}$ est déterminé de manière unique par la compatibilité avec le global (grâce au \lemref{lem-comme-poincare}
 et au fait que le niveau $N^{v}$ en dehors de $v$ est arbitraire). On en déduit  qu'il est central (puisque les opérateurs d'excursion globaux commutent avec les opérateurs de Hecke en $v$), compatible avec la limite sur $s$ et $m$. 
 On note $\mf z_{I,f,(\gamma_{i})_{i\in I}}\in  \wh{\mf Z_{\mc O_{E}}}$ l'élément obtenu par limite sur $s$ et $m$. On en déduit que $\mf z_{I,f,(\gamma_{i})_{i\in I}}$  
  vérifie toutes les propriétés du  \thmref{thm-mfz}, sauf la propriété (vi) 
(car les mêmes propriétés sont vérifiées par les opérateurs d'excursion globaux, d'après \cite{coh}). \cqfd

\begin{rem}\label{rem-supports1} 
En fait  le support de $\mf z_{m,s,I,f,(\gamma_{i})_{i\in I}}$  est born\'e ind\'ependamment de $s$ d'apr\`es la \remref{rem-support-para5}. 
En revanche 
 le support de $\mf z_{m,s,I,f,(\gamma_{i})_{i\in I}}$ augmente avec $m$
(puisqu'à la limite on obtient un élément du centre de Bernstein).  \end{rem}

\begin{rem} La \propref{prop-mfz} généralise le lemme 6.11 de 
\cite{coh}. Plus précisément le lemme 6.11 correspond au cas 
particulier où en notant $V$ une représentation de $\wh G$,  $I=\{1,2\}$, 
$f(x_{1},x_{2})=\on{Tr}_{V}(x_{1}x_{2}^{-1})$, $(\gamma_{1},\gamma_{2})=(\Frob_{v},1)$, et enfin $U_{m}=G(\mc O_{K})$ 
(en oubliant la condition qu'il doit être d'ordre premier à $\ell$, que l'on a imposée ici  par commodité plus que par nécessité).   En fait dans ce cas particulier la proposition  6.2  de \cite{coh} montre que $\mf z_{ I,f,(\gamma_{i})_{i\in I}}$
 est égal en niveau sphérique à l'élément $h_{V,v}$ de l'algèbre de Hecke sphérique associé par l'isomorphisme de Satake à $V$.  
Bien évidemment on ne possède pas en général une telle formule explicite 
pour $\mf z_{ I,f,(\gamma_{i})_{i\in I}}$.  En niveau Iwahori on devrait pouvoir faire le lien avec les fonctions traces de Frobenius sur les faisceaux pervers 
construits par Gaitsgory dans \cite{ga-iwahori}. Cependant cela ne fournirait pas de résultat supplémentaire dans le cas déployé car le centre en niveau Iwahori est isomorphe au centre en niveau sphérique. 
\end{rem}

Le reste de ce paragraphe est consacré à la preuve de  l'assertion (vi)
 du  \thmref{thm-mfz}. 

En général soient 
$N,I,W,x,\xi,\delta$ comme avant \eqref{excursion-etendus}. 
Alors en notant $\mc C_{P}^{G, \nu}:
C_{c}(\Bun_{G,N}(\Fq)/\Xi, \mc O_{E})\to C_{c}(\Bun_{M,N}^{',\nu}(\Fq)/\Xi,\mc O_{E})$ 
le morphisme terme constant (comme dans \cite{cong}) on a un diagramme commutatif
  \begin{gather*}\xymatrixcolsep{5pc}
   \xymatrix{
 C_{c}(\Bun_{G,N}(\Fq)/\Xi,\mc O_{E}/\lambda_{E}^{s}\mc O_{E})
    \ar[r]^{S^{G}_{I,f,(\gamma_{i})_{i\in I}}}   \ar[d]^{\prod  \mc C_{P}^{G, \nu}}& C_{c}(\Bun_{G,N}(\Fq)/\Xi,\mc O_{E}/\lambda_{E}^{s}\mc O_{E})\ar[d]^{\prod \mc C_{P}^{G, \nu}}
   \\   
  \prod    C_{c}(\Bun_{M,N}^{',\nu}(\Fq)/\Xi,\mc O_{E}/\lambda_{E}^{s}\mc O_{E})  \ar[r]^{S^{M}_{I,f,(\gamma_{i})_{i\in I}}} &  
\prod  C_{c}(\Bun_{M,N}^{',\nu}(\Fq)/\Xi,\mc O_{E}/\lambda_{E}^{s}\mc O_{E})      } \end{gather*}
où $S^{G}_{I,f,(\gamma_{i})_{i\in I}}$ et $S^{M}_{I,f,(\gamma_{i})_{i\in I}}$ sont les opérateurs  introduits ci-dessus (les exposants $G$ et $M$ indiquent que le premier est pour $G$ et le second pour $M$). 
Dans \cite{cong}, Cong Xue a défini les morphismes terme constant sur la cohomologie à support compact des champs de chtoucas, la commutativité du diagramme ci-dessus résulte alors de leur compatibilité avec les morphismes de création  et d'annihilation. 

\noindent{\bf Démonstration de l'assertion (vi) du   \thmref{thm-mfz}  en admettant la \propref{prop-mfz}.}
On applique la commutativité du diagramme ci-dessus  avec un $I$-uplet $(\gamma_{i})_{i\in I}$ d'éléments de  $\on{Weil}(\ov K/K)$, $N=N^{v}+mv$. 
Alors les opérateurs $S^{G}_{I,f,(\gamma_{i})_{i\in I}}$ et $S^{M}_{I,f,(\gamma_{i})_{i\in I}}$ sont  égaux à la convolution à droite par des éléments 
 $\mf z^{G}_{m,s,I,f,(\gamma_{i})_{i\in I}}$ et  $\mf z^{M}_{m,s,I,f,(\gamma_{i})_{i\in I}}$ des centres de Bernstein pour $G$ et $M$. On en déduit le diagramme commutatif
 
  \begin{gather*} \xymatrixcolsep{5pc}
   \xymatrix{
 C_{c}(\Bun_{G,N}(\Fq)/\Xi,\mc O_{E}/\lambda_{E}^{s}\mc O_{E})
    \ar[r]^{\mf z^{G}_{I,f,(\gamma_{i})_{i\in I}}}   \ar[d]^{\prod  \mc C_{P}^{G, \nu}}& C_{c}(\Bun_{G,N}(\Fq)/\Xi,\mc O_{E}/\lambda_{E}^{s}\mc O_{E})\ar[d]^{\prod \mc C_{P}^{G, \nu}}
   \\   
  \prod    C_{c}(\Bun_{M,N}^{',\nu}(\Fq)/\Xi,\mc O_{E}/\lambda_{E}^{s}\mc O_{E})  \ar[r]^{\mf z^{M}_{I,f^{M},(\gamma_{i})_{i\in I}}} &  
\prod  C_{c}(\Bun_{M,N}^{',\nu}(\Fq)/\Xi,\mc O_{E}/\lambda_{E}^{s}\mc O_{E})      } \end{gather*}

Quitte à diminuer le Levi $M$,   
il suffit de montrer (vi) pour $\tau$ supercuspidale. On tord $\tau$ par un caractère entier générique de $M(\mb A)$ trivial sur $M(F_{X})$ et sur $\Xi$. 
On choisit $m$ tel que $(\on{Ind}_{P(K)}^{G(K)}\tau)^{U_{m}}\neq 0$. 
On choisit $N^{v}$ comme dans  le \lemref{lem-comme-poincare} appliqué à $M$ 
(et à un sous-groupe discret cocompact de $Z_{M}(F_{X})\backslash Z_{M}(\mb A)$ contenant $\Xi$), de sorte que  
 la contragrédiente $\check \tau$ de  $\tau$ apparaît comme facteur en $v$ d'une représentation cuspidale pour $M$ en niveau $N$, et la série d'Eisenstein associée à cette dernière représentation (dont le facteur en $v$ est $(\on{Ind}_{P(K)}^{G(K)}\check \tau)^{U_{m}}$) 
est non nulle, donc a  un produit scalaire non nul avec 
$C_{c}(\Bun_{G,N}(\Fq)/\Xi,E)$. 
D'après la proposition 2.11 de \cite{bernstein-deligne}, 
 $\mf z^{G}_{I,f,(\gamma_{i})_{i\in I}}$ agit sur  $\on{Ind}_{P(K)}^{G(K)}\tau$ par un scalaire, que l'on note  $\lambda$. Comme 
  $\on{Ind}_{P(K)}^{G(K)}\check \tau$ est la représentation contragrédiente de $\on{Ind}_{P(K)}^{G(K)}\tau$, ${}^{t}\mf z^{G}_{I,f,(\gamma_{i})_{i\in I}}$ agit sur 
    $\on{Ind}_{P(K)}^{G(K)}\check \tau$ par $\lambda$. 
 Par l'adjonction entre séries d'Eisenstein et morphismes terme constant, 
 et la commutativité du diagramme ci-dessus, 
 $\lambda$ est aussi le scalaire par lequel ${}^{t}\mf z^{M}_{I,f^{M},(\gamma_{i})_{i\in I}}$ agit sur $\check \tau$. C'est donc aussi le scalaire par lequel 
  $\mf z^{M}_{I,f^{M},(\gamma_{i})_{i\in I}}$ agit sur $\tau$ et on a fini la preuve de (vi). 
\cqfd

\section{Chtoucas restreints}\label{para-cht-res}

On se donne un corps fini $\mf k$ de cardinal premier à $\ell$,  $Y$ une courbe
lisse (non nécessairement projective) sur $\mf k$, et $v$ un point de $Y$ défini sur $\mf k$.  

\begin{rem}
Dans la comparaison entre local et global on aura un morphisme $Y\to X$ 
induisant un isomorphisme de $v$ sur son image. Autrement dit 
le rôle de $Y$ sera celui  d'une localisation hensélienne de $X$ au voisinage de $v$, avec la condition supplémentaire que $\mf k=k_{v}$ et que $Y$ est définie sur $\mf k$ (pour localiser $X$ au sens hensélien en $v$ il suffit de considérer de telles courbes $Y$ et cette condition supplémentaire simplifiera beaucoup la définition des chtoucas restreints sur $Y$). De plus  $\mc O_{K}$ est l'anneau des fonctions sur le    complété de $Y$ en $v$. 
On doit donc considérer $Y$ comme  une algébrisation de $\mc O_{K}$ qui permettra d'appliquer la théorie des cycles proches de \cite{orgogozo}. 
\end{rem}

\begin{notation}
Dans toute la suite on considérera des morphismes de Frobenius relatifs à $Y$ sur $\mf k$, ou bien relatifs à une courbe $X$ sur $\Fq$ pour le lien avec le global. Pour les distinguer, on convient que 
\begin{itemize}
\item
pour tout schéma $S$ sur $\Fq$, {\it resp. } $\mf k$,  
\item et pour tout objet sur 
$X\times_{  \Fq } S$, {\it resp. } $Y\times_{ \mf k } S$ (comme un $G$-torseur, un diviseur), 
\end{itemize}
 on notera par un exposant $\tau $, {\it resp. } $\sigma$ à gauche l'image inverse de cet objet par $\Id_{X}\times \Frob _{S}$, {\it resp. }  $\Id_{Y}\times \Frob _{S}$, où $ \Frob _{S}$ est relatif à $\Fq$, {\it resp. } $\mf k$. 
\end{notation}

Pour l'instant on ne considère que la courbe $Y$ sur $\mf k$ donc les morphismes de Frobenius sont  relatifs à $\mf k$ et on utilise seulement la lettre $\sigma$. 

Si $y:S\to Y$ est un morphisme sur $\mf k$ on considère aussi $y$ comme un diviseur (de Cartier) effectif de $Y\times S$ relativement à $S$. En revanche on  note  $\Gamma_{ y}$  le graphe de $y$ considéré en tant que sous-schéma fermé de $Y\times S$ (car la notation simplifiée $y$ pour ce graphe aurait créé des confusions ensuite). 
 
Pour tout diviseur $D$ de $Y\times S$ (effectif et relatif sur $S$),  on note $\tav D=(\Id_{Y}\times \Frob_{S}) ^{*}(D)$.  Alors   $\tav y$ désigne $y\circ \Frob_{S}$. En effet 
$ \tav y $ est le diviseur de $Y\times S$ image inverse de 
$y$ par $\Id_{Y}\times \Frob_{S}$, qui est le diviseur associé à $y\circ \Frob_{S}$. 

Comme  $S$ n'est pas supposé parfait, $\tuv{-1} y$ n'a pas de sens en tant que morphisme de $S$ vers $Y$. Pour tout entier $a\in \N$ on {\it définit }
$ \tuv{-a} y $ comme le diviseur de $Y\times S$ (effectif et relatif sur $S$)  image inverse de  $y$   par $(\Frob_{Y}\times \Id_{S})^{a}$. Plus généralement, 
pour tout diviseur $D$ de $Y\times S$ (effectif et relatif sur $S$),  on  définit 
$ \tuv{-a} D $ comme le diviseur de $Y\times S$ (effectif et relatif sur $S$) image inverse de  $D$   par $(\Frob_{Y}\times \Id_{S})^{a}$. 
On remarque que $\tav (\tuv{-1} D)=\tuv{-1} (\tav D)=(\sharp \mf k)D$ en tant que diviseur de $Y\times S$. 
On remarque aussi que, en notant $\deg$ le degré relatif, 
 $\deg( \tuv{-a} D )=(\sharp \mf k)^{a}\deg(D)$. 

On en déduit par exemple, que si on se donne un   ensemble fini $J$ , des morphismes  $y_{j}: S\to Y$, des entiers $(a_{j})_{j\in J}\in \N^{J}$ et   $(n_{j})_{j\in J}\in \N^{J}$, 
le diviseur   $\sum_{j\in J} n_{j} \tuv{-a_{j}} y_j$  de $Y\times S$ est effectif et relatif  sur $S$ de degré relatif $\sum  n_{j}(\sharp \mf k)^{a_{j}} $. 
   
   Pour tout   diviseur $Q$ de $Y\times S$ effectif et relatif  sur $S$  de degré relatif $d$, 
    on note 
  $\Gamma_{Q}$ le sous-schéma fermé de 
  $Y\times S$ qui lui est associé, et 
    $G_{Q}$ le schéma en groupes lisse sur $S$ (de dimension $d\dim G$)  obtenu par restriction à la Weil de $G$ de $\Gamma_{Q}$ à $S$. 
    Comme  $\Gamma_{Q}$ est localement libre sur $S$ (de rang $d$) cette restriction \`a la Weil est bien d\'efinie. 
    On a le lemme suivant. 
    
    \begin{lem}\label{lemme-restr-weil}
   Le foncteur de restriction \`a la Weil  associ\'e \`a la projection $\Gamma_{Q}\to S$ donne  une équivalence entre 
    \begin{itemize}
    \item $G$-torseurs sur $\Gamma_{Q}$
     \item $G_{Q}$-torseurs sur $S$. 
    \end{itemize}
    \end{lem}
   Dans la suite on identifie ces deux notions de fa\c con implicite. En particulier si on a une inclusion de diviseurs comme ci-dessus $Q'\subset Q$, 
   et un $G_{Q}$-torseur $\mc G$,   on notera $\restr{\mc G}{Q'}$ le  $G_{Q'}$-torseur qui est le quotient de $\mc G$ \'egal \`a $\mc G \times ^{G_{Q}} G_{Q'}$. 

\dem  Comme $G$, resp.    $G_{Q}$ sont lisses sur $\Gamma_{Q}$, resp. sur $S$, un $G$-torseur, resp. un $G_{Q}$-torseur est localement trivial pour la topologie lisse sur $\Gamma_{Q}$, resp. sur $S$. Le lemme résulte du fait qu'un $G$-torseur sur $\Gamma_{Q}$ est toujours  localement trivial relativement à la topologie lisse sur $S$. En effet, s'il est trivialisé par $Z\to \Gamma_{Q}$ morphisme lisse surjectif, il l'est aussi par la restriction à la Weil $\wt Z$  de $Z$ par le morphisme fini de  $\Gamma_{Q}$ à $S$ (en tirant en arrière la trivialisation sur $Z$ par le morphisme tautologique $\wt Z\times _{Y}\Gamma_{Q}\to Z$). Or il résulte de 
la proposition A.5.2 (4)  page 506 de 
\cite{pseudo-red} ou de la proposition 6 pages 195-196 de 
\cite{neron}
que le morphisme $\wt Z\to Y$ est lisse. 
\cqfd

On   rappelle maintenant  la définition  de la grassmannienne affine de Beilinson-Drinfeld, en reprenant les notations de  \cite{coh}.  

\begin{defi} Soit $I$ un ensemble fini et $(I_{1},...,I_{k})$ une partition ordonnée de $I$.  La grassmannienne  affine  de Beilinson-Drinfeld est l'ind-schéma  $\mr{Gr}_{Y,I }^{(I_{1},...,I_{k})}$ sur $Y^{I}$ dont les $S$-points, pour tout 
$\mf k$-schéma $S$,  classifient la donnée de  
\begin{gather}\label{formule-rem-grassm-intro2}
\big((y_{i})_{i\in I}, \mc G_{0} \xrightarrow{\phi_{1}}  
\mc G_{1}\xrightarrow{\phi_{2}}
\cdots\xrightarrow{\phi_{k-1}}  \mc G_{k-1} \xrightarrow{ \phi_{k}}   \mc G_{k}   \isor{\theta} G_{Y\times S} \big)  \end{gather}
où 
\begin{itemize}
\item les $y_{i}:S\to Y$ sont des morphismes, 
\item les $\mc G_{i}$ sont des $G$-torseurs sur $Y\times S$, 
\item pour tout $j$,  $\phi_{j}$ est un isomorphisme sur $(Y\times S)\sm(\bigcup_{i\in I_{j}}\Gamma_{y_i})$,    
\item $\theta$ est une trivialisation de $\mc G_{k}$. 
\end{itemize}

 D'après Beauville-Laszlo \cite{BL} (voir aussi le premier paragraphe   de \cite{coh} pour des références complémentaires  dans \cite{hitchin}), 
   $\mr{Gr}_{Y,I }^{(I_{1},...,I_{k})}$ peut aussi être défini comme l'ind-schéma dont les $S$-points classifient 
       \begin{gather}\label{formule-rem-grassm-intro-loc}\big((y_{i})_{i\in I}, \mc G_{0} \xrightarrow{\phi_{1}}  
\mc G_{1}\xrightarrow{\phi_{2}}
\cdots\xrightarrow{\phi_{k-1}}  \mc G_{k-1} \xrightarrow{ \phi_{k}}   \mc G_{k}\isor{\theta} G_{\Gamma_{\sum \infty y_i}} \big)  \end{gather}
 où les $\mc G_{i}$ sont des $G$-torseurs  sur le voisinage formel 
  $\Gamma_{\sum \infty y_i}$ de la réunion des graphes des $y_{i}$ dans
 $  Y\times S$, $\phi_{i}$ est un isomorphisme sur $\Gamma_{\sum \infty y_i}\sm(\bigcup_{i\in I_{j}}\Gamma_{y_i})$ (au sens expliqué dans  la notation 1.7 de \cite{coh}) et $\theta$ est une trivialisation de $\mc G_{k}$. 

 La restriction à la Weil $G_{\sum \infty y_i}$ de $G$ de $\Gamma_{\sum \infty y_i}$ à $S$ agit donc sur  $\mr{Gr}_{Y,I }^{(I_{1},...,I_{k})}$    par changement de la trivialisation $\theta$.

Soit $\underline \omega=(\omega_{i})_{i\in I}$ une famille de  poids dominants de $\wh G$. 
 La strate  fermée   $\mr{Gr}_{Y,I,\underline \omega}^{(I_{1},...,I_{k})}$ est le  sous-schéma fermé réduit de 
  $\mr{Gr}_{Y,I }^{(I_{1},...,I_{k})}$ dont les points géométriques vérifient  la condition que 
  la position relative de  $\mc G_{j-1}$ par rapport à  $\mc G_{j}$ en  $y_{i}$ (pour  $i\in I_{j}$) est  bornée par le  copoids dominant de  $G$ correspondant au  poids dominant $\omega_{i}$ de $\wh G$. Plus précisément 
  au-dessus de l'ouvert $\mc U$ de $Y^{I}$ où les $y_{i}$ sont deux à deux distincts, $\mr{Gr}_{Y,I }^{(I_{1},...,I_{k})}$ est un produit (au-dessus de  $\mc U$) de grassmanniennes affines usuelles et 
  \begin{itemize}
  \item on définit 
  la restriction de $\mr{Gr}_{Y,I,\underline \omega}^{(I_{1},...,I_{k})}$ au-dessus de $\mc U$  comme le produit des strates fermées réduites habituelles (définies comme les adhérences de Zariski des $G(\mc O)$-orbites et notées $\ov {\mr{Gr}_{ \omega_{i}}}$ dans  \cite{mv,brav-gaitsgory}), 
  \item puis 
    on   définit  $\mr{Gr}_{Y,I,\underline \omega}^{(I_{1},...,I_{k})}$ comme l'adhérence de Zariski (dans $\mr{Gr}_{Y,I }^{(I_{1},...,I_{k})}$) de 
  sa restriction au-dessus de $\mc U$. 
  \end{itemize}
    \end{defi}
    
  \begin{notation}  \label{notation-satake}
  Soit $\Lambda$ un anneau de coefficients, égal à $E$, $\mc O_{E}$ ou 
    $\mc O_{E}/\lambda_{E}^{s}\mc O_{E}$ (le troisième cas étant celui qui nous concerne   dans cet article). 
    On rappelle que   l'équivalence de Satake géométrique (montrée dans \cite{mv}, avec des compléments dans \cite{ga-de-jong,baumann-riche}, et rappelée dans le théorème 1.17 de \cite{coh}) fournit 
un faisceau pervers \footnote{On considère ici la t-structure perverse usuelle, comme dans \cite{mv} et  \cite{ga-de-jong}. Elle n'est pas stable par dualité de Verdier lorsque $\Lambda$ n'est pas un corps mais cela ne nous pose pas de problème. On renvoie à la note en bas de la page 160 de \cite{ga-de-jong} pour un commentaire au sujet de cette t-structure perverse.}
 $G_{\sum \infty y_{i} }$-équivariant  $\mc S_{I,W}^{\Lambda, (I_{1},...,I_{k})}$ sur 
$\mr{Gr}_{Y,I }^{(I_{1},...,I_{k})}$ pour toute représentation $\Lambda$-linéaire de type fini 
$W$ de $\wh G^{I}$. On note $\mr{Gr}_{Y,I,W }^{(I_{1},...,I_{k})}$ le support 
de $\mc S_{I,W}^{\Lambda, (I_{1},...,I_{k})}$. C'est une réunion finie de strates fermées 
$\mr{Gr}_{Y,I,\underline \omega}^{(I_{1},...,I_{k})}$ comme ci-dessus. 
\end{notation}

\begin{rem} \label{rem-W-reg} D'après l'appendice B de \cite{ga-de-jong}, et comme on l'a rappelé dans le premier paragraphe de \cite{coh} le faisceau pervers $\mc S_{I,W}^{\Lambda, (I_{1},...,I_{k})}$  est obtenu de la manière suivante. 
  Soit $\mc R$ l'algèbre des fonctions régulières sur $\wh G$ (à coefficients dans $\Lambda$). Munie de l'action (\`a droite) de   $\wh G$  par multiplication à gauche c'est une limite inductive de représentations de $\wh G$ sur des $\Lambda$-modules libres de type fini. Alors 
$\boxtimes_{i\in I} \mc R$  est l'algèbre des fonctions régulières sur $\wh G^{I}$. 
 La construction de  $\mc S_{I,W}^{ \Lambda,(I_{1},...,I_{k})}$ est réalisée  par la formule 
  $$\mc S_{I,W}^{ \Lambda,(I_{1},...,I_{k})}=
\Big( \mc S_{I,\boxtimes _{i\in I}\mc R }^{\Lambda,(I_{1},...,I_{k})} \otimes W\Big)^{\wh G^{I}}
,$$
où les invariants (d\'eriv\'es) sont pris par l'action diagonale de $\wh G^{I}$, 
et où l'action de $\wh G^{I}$ sur $\mc S_{I,\boxtimes _{i\in I}\mc R }^{\Lambda,(I_{1},...,I_{k})} $ provient de l'action de $\wh G$ par multiplication à droite sur 
$\mc R$.  
Par conséquent pour toutes les propriétés de ces faisceaux et des faisceaux associés sur les champs de chtoucas restreints, et notamment pour l'action des morphismes de Frobenius partiels que l'on construira plus tard, il suffit de considérer des 
  représentations $W$ de $\wh G ^{I}$ 
de la forme $\boxtimes_{i\in I} W_{i}$, où 
chaque $W_{i}$ est une 
 représentation de $\wh G$ sur un $\Lambda$-module libre de type fini 
  (en effet $\boxtimes _{i\in I}\mc R$ est une limite inductive de telles représentations). 
\end{rem}

  \begin{rem} Lorsque $W$ est une représentation $E$-linéaire de dimension finie  
  de $\wh G^{I}$,  $\mr{Gr}_{Y,I,W}^{(I_{1},...,I_{k})}$ est la réunion des $\mr{Gr}_{Y,I,\underline \omega}^{(I_{1},...,I_{k})}$ pour $\underline \omega$ poids dominant d'un constituant irréductible de $W$. Lorsque $W$ est une représentation 
  $\mc O_{E}/\lambda_{E}^{s}\mc O_{E}$-linéaire de type fini, $\mr{Gr}_{Y,I,W}^{(I_{1},...,I_{k})}$ est sans doute plus compliqué à décrire mais nous n'avons pas besoin de le connaître explicitement. 
    \end{rem}
  
  L'action de  $G_{\sum \infty y_i}$   sur  $\mr{Gr}_{Y,I }^{(I_{1},...,I_{k})}$  préserve  $\mr{Gr}_{Y,I,W}^{(I_{1},...,I_{k})}$. 
De plus  si les entiers $n_{i}$ sont assez grands en fonction de $W$, l'action de $G_{\sum \infty y_i}$ sur $\mr{Gr}_{Y,I,W}^{(I_{1},...,I_{k})}$ se factorise par 
  $G_{\sum n_{i} y_i}$.

 On aura besoin du résultat un peu plus fort suivant. 
 Pour toute famille de  morphismes $y_{i}:S\to Y$ on note $\mr{Gr}_{Y,I,W}^{(I_{1},...,I_{k}), (y_{i})_{i\in I}}=\mr{Gr}_{Y,I,W}^{(I_{1},...,I_{k})}\times_{Y^{I}}S$. 
  Pour tout diviseur $Q$ de $Y\times S$ effectif et relatif  sur $S$, on note 
   $\mr{Gr}_{Y,I,W}^{(I_{1},...,I_{k}), (y_{i})_{i\in I}, \mr{triv}(Q)}$ le $G_{Q}$-torseur  sur 
    $\mr{Gr}_{Y,I,W}^{(I_{1},...,I_{k}), (y_{i})_{i\in I}}$ dont la fibre en un point \eqref{formule-rem-grassm-intro2} est la restriction \`a la Weil  de $\restr{\mc G_{0}}{Q}$ (comme dans le \lemref{lemme-restr-weil}). De plus il est muni d'une action de $G_{Q+\sum \infty y_i}$. 

En effet on a la variante suivante de Beauville-Laszlo:    $\mr{Gr}_{Y,I }^{(I_{1},...,I_{k}), (y_{i})_{i\in I}}$   classifie les        \begin{gather}\label{formule-rem-grassm-intro-loc}\big( \mc G_{0} \xrightarrow{\phi_{1}}  
\mc G_{1}\xrightarrow{\phi_{2}}
\cdots\xrightarrow{\phi_{k-1}}  \mc G_{k-1} \xrightarrow{ \phi_{k}}   \mc G_{k}\isor{\theta} G_{Q+\Gamma_{\sum \infty y_i}} \big)  \end{gather}
 où les $\mc G_{i}$ sont des $G$-torseurs  sur  
  $\Gamma_{Q+\sum \infty y_i}$, $\phi_{i}$ est un isomorphisme sur $\Gamma_{Q+\sum \infty y_i}\sm(\bigcup_{i\in I_{j}}\Gamma_{y_i})$ et $\theta$ est une trivialisation de $\mc G_{k}$. De plus $\mr{Gr}_{Y,I,W}^{(I_{1},...,I_{k}), (y_{i})_{i\in I}, \mr{triv}(Q)}$  est le $G_{Q}$-torseur  sur 
    $\mr{Gr}_{Y,I,W}^{(I_{1},...,I_{k}), (y_{i})_{i\in I}}$ dont la fibre est la restriction \`a la Weil  de $\restr{\mc G_{0}}{Q}$. 
 Donc $G_{Q+\sum \infty y_i}$  agit  sur  $\mr{Gr}_{Y,I,W}^{(I_{1},...,I_{k}), (y_{i})_{i\in I}, \mr{triv}(Q)}$  par changement de la trivialisation $\theta$.

  \begin{lem}\label{lem-n_i-W}
    Si les entiers $n_{i}$ sont assez grands en fonction de $W$, alors pour tout diviseur $Q$ de $Y\times S$ effectif et relatif  sur $S$, 
    l'action de $G_{Q+\sum \infty y_i}$ sur l'espace total du $G_{Q}$-torseur 
     $\mr{Gr}_{Y,I,W}^{(I_{1},...,I_{k}), (y_{i})_{i\in I}, \mr{triv}(Q)}$  
     se factorise par $G_{Q+\sum n_{i} y_i}$. \cqfd
        \end{lem}

\begin{construction}
Pour tout $G_{Q+\sum n_{i} y_i}$-torseur $\mc G$ sur   $S$, et pour toute section $z$ de 
\begin{gather}\label{lem-n_i-W1}\mr{Gr}_{Y,I,W}^{(I_{1},...,I_{k}), (y_{i})_{i\in I}}\times^{G_{ \sum n_{i} y_i}}\restr{\mc G}{\sum n_{i} y_i}\end{gather}   sur $S$ on note 
$a_{z}(\mc G)$ le $G_{Q}$-torseur sur $S$ dont l'espace total est la fibre de 
\begin{gather}\label{lem-n_i-W2}\mr{Gr}_{Y,I,W}^{(I_{1},...,I_{k}), (y_{i})_{i\in I}, \mr{triv}(Q)}\times^{G_{Q+\sum n_{i} y_i}}\mc G\end{gather} 
au-dessus de $z$. 
En effet \eqref{lem-n_i-W2} est un $G_{Q}$-torseur sur \eqref{lem-n_i-W1} car \eqref{lem-n_i-W1} est \'egal \`a 
\begin{gather}\nonumber \mr{Gr}_{Y,I,W}^{(I_{1},...,I_{k}), (y_{i})_{i\in I}}\times^{G_{Q+\sum n_{i} y_i}}\mc G\end{gather} 

En particulier soit $   \mc G_{k} $ un $G$-torseur sur $Y\times S$ et 
\begin{gather}\nonumber
\big( (y_{i})_{i\in I}, ( \mc G_{0} \xrightarrow{\phi_{1}}  
\mc G_{1}\xrightarrow{\phi_{2}}
\cdots\xrightarrow{\phi_{k-1}}  \mc G_{k-1} \xrightarrow{ \phi_{k}}  \mc G_{k})  \big) \end{gather}  comme dans \eqref{formule-rem-grassm-intro2} (sans $\theta$) tel qu'en restreignant les $ \mc G_{i}$ et $\phi_{i} $ à 
$\Gamma_{\sum \infty y_i}$ on obtienne un 
$S$-point $z$ de 
$\mr{Gr}_{Y,I,W }^{(I_{1},...,I_{k}), (y_{i})_{i\in I}}\times^{G_{\sum \infty y_{i}}}\restr{ \mc G_{k}}{\sum \infty y_{i}} $. 
Alors 
en posant $\mc G=\restr{\mc G_{k}}{Q+\sum n_{i}y_{i}}$, 
$a_{z}(\mc G)$ correspond par le \lemref{lemme-restr-weil} à $\restr{ \mc G_{0}}{Q}$. 
En d'autres termes $a_{z}(\mc G)$ est la modification de $\mc G$ par $z$ et  le lemme dit que cette modification est possible avec des niveaux finis si on accepte une perte de niveau de $\sum_{i\in I}n_{i}y_{i}$. 
\end{construction}

\begin{notation}\label{notation-niW}
Dans toute la suite on notera $(n_{i}^{W})_{i\in I}$ des entiers tels que  le lemme précédent soit vrai pour  tout
$I$-uplet $(n_{i})_{i\in I}$ vérifiant 
\begin{gather}\label{ni-niW}n_{i}\geq n_{i}^{W}\text{ \  pour tout \ }i\in I.
\end{gather}
\end{notation}

Soit $I$ un ensemble fini et $(I_{1},...,I_{k})$ une partition ordonnée de $I$.  
Soit  $W$ une  représentation $\Lambda$-linéaire de type fini   
  de $\wh G^{I}$. Soit  $(n_{i})_{i\in I}\in \N^{I}$   vérifiant \eqref{ni-niW}. Soit $n\in \N$ (qui indiquera le niveau en $v$). 
  
  Dans la définition suivante, l'entier $r\in \N$ indique une ``réserve'' de données au-dessus des images inverses des pattes par les puissances de $\Frob_{Y}$,    qui permettra de définir les morphismes de Frobenius partiels sur les champs de chtoucas restreints. 
  
 \begin{defi}\label{defi-Chr} 
 On définit     $\Chr_{Y,I,W,(n_{i})_{i\in I},r} ^{nv,(I_{1},...,I_{k})}$ comme le champ d'Artin    sur $Y^{I}$  dont les points sur 
 un schéma   $S$ sur  $\mf k$  classifient  les données    \begin{gather}\label{donnee-chtouca-res}\big( (y_i)_{i\in I}, \mc G , z ,   \theta   
\big)
\end{gather}
avec 
 \begin{itemize}
\item $y_i:S\to Y  $ un $\mf k$-morphisme, pour $i\in I$, 
\item $\mc G $ est un 
$G_{nv+\sum n_{i}y_{i} + \sum n_{i}\tuv{-1}y_{i} + \cdots+ \sum n_{i}\tuv{-r}y_{i}}$-torseur sur $S$, c'est-à-dire un 
$G$-torseur sur 
$\Gamma_{nv+\sum n_{i}y_{i} + \sum n_{i}\tuv{-1}y_{i} + \cdots+ \sum n_{i}\tuv{-r}y_{i}}$, 
\item  $z$ est 
un $S$-point de 
 $\mr{Gr}_{Y,I,W}^{(I_{1},...,I_{k}), (y_{i})_{i\in I}}\times^{G_{\sum n_{i}y_{i} }} \restr{\mc G }{\sum n_{i}y_{i} }$ 
\item en notant  $a_{z}(\mc G)$   le 
$G_{nv + \sum n_{i}\tuv{-1}y_{i} + \cdots+ \sum n_{i}\tuv{-r}y_{i}}$-torseur sur $S$ 
obtenu en appliquant la construction précédente à $Q=nv + \sum n_{i}\tuv{-1}y_{i} + \cdots+ \sum n_{i}\tuv{-r}y_{i}$, 
$$\theta:\restr{\tav (a_{z}(\mc G)) }{nv + \sum n_{i} y_{i} + \cdots+ \sum n_{i}\tuv{-(r-1)}y_{i}}\isom \restr{\mc G }{nv + \sum n_{i} y_{i} + \cdots+ \sum n_{i}\tuv{-(r-1)}y_{i}}$$ est un isomorphisme de $G_{nv + \sum n_{i} y_{i} + \cdots+ \sum n_{i}\tuv{-(r-1)}y_{i}}$-torseurs sur $S$. 
 \end{itemize}
 \end{defi}

\begin{rem}
Comme  
$\tav (\tuv{-j} y_{i})=(\sharp \mf k)(\tuv{-(j-1)}y_{i})$ pour tout $j\geq 1$, 
$\tav (a_{z}(\mc G))$ est défini en fait sur $\Gamma_{nv+(\sharp \mf k)( \sum n_{i} y_{i} + \cdots+ \sum n_{i}\tuv{-(r-1)}y_{i} )}$ et l'isomorphisme $\theta$ ci-dessus concerne sa restriction à $\Gamma_{nv+  \sum n_{i} y_{i} + \cdots+ \sum n_{i}\tuv{-(r-1)}y_{i}  }$. 
\end{rem}

\begin{rem}
Deux  variantes de la définition précédente seraient possibles pour définir les champs de chtoucas restreints: 
\begin{itemize}
\item [] 1) pour éviter les puissances négatives de $\sigma$, relever 
les $y_{i}$ par $\Frob_{Y}^{r}$ en des points   $z_{i}$, et demander que $\mc G$ soit un 
$G_{nv+\sum n_{i}z_{i} + \sum n_{i}\tav z_{i} + \cdots+ \sum n_{i}\tuv{r}z_{i}}$-torseur sur $S$, que la modification ait lieu en $\sum n_{i}\tuv{r}z_{i}$, et que l'isomorphisme $\theta$ soit défini sur $nv + \sum n_{i} \tav z_{i}  + \cdots+ \sum n_{i}\tuv{r}z_{i} $, 
\item [] 2) garder les points $y_{i}$ mais demander que $\mc G$ soit un 
\newline $G_{nv+q^{r}\sum n_{i}y_{i} + q^{r-1}\sum n_{i}\tuv{-1}y_{i} + \cdots+ \sum n_{i}\tuv{-r}y_{i}}$-torseur, que la modification ait lieu en $q^{r} \sum n_{i}y_{i}$, et que l'isomorphisme   $\theta$ soit défini sur 
$\Gamma_{nv+  q^{r}\sum n_{i} y_{i} + \cdots+ q\sum n_{i}\tuv{-(r-1)}y_{i}  }$. 
\end{itemize}
Les variantes 1) et 2) ci-dessus sont en fait très semblables (on passe de 2) à 1) en relevant les $y_{i}$ par les $z_{i}$ et en multipliant tous les $n_{i}$ par $q^{r}$). 
La   définition ci-dessus est peut-être moins naturelle  que ces deux variantes  mais elle présente l'avantage de raccourcir les formules.
\end{rem}

Pour faciliter la compréhension de la 
  définition précédente et des constructions qui vont suivre, 
  on introduit dès à présent  les morphismes de restriction 
depuis les champs de chtoucas globaux. Le lecteur peut supposer d'abord $\deg(v)=1$ car le cas où $\deg(v)>1$ est plus compliqué. 

On introduit donc ici le cadre global: un corps fini $\Fq$ inclus dans $\mf k$, une courbe projective lisse géométriquement irréductible $X$ sur $\Fq$, 
 un groupe réductif $G$ sur $X$ (supposé déployé jusqu'au dernier paragraphe),   un point fermé $v$ de 
$X$ tel que $\mf k$ est égal au corps résiduel $k_{v}$, et enfin 
un niveau $N^{v}$  sur $X$ en dehors de $v$. 
 On se donne de plus un morphisme étale de schémas sur $\Fq$ 
 \begin{gather}\label{mor-Y-X}
 \pi:Y\to X\sm N^{v}\end{gather} 
 tel que 
 \begin{gather}\label{cond-X-Y-v}
 \text{l'image inverse du point fermé } v \text{ de } X \text{ consiste en l'{\it unique} 
  point ferm\'e } v \text{ de }  Y\\ \nonumber \text{ et le morphisme sur les corps résiduels est l'égalité }
  k_{v}=\mf k \text{ fixée ci-dessus.}\end{gather}    
 
 On rappelle que $Y$ est une courbe lisse sur $\mf k$. 
 
 On convient que $Y^{I}$ désignera toujours le produit sur $\mf k$ (et non sur $\Fq$). 
On note $U_{I,r}$ l'ouvert de $Y^{I}$ formé des $(y_{i})_{i\in I}$ tels que  
les $\pi(y_{i})$ vérifient la condition suivante: 
\begin{gather} \label{cond-UIr} \forall i,i'\in I, \forall j\in \{1,...,r\deg(v)\}  \text{ non multiple de }\deg(v) ,  \tu{j} (\pi(y_{i}))\neq \pi(y_{i'}). \end{gather} 
 Cet ouvert $U_{I,r}$ contient $\Delta(v)$ où $\Delta:Y\to Y^{I}$ désigne le morphisme diagonal et cela nous suffira   car c'est en ce point que nous prendrons les cycles proches. 

\begin{rem}\label{rem-deux-div-X-Y}
Si $S$ est un $\mf k$-schéma et $(y_{i})_{i\in I}:S\to Y^{I}$ un morphisme,  alors 
  on note encore $\pi$ par abus le morphisme 
  $$Y\times_{\mf k}S
\to (X\sm N^{v})\times _{\Fq} S,$$ et on remarque que $\pi$ fournit un isomorphisme de    schémas finis sur $S$ de 
   \begin{gather}\label{div-yi-mfk}\Gamma_{nv+\sum n_{i}y_{i} + \sum n_{i}\tuv{-1}y_{i} + \cdots+ \sum n_{i}\tuv{-r}y_{i}} \subset Y\times_{\mf k}S,  \end{gather}
 vers 
  \begin{gather}\label{div-yi-Fq}\Gamma_{n(\pi \circ v)+\sum n_{i}(\pi\circ y_{i} )+ \sum n_{i}\tu {-\deg(v)}(\pi\circ y_{i}) + \cdots+ \sum n_{i}\tu{-\deg(v)r}(\pi\circ y_{i})} \subset (X\sm N^{v})\times_{\Fq}S.  \end{gather}
  Dans \eqref{div-yi-mfk} et \eqref{div-yi-Fq},  $v$ désigne le morphisme 
  \begin{gather}
  \label{v-mor-S-Y}
  S\to \on{Spec}(\mf k)\xrightarrow{v} Y.\end{gather} 
\end{rem}

\begin{defi}
Le morphisme de restriction   
\begin{gather}\label{morp-glob-restr}
 \mc R_{ I,W,(n_{i})_{i\in I},r} ^{nv,(I_{1},...,I_{k})}: 
\Cht_{N^{v}, I,W} ^{(I_{1},...,I_{k})}\times _{(X\sm N^{v})^{I}}U_{I,r}\to 
\Chr_{Y,I,W,(n_{i})_{i\in I},r} ^{nv,(I_{1},...,I_{k})}\end{gather}
associe à un $S$-point 
\begin{gather}\nonumber 
\Big(\big( (x_i)_{i\in I},    ( \mc G_{0}  \xrightarrow{\phi_{1}}   \mc G_{1}  \xrightarrow{\phi_{2}}
\cdots\xrightarrow{\phi_{k-1}}   \mc G_{k-1}  \xrightarrow{ \phi_{k}}  \mc G_{k}\isor{\lambda} \ta  \mc G_{0})       
\big), \psi, (y_i)_{i\in I}\Big)
\end{gather}
le chtouca restreint $\big( (y_i)_{i\in I}, \mc G , z ,   \theta   
\big)
$ avec 
\begin{itemize}
\item 
$\mc G$ égal à la restriction de $ \pi^{*}\big(\tu{\deg(v)} \mc G_{0}\big)=\tav \big(\pi^{*} \mc G_{0}\big)  $ à \eqref{div-yi-mfk}, 
\item comme construction auxiliaire, 
$$\alpha: \restr{ \pi^{*} \mc G_{k}}{nv+\sum n_{i}y_{i} + \sum n_{i}\tuv{-1}y_{i} + \cdots+ \sum n_{i}\tuv{-r}y_{i}}\isom \mc G$$ l'isomorphisme obtenu par restriction de $$\tu{\deg(v)-1} \big( \lambda\phi_{k}...\phi_{1}\big)... \ta\big( \lambda\phi_{k}...\phi_{1}\big)\lambda$$ (c'est un isomorphisme
 grâce aux conditions \eqref{cond-X-Y-v} et  \eqref{cond-UIr}), 
\item 
$z$ associé à la restriction de $$(  \pi^{*}\mc G_{0}  \xrightarrow{\phi_{1}}   \pi^{*} \mc G_{1}  \xrightarrow{\phi_{2}}
\cdots\xrightarrow{\phi_{k-1}}   \pi^{*} \mc G_{k-1}  \xrightarrow{ \phi_{k}}  \pi^{*} \mc G_{k})$$ à \eqref{div-yi-mfk}, plus précisément $z$ est un $S$-point de
$$\mr{Gr}_{Y,I,W}^{(I_{1},...,I_{k}), (y_{i})_{i\in I}}\times^{G_{\sum n_{i}y_{i}}} \restr{\pi^{*}\mc G_{k}}{\sum n_{i}y_{i}}\isom \mr{Gr}_{Y,I,W}^{(I_{1},...,I_{k}), (y_{i})_{i\in I}}\times^{G_{\sum n_{i}y_{i}}} \restr{\mc G }{\sum n_{i}y_{i}}
$$
où l'isomorphisme vient de $\alpha$ ci-dessus, si bien  que \begin{gather}\label{oz-de-g}
a_{z}(\mc G)=\restr{ \pi^{*} \mc G_{0} }{nv + \sum n_{i}\tuv{-1}y_{i} + \cdots+ \sum n_{i}\tuv{-r}y_{i}}, \end{gather}
\item et enfin $$\theta:\restr{\tav(a_{z}(\mc G))}{nv+  \sum n_{i} y_{i} + \cdots+ \sum n_{i}\tuv{-(r-1)}y_{i} }\isom \restr{ \mc G }{nv+  \sum n_{i} y_{i} + \cdots+ \sum n_{i}\tuv{-(r-1)}y_{i} } $$ résultant 
de \eqref{oz-de-g} et 
de la définition de $\mc G$ comme restriction de $\tav \big(\pi^{*} \mc G_{0}\big)  $ à \eqref{div-yi-mfk}.  
\end{itemize}
\end{defi}

Lorsque $I$ est vide, 
$\Chr_{Y,\emptyset,\mbf 1,()_{i\in \emptyset},r} ^{nv,(\emptyset)}
=\bullet/G(\mc O_{nv})$. Plus généralement  on a 
\begin{gather}\label{mor-GON}
\restr{\Chr_{Y,I,W,(n_{i})_{i\in I},r}^{nv,(I_{1},...,I_{k})}}{(Y\sm v)^{I}}=
\restr{\Chr_{Y,I,W,(n_{i})_{i\in I},r}^{\emptyset,(I_{1},...,I_{k})}}{(Y\sm v)^{I}} 
\times \bullet/G(\mc O_{nv}). 
\end{gather} 
On remarque d'ailleurs que $\Cht_{N^{v}+nv,I,W}^{(I_{1},...,I_{k})}\times _{(X\sm N^{v}\sm v)^{I}}(U_{I,r}\cap (Y\sm v)^{I})$,  qui 
est un $G(\mc O_{nv})$-torseur sur $ \Cht_{N^{v},I,W}^{(I_{1},...,I_{k})} \times _{(X\sm N^{v} )^{I}}(U_{I,r}\cap (Y\sm v)^{I})$, 
s'obtient par image inverse du $G(\mc O_{nv})$-torseur 
$ \bullet\to \bullet/G(\mc O_{nv})$, ou encore du $G(\mc O_{nv})$-torseur 
\begin{gather}\label{GOn-torsur-ChtR} \restr{\Chr_{Y,I,W,(n_{i})_{i\in I},r}^{\emptyset,(I_{1},...,I_{k})}}{(Y\sm v)^{I}}
\to 
\restr{\Chr_{Y,I,W,(n_{i})_{i\in I},r}^{nv,(I_{1},...,I_{k})}}{(Y\sm v)^{I}}. 
\end{gather}
C'est pourquoi on a mis l'indice $nv$   en haut pour les champs de chtoucas restreints alors que l'indice $N$  est en bas pour les champs de chtoucas globaux.

\noindent {\bf Notation.} 
Pour tout anneau de coefficients $\Lambda$ comme dans la notation \ref{notation-satake}, on  note $\mc L_{I,W,(n_{i})_{i\in I},r}^{nv, \Lambda, (I_{1},...,I_{k})}$ le $\Lambda$-faisceau lisse  sur 
$\restr{\Chr_{Y,I,W,(n_{i})_{i\in I},r}^{nv,(I_{1},...,I_{k})}}{(Y\sm v)^{I}}$
égal à l'image directe de $\Lambda$ par le morphisme fini donnant le  $G(\mc O_{nv})$-torseur \eqref{GOn-torsur-ChtR}. Autrement dit c'est le faisceau associé à la représentation régulière de $G(\mc O_{nv})$ à coefficients dans $\Lambda$.

On remarque que lorsque $r=0$ et $n=0$, 
 $$\Chr_{Y,I,W,(n_{i})_{i\in I},0} ^{\emptyset,(I_{1},...,I_{k})}=
\mr{Gr}_{Y,I,W}^{(I_{1},...,I_{k})}/G_{\sum n_{i}y_{i} }.$$
En particulier la proposition suivante généralise la proposition 2.8 de \cite{coh}. 

\begin{prop}\label{prop-lisse1} 
Le morphisme de restriction   \eqref{morp-glob-restr}
 est lisse de dimension $q^{r\deg(v)}(\sum n_{i})\dim G$ sur $Y^{I}$. 
\end{prop}
\dem 
Comme le morphisme $ \Cht_{N^{v}, I,W} ^{(I_{1},...,I_{k})}\to  \restr{\Cht_{  I,W} ^{(I_{1},...,I_{k})}}{(X\sm N^{v})^{I}}$ est étale il suffit de montrer l'énoncé avec $N^{v}$ vide, ce qu'on suppose maintenant pour alléger les notations  
  (mais on pourrait
au contraire prendre $N^{v}$ assez grand en fonction d'une troncature de Harder-Narasimhan pour avoir des schémas lisses au lieu de champs lisses dans ce qui suit). On fixe un schéma $S$ sur $U_{I,r}$ et un $S$-point \eqref{donnee-chtouca-res} de $\Chr_{Y,I,W,(n_{i})_{i\in I},r} ^{nv,(I_{1},...,I_{k})}$. 
Pour raccourcir les formules on pose $\mc G'=a_{z}(\mc G)$. 
On veut montrer que la fibre de \eqref{morp-glob-restr} au-dessus de \eqref{donnee-chtouca-res} est lisse sur $S$. 
 Elle est l'égalisateur de deux morphismes $a,b $ 
(de champs sur $S$) : 
\begin{gather*} \Bun_{G,nv + \sum n_{i} (\pi \circ y_{i}) +\sum n_{i}\tu{-1}(\pi \circ y_{i}) + \cdots+ \sum n_{i}\tu{-r\deg(v)}(\pi \circ y_{i})}
\times \big(\mc G\times \tu{r-1}\mc G' \times ...\times  \ta\mc G'\big) \to  \\ 
\Bun_{G,nv + \sum n_{i} (\pi \circ y_{i}) +\sum n_{i}\tu{-1}(\pi \circ y_{i})+ \cdots+ \sum n_{i}\tu{-(r\deg(v)-1)}(\pi \circ y_{i})}
\times \big(\mc G\times \tu{r-1}\mc G' \times ...\times  \ta\mc G'\big) . \end{gather*}
Avant de définir $a,b$ on doit expliquer les notations. D'abord $v$ apparaît ici comme place de $X$ et comme $S$ est un schéma sur $\mf k$, le niveau $nv$ équivaut en fait à un niveau $n(\pi\circ v)+ n\ta (\pi\circ v)+ ... + n\tu{(r-1)} (\pi\circ v)$ en considérant $v$ comme le  morphisme \eqref{v-mor-S-Y} 
(on   utilise   pour $\pi\circ v$ des puissances positives de $\tau$ plutôt que des puissances négatives qui épaissiraient les points).  
Grâce aux conditions  \eqref{cond-X-Y-v} et  \eqref{cond-UIr}, les diviseurs (de 
$(X\sm N^{v})\times_{\Fq}S$ )  
\begin{gather}\label{div-0}  n  (\pi\circ v)+ \sum n_{i} (\pi \circ y_{i}) + ...+ \sum n_{i}\tu{- r  \deg(v)}(\pi \circ y_{i}),   \text{ et } 
\\  \label{div-j}
  n\tu{ r -j } (\pi\circ v)+ \sum n_{i}\tu{-j}(\pi \circ y_{i}) + ...+ \sum n_{i}\tu{-(r-1)\deg(v)-j}(\pi \circ y_{i})  \end{gather}
pour $j\in \{1,...,r-1\}$,  sont deux à deux disjoints, et permettent de tordre respectivement par 
$\mc G$ et $\tu{r-j}(\mc G') $ pour $j\in \{1,...,r-1\}$. 
On définit   $a$ comme 
  le morphisme d'oubli, qui est lisse de dimension 
$q^{r\deg(v)} (\sum n_{i}) \dim G$. 
On définit $b $  
comme  la composée 
\begin{gather*}
\Bun_{G,nv+\sum n_{i} (\pi \circ y_{i})+ \sum n_{i}\tu{-1}(\pi \circ y_{i}) + \cdots+ \sum n_{i}\tu{-r\deg(v)}(\pi \circ y_{i})}
\times \big(\mc G\times \tu{r-1}\mc G' \times ...\times  \ta\mc G'\big)\to 
\\ 
\Bun_{G,nv+\sum n_{i}\tu{-1}(\pi \circ y_{i}) + \cdots+ \sum n_{i}\tu{-r\deg(v)}(\pi \circ y_{i})} 
\times \big(\mc G' \times \tu{r-1}\mc G' \times ...\times  \ta\mc G'\big)\xrightarrow{\Frob_{/ X^{I}}}
    \\ \Bun_{G,nv + q\sum n_{i} (\pi \circ y_{i}) + \cdots+ q\sum n_{i}\tu{-(r\deg(v)-1)}(\pi \circ y_{i})}   
\times \big( \tu{r}\mc G' \times \tu{r-1}\mc G' \times ...\times  \ta\mc G'\big)
 \to \\
\Bun_{G,nv +  \sum n_{i} (\pi \circ y_{i}) + \cdots+  \sum n_{i}\tu{-(r\deg(v)-1)}(\pi \circ y_{i})} \times \big( \tu{r}\mc G' \times \tu{r-1}\mc G' \times ...\times  \ta\mc G'\big) \to \\
\Bun_{G,nv +  \sum n_{i} (\pi \circ y_{i}) + \cdots+  \sum n_{i}\tu{-(r\deg(v)-1)}(\pi \circ y_{i})}\times  \big(  \mc G  \times \tu{r-1}\mc G' \times ...\times  \ta\mc G'\big)
.
\end{gather*}
où 
\begin{itemize}
\item le premier morphisme est la modification par $z$, 
qui a lieu uniquement sur le diviseur \eqref{div-0} (et change donc la torsion par $\mc G$ en une torsion par $\mc G'$ sur ce diviseur, avec une  perte de niveau $\sum n_{i}(\pi \circ y_{i})$), 
\item le deuxième provient 
  du  morphisme de Frobenius  du champ $\Bun_{G,nv+\sum n_{i}\tu{-1}x_{i} + \cdots+ \sum n_{i}\tu{-r\deg(v)}x_{i}}
$ relatif à $X^{I}$ (plus une permutation circulaire de $\tu{r}\mc G'$, ..., $\ta\mc G'$), donc la différentielle de ce morphisme    est nulle dans les fibres au-dessus de $S$,  
  \item  le troisième est
  l'oubli (qui enlève les facteurs $q$ inutiles dans les diviseurs), 
   \item  le quatrième   change, grâce à  l'isomorphisme $\theta$ du chtouca restreint,  la torsion par 
   $\tu{r}\mc G' $ en une torsion par $\mc G$ 
    qui, grâce à la \remref{rem-deux-div-X-Y}, est défini  sur  \eqref{div-yi-Fq}
    et que l'on restreint à  $\Gamma_{n(\pi \circ v)+\sum n_{i}(\pi\circ y_{i} )+ \sum n_{i}\tu {-\deg(v)}(\pi\circ y_{i}) + \cdots+ \sum n_{i}\tu{-\deg(v)(r-1)}(\pi\circ y_{i})} $.   \end{itemize}
  Donc $b$  a aussi la propriété que sa différentielle est nulle dans les fibres au-dessus de $S$.  
  On conclut par la lissité de l'égalisateur entre deux morphismes $a,b$ de champs lisses sur $S$, tels que $a$  est lisse et $b$ a une différentielle nulle le long des fibres au-dessus de $S$. 
   \cqfd

On se place de nouveau dans un cadre purement local. 

\begin{lem} \label{lem-oubli-chtR-Gr} 
On suppose $n\geq n' $, $n_{i}\geq n'_{i}$ pour tout $i\in I$ et $r\geq r'$. Alors
 le morphisme évident 
\begin{gather}\label{oubli-chtR-ChtR}\Chr_{Y,I,W,(n_{i})_{i\in I},r} ^{nv,(I_{1},...,I_{k})} 
\to \Chr_{Y,I,W,(n'_{i})_{i\in I},r'} ^{n' v,(I_{1},...,I_{k})} 
\end{gather}
qui à \eqref{donnee-chtouca-res}
associe  \begin{gather}\label{donnee-chtouca-res-'}\big( (y_i)_{i\in I}, \restr{\mc G }{n' v+\sum n'_{i}y_{i} + \sum n'_{i}\tuv{-1}y_{i} + \cdots+ \sum n'_{i}\tuv{-r'}y_{i}}, z ,  
  \restr{\theta}{n' v + \sum n'_{i} y_{i} + \cdots+ \sum n'_{i}\tuv{-(r'-1)}y_{i}}   
\big)
\end{gather}
  est lisse (de dimension $\big( q^{r\deg(v)}(\sum n_{i})-q^{r'\deg(v)}(\sum n'_{i})\big) \dim G$). \end{lem}
    \dem
 Ce morphisme est la composée (de droite à gauche) 
 \begin{itemize}
 \item du morphisme lisse représentable qui consiste à restreindre $\theta$
   $$\text{de \  }\Gamma_{nv + \sum n_{i} y_{i} + \cdots+ \sum n_{i}\tuv{-(r-1)}y_{i}}\text{ \  à \ 
 } \Gamma_{n' v + \sum n'_{i} y_{i} + \cdots+ \sum n'_{i}\tuv{-(r'-1)}y_{i}},$$  
  et dont les fibres sont isomorphes au schéma lisse  
  $$\ker \big(G_{nv + \sum n_{i} y_{i} + \cdots+ \sum n_{i}\tuv{-(r-1)}y_{i}}\to G_{n' v + \sum n'_{i} y_{i} + \cdots+ \sum n'_{i}\tuv{-(r'-1)}y_{i}}\big)$$ tordu à gauche et à droite par $\tav (a_{z}(\mc G)) $ et $\mc G $,   
 \item du morphisme lisse  consistant à restreindre $\mc G $ 
  $$\text{de \  }\Gamma_{nv+\sum n_{i}y_{i} + \sum n_{i}\tuv{-1}y_{i} + \cdots+ \sum n_{i}\tuv{-r}y_{i}} \text{ \  à \ 
 } \Gamma_{n' v+\sum n'_{i}y_{i} + \sum n'_{i}\tuv{-1}y_{i} + \cdots+ \sum n'_{i}\tuv{-r'}y_{i}},$$ et dont les fibres sont des gerbes sur le schéma en groupes lisse
   $$\ker(G_{nv+\sum n_{i}y_{i} + \sum n_{i}\tuv{-1}y_{i} + \cdots+ \sum n_{i}\tuv{-r}y_{i}}  \to G_{n' v+\sum n'_{i}y_{i} + \sum n'_{i}\tuv{-1}y_{i} + \cdots+ \sum n'_{i}\tuv{-r'}y_{i}})$$  tordu par $\mc G$.    \cqfd
 \end{itemize}

En particulier, en prenant $n'=0$ et $r'=0$, on voit que  le morphisme 
\begin{gather}\label{oubli-chtR-Gr}\Chr_{Y,I,W,(n_{i})_{i\in I},r} ^{nv,(I_{1},...,I_{k})} 
\to  \mr{Gr}_{Y,I,W}^{(I_{1},...,I_{k})}/G_{\sum n_{i}y_{i} }\end{gather}
qui à \eqref{donnee-chtouca-res}
associe \begin{gather}
\big( (y_i)_{i\in I}, \restr{\mc G}{\sum n_{i}y_{i}},  z  \big)
\end{gather}
  est lisse.

\noindent {\bf Notation.} 
On définit $\mc S_{ I,W,(n_{i})_{i\in I},r} ^{nv,\Lambda,(I_{1},...,I_{k})}$ comme l'image inverse de $\mc S_{I,W}^{\Lambda, (I_{1},...,I_{k})}$ 
(introduit dans la notation \ref{notation-satake}) par le morphisme lisse  \eqref{oubli-chtR-Gr}. C'est donc un faisceau pervers à un décalage près sur $\Chr_{Y,I,W,(n_{i})_{i\in I},r} ^{nv,(I_{1},...,I_{k})} $. 

\begin{rem}\label{prod-GrI(I)}
On suppose  $W$ de la forme $\boxtimes_{i\in I} W_{i}$, où 
chaque $W_{i}$ est une 
 représentation de $\wh G$ sur un $\Lambda$-module libre de type fini. 
On déduit du lemme 1.16 et du théorème 1.17 de  \cite{coh} (rappelant l'équivalence de Satake géométrique) que   le morphisme évident 
\begin{gather}\label{oubli-chtR-GrI(I)}\Chr_{Y,I,W,(n_{i})_{i\in I},r} ^{nv,(I_{1},...,I_{k})} 
\to  \prod_{j=1}^{k} \mr{Gr}_{Y,I_{j},\boxtimes_{i\in I_{j}} W_{i}}^{(I_{j})}/G_{\sum_{i\in I_{j}} n_{i}y_{i} }\end{gather}
est lisse et que    $\mc S_{ I,W,(n_{i})_{i\in I},r} ^{nv,\Lambda,(I_{1},...,I_{k})}$ est égal à l'image inverse de 
$\boxtimes \mc S_{I_{j},\boxtimes_{i\in I_{j}} W_{i}}^{\Lambda, (I_{j})}$ par le morphisme \eqref{oubli-chtR-GrI(I)}.  
\end{rem}

On rappelle que 
$\mc S_{ I,W,(n_{i})_{i\in I},r} ^{nv,\Lambda, (I_{1},...,I_{k})}$
est un faisceau pervers   sur $\Chr_{Y,I,W,(n_{i})_{i\in I},r} ^{nv,(I_{1},...,I_{k})}$, qui est ULA sur $(Y\sm v)^{I}$ (en fait même sur $Y^{I}$ dans le cas non déployé considéré ici) et que 
$\mc L_{I,W,(n_{i})_{i\in I},r}^{nv, \Lambda, (I_{1},...,I_{k})}$ (d\'efini avant la \propref{prop-lisse1}) est un faisceau lisse (en $\Lambda$-modules localement libres) sur $\restr{\Chr_{Y,I,W,(n_{i})_{i\in I},r} ^{nv,(I_{1},...,I_{k})}}{(Y\sm v)^{I}}$.  

Pour abréger tout en gardant les indices nécessaires on note 
$$ \mc S\mc L_{ I,W,(n_{i})_{i\in I},r} ^{nv, \Lambda,(I_{1},...,I_{k})}=
\mc S_{ I,W,(n_{i})_{i\in I},r} ^{nv,\Lambda,(I_{1},...,I_{k})}\otimes_{\Lambda} \mc L_{I,W,(n_{i})_{i\in I},r}^{nv, \Lambda, (I_{1},...,I_{k})}$$
qui est donc un faisceau pervers  (ULA sur $(Y\sm v)^{I} $) sur $\restr{\Chr_{Y,I,W,(n_{i})_{i\in I},r} ^{nv,(I_{1},...,I_{k})}}{(Y\sm v)^{I}}$. 

\noindent {\bf Notation.} Soit $(m_{i})_{i\in I}\in \N^{I}$. On introduit le morphisme 
$$\Delta_{(m_{i})_{i\in I}} :Y\to Y^{I},  \ \ y\mapsto (\Frob_{Y}^{m_{i}}(y))_{i\in I}. $$
En particulier $\Delta_{(0)_{i\in I}}:Y\to Y^{I}$ est le morphisme diagonal, et 
en général $\Delta_{(m_{i})_{i\in I}}$ est son image par le morphisme de Frobenius partiel $\prod_{i\in I} \Frob_{Y,\{i\}}^{m_{i}}$.

Donc \begin{gather}\label{deltabig} \Delta_{(m_{i})_{i\in I}}^{*}\Big(\mc S\mc L_{ I,W,(n_{i})_{i\in I},r} ^{nv, \Lambda,(I_{1},...,I_{k})}
\Big)\end{gather}
est un faisceau pervers sur 
$ \Chr_{Y,I,W,(n_{i})_{i\in I},r} ^{nv,(I_{1},...,I_{k})}\times _{Y^{I},\Delta_{(m_{i})_{i\in I}}} (Y\sm v) $ où la notation indique que le produit fibré est pris pour le morphisme $\Delta_{(m_{i})_{i\in I}}: Y\sm v\to Y^{I}$. 

Dans toute la suite on fixe un point géométrique algébrique 
$\ov v$ au-dessus de $v$. 

Le formalisme des cycles proches sur un trait   au sens de SGA7 (et \cite{laumon-moret-bailly} 18.4 pour les champs d'Artin), appliquée à \eqref{deltabig},  fournit donc 
$$R\Psi_{(m_{i})_{i\in I}} \mc S\mc L_{ I,W,(n_{i})_{i\in I},r} ^{nv, \Lambda,(I_{1},...,I_{k})}\in D_{c}^{b}(\restr{\Chr_{Y,I,W,(n_{i})_{i\in I},r} ^{nv,(I_{1},...,I_{k})}}{\Delta(\ov v)},\Lambda).$$ 

\begin{rem} 
Le
  \lemref{lem-oubli-chtR-Gr} montre que 
  pour tout  
  $(n_{i})_{i\in I}\in \N^{I}$ vérifiant \eqref{ni-niW},   pour tout $r\geq 0$,  le morphisme d'oubli  
\begin{gather}\label{mor-oubli-r-0}\Chr_{Y,I,W,(n_{i})_{i\in I},r}^{nv,(I_{1},...,I_{k})}\to \Chr_{Y,I,W,(n_{i}^{W})_{i\in I},0}^{nv,(I_{1},...,I_{k})}\end{gather} est lisse. Donc 
 on voit que pour tout $(m_{i})_{i\in I}\in \N^{I}$, 
$R\Psi_{(m_{i})_{i\in I}} \mc S\mc L_{ I,W,(n_{i})_{i\in I},r} ^{nv, \Lambda,(I_{1},...,I_{k})}$ est l'image inverse de $R\Psi_{(m_{i})_{i\in I}} \mc S\mc L_{ I,W,(n_{i}^{W})_{i\in I},0} ^{nv, \Lambda,(I_{1},...,I_{k})}$ par la 
 restriction du morphisme  \eqref{mor-oubli-r-0} en la fibre au-dessus de $\Delta(\ov v)$. 
\end{rem}
 
 On va introduire  les morphismes de Frobenius partiels, qui diminuent $r$ de $1$. C'est d'ailleurs pour qu'ils soient définis et puissent être itérés que l'on a introduit l'entier $r$. D'après la \remref{rem-W-reg}, on peut supposer $W$ de la forme $\boxtimes_{i\in I} W_{i}$, où 
chaque $W_{i}$ est une 
 représentation de $\wh G$ sur un $\Lambda$-module libre de type fini.    
 On construit   le morphisme 
 $$\on{Fr\mc R}_{Y,I_{1}}:\Chr_{Y,I,W,(n_{i})_{i\in I},r} ^{nv,(I_{1},...,I_{k})}\to 
 \Chr_{Y,I,W,(n_{i})_{i\in I},r-1} ^{nv,(I_{2},...,I_{k},I_{1})}. $$
 L'idée est la suivante: partant d'un $S$-point 
 \eqref{donnee-chtouca-res} de l'espace de départ, 
 on lui associe le $S$-point   $\big( (y'_i)_{i\in I},  \mc G'  , z' ,   \theta '   
\big)$ de l'espace d'arrivée 
où 
 \begin{itemize}
 \item $y'_{i}=\Frob_{Y}(y_{i})$ si $i\in I_{1}$ et $y'_{i}=y_{i}$ sinon, 
 \item on définit le point $z_{2,...,k}$ de     $\mr{Gr}_{Y,I\sm I_{1},\boxtimes_{i\in I\sm I_{1}} W_{i}}^{(I_{2},...,I_{k})}\times_{G_{\sum n_{i}y_{i} }} \restr{\mc G }{\sum n_{i}y_{i} }$ 
 obtenu en ``tronquant'' $z$, c'est-à-dire en oubliant la modification la plus à gauche, en les pattes indexées par $I_{1}$, 
 \item on introduit   le $G_{nv+\sum_{i\in I_{1}} n_{i}y_{i} + \sum_{i\in I }  n_{i}\tuv{-1}y_{i} + \cdots+ \sum_{i\in I }  n_{i}\tuv{-r}y_{i}}$-torseur  $a_{z_{2,...,k}}(\mc G)$ et on pose 
 $$\mc G'=\restr{ \tav(a_{z_{2,...,k}}(\mc G))}{nv+\sum_{i\in I } n_{i}y'_{i} + \sum_{i\in I }  n_{i}\tuv{-1}y'_{i} + \cdots+ \sum_{i\in I }  n_{i}\tuv{-(r-1)}y'_{i}}, $$
 \item on note $z_{1}$ le $S$-point de $\mr{Gr}_{Y,  I_{1},\boxtimes_{i\in I_{1}} W_{i}}^{(I_{1})}\times_{G_{\sum_{i\in I_{1} } n_{i}y_{i} }} \restr{a_{z_{2,...,k}}(\mc G)}{\sum_{i\in I_{1}} n_{i}y_{i} }$ obtenu en ``tronquant'' $z$ de l'autre côté, 
\item on définit $z'$ comme le  $S$-point de $\mr{Gr}_{Y,  I ,W}^{(I_{2},...,I_{k}, I_{1})}\times_{G_{\sum_{i\in I  } n_{i}y'_{i} }} \restr{ \mc G' }{\sum_{i\in I }   n_{i}y'_{i} }$ obtenu en recollant $\tav  z_{1}$ et $z_{2,...,k}$ par $\theta$, 
\item on définit enfin $\theta'$ de fa\c con naturelle. 
 \end{itemize}
  
   Les détails de la construction sont laissés au lecteur car les notations sont très compliquées alors que l'idée est très simple. Un cas particulier est rédigé dans la preuve du lemme 6.11 de \cite{coh}. On peut résumer la construction par le dessin
  $$   \xymatrix{
a_{z}(\mc G) &  a_{z_{2,...,k}}(\mc G)  \ar[l]^{a_{z_1} }&& \mc G \ar[ll]^-{a_{z_{2,...,k}} } \ar@/_2pc/[lll]^{a_{z}} & \ar[l]_{a_{(\tav z_1)} } \mc G'
 \ar@/^2pc/[lll]^{a_{z'}}  }$$
  où l'on voit que l'on a coupé $a_{z}$ en deux morceaux $a_{z_1}$, 
  $a_{z_{2,...,k}} $, que l'on a appliqué  le Frobenius au premier, et qu'on les a recollés dans l'autre sens pour obtenir $a_{z'}$. 
  
   On compare de nouveau avec le cadre global. 
 Les morphismes de restrictions entrelacent 
le morphisme $\on{Fr\mc R}_{Y,I_{1}}$ ci-dessus  avec $\on{Fr}_{X,I_{1}}^{\deg(v)}$, où $\on{Fr}_{X,I_{1}}$ est le morphisme de Frobenius partiel    sur les champs de  chtoucas globaux défini dans \cite{coh} (la condition \eqref{cond-UIr} assure que les pattes dont on doit changer l'ordre pour itérer $\deg(v)$ fois   
$\on{Fr}_{X,I_{1}}$ sont   disjointes). 

Les trois isomorphismes suivants seront sur $\restr{\Chr_{Y,I,W,(n_{i})_{i\in I},r} ^{nv,(I_{1},...,I_{k})}}{(Y\sm v)^{I}}$. 

 On a d'abord un isomorphisme canonique évident 
 $$\on{Fr\mc R}_{Y,I_{1}}^{*}(\mc L_{I,W,(n_{i})_{i\in I},r-1}^{nv, \Lambda, (I_{2},...,I_{k},I_{1})})\isom 
 \mc L_{I,W,(n_{i})_{i\in I},r}^{nv, \Lambda,(I_{1},...,I_{k}) }. $$
 Grâce à la \remref{prod-GrI(I)} on construit   un isomorphisme canonique 
 $$\on{Fr\mc R}_{Y,I_{1}}^{*}(
    \mc S_{ I,W,(n_{i})_{i\in I},r-1} ^{nv,\Lambda,(I_{2},...,I_{k},I_{1})})\isom 
     \mc S_{ I,W,(n_{i})_{i\in I},r} ^{nv,\Lambda,(I_{1},...,I_{k})}.$$
 Le produit tensoriel des deux isomorphismes précédents donne   un isomorphisme canonique 
 $$\on{Fr\mc R}_{Y,I_{1}}^{*}(\mc S\mc L_{ I,W,(n_{i})_{i\in I},r-1} ^{nv, \Lambda,(I_{2},...,I_{k},I_{1})})\isom 
    \mc S\mc L_{ I,W,(n_{i})_{i\in I},r} ^{nv, \Lambda,(I_{1},...,I_{k})}.$$

 On en déduit un isomorphisme canonique 
 \begin{gather}
 \label{corresp-coho-frob-partiels}  
 F\mc R_{I_{1}}: \Big(\restr{\on{Fr\mc R}_{Y,I_{1}}}{\Delta(\ov v)}\Big)^{*}
 \Big(R\Psi_{(m_{i}+\chi_{I_{1}}(i))_{i\in I}} \mc S\mc L_{ I,W,(n_{i})_{i\in I},r-1} ^{nv, \Lambda,(I_{2},...,I_{k},I_{1})} \Big) \\ \nonumber \isom 
 R\Psi_{(m_{i})_{i\in I}} \mc S\mc L_{ I,W,(n_{i})_{i\in I},r} ^{nv, \Lambda,(I_{1},...,I_{k})}\end{gather}
 sur $\restr{\Chr_{Y,I,W,(n_{i})_{i\in I},r} ^{nv,(I_{1},...,I_{k})}}{\Delta(\ov v)}$. 

Comme dans \cite{coh}, la seule raison pour laquelle on   introduit des partitions $(I_{1},...,I_{k})$ arbitraires au lieu de la partition grossière $(I)$ est que cela est   nécessaire pour définir les morphismes de Frobenius partiels. On a montré dans \cite{coh} que les faisceaux de cohomologie des champs de chtoucas globaux ne dépendent pas du choix de la partition. 
  L'énoncé ci-dessous est analogue pour les chtoucas restreints. Il est  absolument nécessaire pour pouvoir itérer les morphismes de Frobenius partiels. 
  
  \begin{lem} \label{lem-petit}
  Le morphisme d'oubli des modifications intermédiaires
  $\pi: \Chr_{Y,I,W,(n_{i})_{i\in I},r} ^{nv,(I_{1},...,I_{k})}\to 
  \Chr_{Y,I,W,(n_{i})_{i\in I},r} ^{nv,(I)}$ est représentable et propre  et il induit un isomorphisme
  $$R\pi_{*}\big( \mc S\mc L_{ I,W,(n_{i})_{i\in I},r} ^{nv, \Lambda,(I_{1},...,I_{k})} \big) \isom  \mc S\mc L_{ I,W,(n_{i})_{i\in I},r} ^{nv, \Lambda,(I)}. $$
  \end{lem}
 
 L'image du  morphisme ci-dessus par les cycles proches est donc 
 \begin{gather}\label{oubli-cycles-proches}R\pi_{*}\Big( R\Psi_{(m_{i})_{i\in I}} \mc S\mc L_{ I,W,(n_{i})_{i\in I},r} ^{nv, \Lambda,(I_{1},...,I_{k})}\Big)  \isom  R\Psi_{(m_{i})_{i\in I}} \mc S\mc L_{ I,W,(n_{i})_{i\in I},r} ^{nv, \Lambda,(I)}. \end{gather}

 Pour la suite on notera que les morphismes 
  \eqref{corresp-coho-frob-partiels} et  \eqref{oubli-cycles-proches}  
 fournissent des correspondances cohomologiques entre faisceaux de cycles proches sur la fibre spéciale de champs de chtoucas restreints. 
 
\section{Une construction utilisant les cycles proches sur une base générale 
(de Deligne, Gabber, Laumon, Orgogozo)} 

Pour montrer la \propref{prop-mfz} on utilisera comme boîte noire la \propref{prop-boite-noire} ci-dessous. L'intérêt de procéder ainsi est que  son énoncé ne fait intervenir que des cycles proches au sens habituel, alors que  sa démonstration repose de fa\c con essentielle sur les cycles proches sur une base générale, et notamment sur les résultats de 
\cite{orgogozo}. On rappelle que ces cycles proches ont été introduits par Deligne, expliqués dans une note écrite par Laumon \cite{laumon-vanishing}
et que les résultats essentiels (comme la constructibilité après modification de la base) ont été obtenus par Orgogozo dans \cite{orgogozo}. On peut également consulter à ce sujet l'article de survey d'Illusie 
 \cite{illusie-vanishing}. On utilise aussi un th\'eor\`eme non publi\'e de Gabber  rappel\'e dans le th\'eor\`eme A.9 de \cite{arnaud-cong}. 

On choisit une bijection $I=\{1,...,l\}$. Autrement dit on choisit un ordre des coordonnées. Seuls des arguments globaux nous permettront   de montrer que le résultat final de la construction ne dépend pas du choix de cet ordre. 

On suppose $X$ et $Y$ comme dans \eqref{mor-Y-X}. 

Pour la commodité des notations (mais sans incidence mathématique) on introduit une uniformisante $t$ au voisinage de $v$ dans $Y$,
et note $\mf k=k_{v}$ le corps résiduel en $v$,  de sorte que le corps local $K=F_{v}$ s'identifie à $  \mf k((t))$. 
On note $K^{1,...,l}= \mf k((t_{1}))...((t_{l}))$ le corps local en dimension supérieure. Donc $\on{Spec}(K^{1,...,l})$ est ``asymptote'' au drapeau 
$v\subset Y\subset ... \subset Y^{l}$ associé au choix de l'ordre des coordonnées. 

On fixe un plongement  \begin{gather}\label{plong-FX-FY-K}
\ov{ F_X}= \ov{ F_Y}\subset \ov K,
\end{gather} d'où des inclusions
$$\on{Weil}(\ov K/K)\subset \on{Weil}(\ov{ F_Y}/F_Y)=\on{Weil}(\eta_{Y}, \ov{   \eta_{Y}})
\subset \on{Weil}(\ov{ F_X}/F_X)=\on{Weil}(\eta_X, \ov{ \eta_X})  $$
en notant $\ov{ \eta_X}=\ov{ \eta_{Y}}$ les  points géométriques de $X$ et $Y$ correspondant à 
$\ov{ F_X}=\ov{ F_Y}$. 
On déduit aussi de \eqref{plong-FX-FY-K} une identification entre $\ov{\Fq}=\ov{\mf k}$ et le corps résiduel de $\ov K$. 

 On note 
   $F_X^{I}$  le corps des fonctions de $X^{I}$, $(F_X^{I})^{\mr{perf}}$ son 
  perfectisé,  $\ov{F_X^{I}}$ une clôture algébrique de $F_X^{I}$ et  
  $\ov{\eta_{X}^{I}}=\on{Spec}(\ov{F_X^{I}})$. 
On rappelle  que dans la remarque 8.7 de \cite{coh}  on avait défini    \begin{gather*}\on{FWeil}(\eta_X^{I},\ov{\eta_X^{I}})=
       \big\{\delta \in \on{Aut}_{\ov\Fq}(\ov{F_X^{I}}), \exists (d_{i})_{i\in I}\in \Z^{I}, \restr{\delta}{(F_X^{I})^{\mr{perf}}}=\prod_{i\in I}(\Frob_{X,\{i\}})^{d_{i}}\big\} 
  .\end{gather*} 
Ce groupe 
         est une  extension de $\Z^{I}$ par   $\on{Ker}(\pi_{1}(\eta_X^{I},\ov{\eta_X^{I}})\to \wh \Z)$. Lorsque  $I$ est un   singleton, il s'identifie au groupe de   Weil   usuel 
   $\on{Weil} (\eta_X,\ov{\eta_X})=\pi_{1}(\eta_X,\ov{\eta_X})\times_{\wh \Z}\Z$
   (mais leurs actions sur $\ov{\eta_X}$ ne sont pas les mêmes: si $\delta$ a pour image $d$ dans $\mathbb Z$, son action sur $\ov F$ en tant qu'élément de $\on{FWeil}(\eta_X^{I},\ov{\eta_X^{I}})$    est par 
$\Frob^{d} \circ \delta$, où $\delta$ est l'action de $\delta$ sur $\ov F$ en tant qu'élément de $\on{Weil} (\eta_X,\ov{\eta_X})$).

         On renvoie au début du chapitre 8 de \cite{coh} pour des rappels sur les flèches de spécialisation. En particulier, pour tout $(d_{i})_{i\in I}\in \Z^{I}$ et tout $(m_{i})_{i\in I}\in \N^{I}$ tel que $(m_{i}+d_{i})_{i\in I}\in \N^{I}$, tout élément $\delta$  de $\on{FWeil}(\eta_X^{I},\ov{\eta_X^{I}})$ d'image $(d_{i})_{i\in I}$ fournit une  flèche de spécialisation   
 \begin{gather}\label{fleche-spec-Frob-eta-eta}  
  \mf{sp}_{\delta}:   \prod_{i\in I}(\Frob_{X,\{i\}})^{m_{i}+d_{i}}\big( \ov{\eta_X^{I}}\big) \xrightarrow{\delta } 
        \prod_{i\in I}(\Frob_{X,\{i\}})^{m_{i} }\big( \ov{\eta_X^{I}}\big)    \end{gather}
(on note que $\delta\mapsto \mf{sp}_{\delta} $ est un anti-morphisme pour la composition).

Comme dans \cite{coh} on choisit une flèche de spécialisation  
  $\mf{sp}: \ov{\eta_X^{I}}\to \Delta(\ov{\eta_X})$. D'après la discussion après 
  la remarque 8.18 dans  \cite{coh}, elle  fournit 
  une inclusion   \begin{gather}\label{incl-F--F-FI-ov}\ov{F_X}\otimes_{\ov\Fq} \cdots 
 \otimes_{\ov\Fq}\ov{F_X}\subset \ov{F_X^{I}} 
  .\end{gather}  Par restriction des automorphismes, on  en déduit   un 
   morphisme  surjectif 
    \begin{gather}\label{morph-Weil-I-intro}\on{FWeil}(\eta_X^{I},\ov{\eta_X^{I}})\to  \on{Weil} (\eta_X,\ov{\eta_X}) ^{I}\end{gather} (dépendant du choix de  $\mf{sp}$).

On fixe un morphisme de  $\ov{\eta_X^{I}}$  vers le point générique $\eta_Y^{I}$ de $Y^{I}$ et on pose  $\ov{\eta_Y^{I}}=\ov{\eta_X^{I}}$ en tant que   point géométrique au-dessus de $\eta_Y^{I}$. 
On 
définit   $\on{FWeil}(\eta_Y^{I},\ov{\eta_Y^{I}})$ de la même fa\c con que ci-dessus pour $X$ 
(en remarquant que  les morphismes de Frobenius partiels sont relatifs à $\mf k=k_{v}$).  
Le choix du plongement \eqref{incl-F--F-FI-ov}   fournit (grâce aux égalités 
$\ov{F_X}=  \ov{F_Y}$ 
 et $\ov{F_X^{I}}=  \ov{F_Y^{I}}$)  
un  morphisme  surjectif 
    \begin{gather}\label{morph-Weil-I-intro-Y}\on{FWeil}(\eta_Y^{I},\ov{\eta_Y^{I}})\to  \on{Weil} (\eta_Y,\ov{\eta_Y}) ^{I}.\end{gather} 

On fixe une clôture algébrique $\ov{K^{1,...,l}}$ de $K^{1,...,l}$ et 
un plongement 
\begin{gather}\label{plong-K-K1..k}\ov K\otimes _{\ov {\mf k}} .... \otimes _{\ov {\mf k}} \ov K \subset \ov{K^{1,...,l}}.  \end{gather}

On définit $\on{FWeil}(\ov{K^{1,...,l}}/K^{1,...,l})$ de fa\c con similaire à 
$\on{FWeil}(\eta_X^{I},\ov{\eta_X^{I}})$. Le plongement \eqref{plong-K-K1..k} 
 fournit un morphisme 
   \begin{gather}\label{morph-Weil-loc}\on{FWeil}(\ov{K^{1,...,l}}/K^{1,...,l})
   \to   (\on{Weil}(\ov K/K))^{I}.  \end{gather}

On fixe en fait un  diagramme   commutatif de plongements  
    \begin{gather} \label{diag-plong} \xymatrix{
 \ov{F_X^{I}}
    \ar@{=}[r]  &  \ov{F_Y^{I}}
    \ar@{^{(}->}[r]  & 
\ov{K^{1,...,l}} 
   \\ \ov{ F_X}\otimes _{\ov {\mf k}} .... \otimes _{\ov {\mf k}} \ov{ F_X} 
      \ar@{^{(}->}[u]   \ar@{=}[r] &  
      \ov{ F_Y}\otimes _{\ov {\mf k}} .... \otimes _{\ov {\mf k}} \ov{ F_Y} 
      \ar@{^{(}->}[u]   \ar@{^{(}->}[r] &  
 \ov K\otimes _{\ov {\mf k}} .... \otimes _{\ov {\mf k}} \ov K 
    \ar@{^{(}->}[u]               } \end{gather}
    où la flèche verticale de gauche est \eqref{incl-F--F-FI-ov}, 
    celle de droite est \eqref{plong-K-K1..k}, et la ligne du bas est la puissance $I$-ième de \eqref{plong-FX-FY-K}. Un tel diagramme existe car 
 il équivaut à prolonger le plongement $    \ov{ F_Y}\otimes _{\ov {\mf k}} .... \otimes _{\ov {\mf k}} \ov{ F_Y}\hookrightarrow  \ov K\otimes _{\ov {\mf k}} .... \otimes _{\ov {\mf k}} \ov K\hookrightarrow \ov{K^{1,...,l}} $ en un plongement 
 $\ov{F_Y^{I}}\hookrightarrow \ov{K^{1,...,l}} $, et cela est toujours possible. 
    
    On déduit de \eqref{diag-plong} 
     un diagramme commutatif de morphismes de groupes de Weil  
   \begin{gather}\label{diag-gpes-Weil} \xymatrix{
 \on{FWeil} (\eta_{X}^{I}, \ov{\eta_{X}^{I}})
     \ar@{->>}[d]  & \on{FWeil} (\eta_{Y}^{I}, \ov{\eta_{Y}^{I}})
     \ar@{->>}[d]    \ar@{_{(}->}[l]_-{\iota_{X,Y}}  & 
\on{FWeil}(\ov{K^{1,...,l}}/K^{1,...,l})  \ar@{_{(}->}[l]_-{\iota_{Y}}    \ar[d] 
   \\ (\on{Weil} (\eta_{X} , \ov{\eta_{X}}))^{I}
     &  (\on{Weil} (\eta_{Y} , \ov{\eta_{Y}}))^{I}  \ar@{_{(}->}[l]
     &  \ 
 (\on{Weil}(\ov K/K))^{I} 
        \ar@{_{(}->}[l]           } \end{gather}

On verra dans le  paragraphe suivant que la flèche verticale de droite (qui était aussi le morphisme \eqref{morph-Weil-loc} ci-dessus) est également surjective. On note 
\begin{gather}\label{iota-X}
\iota_{X}=\iota_{X,Y}\circ \iota_{Y}:\on{FWeil}(\ov{K^{1,...,l}}/K^{1,...,l})\hookrightarrow   \on{FWeil} (\eta_{X}^{I}, \ov{\eta_{X}^{I}}).\end{gather} 
On remarque  que 
$\iota_{X}$ ne dépend pas de $Y$. Il multiplie le degré (à valeurs dans $\Z^{I}$)  par $\deg(v)$. 
  
Soit  
 $m$ un entier tel que 
$ U_{m} $
soit d'ordre premier à $\ell$.   

Dans la proposition suivante on notera simplement $\mc S\mc L$ 
au lieu de $\mc S\mc L_{ I,W,(n_{i})_{i\in I},r} ^{mv, \Lambda,(I_{1},...,I_{k})}$. 
Par ailleurs on notera 
$R\Gamma_{c}(\Cht_{\Delta(\ov v)}^{\leq \mu}/\Xi , R\Psi_{(m_{i})_{i\in I}} \mc S\mc L)$ au lieu de 
\begin{gather}\label{H-Cht-mu}
  R\Gamma_{c}(\restr{\Cht_{N^{v},I,W}^{(I_{1},...,I_{k}),\leq\mu}}{\Delta(\ov v)}/\Xi , 
  R\Psi_{(m_{i})_{i\in I}} \mc S\mc L). \end{gather}
La notation est justifiée par le fait que \eqref{H-Cht-mu} ne dépend pas du choix de la partition $(I_{1},...,I_{k})$. 
De même on notera, pour tout point géométrique $\ov x$ sur $(X\sm N)^{I}$,  
$R\Gamma_{c}(\Cht_{\ov x}^{\leq \mu}/\Xi ,  \mc S\mc L)$ au lieu de 
$$  R\Gamma_{c}(\restr{\Cht_{N^{v},I,W}^{(I_{1},...,I_{k}),\leq\mu}}{\ov x}/\Xi ,   \mc S\mc L). $$

\begin{notation}  \label{assez-croissant}
On rappelle qu'on a choisi une bijection $I=\{1,...,l\}$. 
On dira que $(m_{i})_{i\in I}$ est assez croissant si les $m_{i+1}-m_{i}$ sont assez grands, c'est-à-dire si $m_{i+1}-m_{i}\geq d$ pour $d$ assez grand. 
\end{notation} 

\begin{prop}\label{prop-boite-noire}
Pour tout $\delta \in \on{FWeil}(\ov{K^{1,...,l}}/K^{1,...,l})$
d'image $(d_{i})_{i\in I}\in \Z^{I}$, et pour $(m_{i})_{i\in I}\in \N^{I}$ assez croissant (et tel que $(m_{i}+d_{i})_{i\in I}\in \N^{I}$), on possède  un morphisme $$T_{\delta}: R\Psi_{(m_{i})_{i\in I}} \mc S\mc L\to 
R\Psi_{(m_{i}+d_{i})_{i\in I}} \mc S\mc L$$ dans  
$D_{c}^{b}(\restr{\Chr_{Y,I,W,(n_{i})_{i\in I},0} ^{mv,(I_{1},...,I_{k})} }{\Delta(\ov v)},\Lambda)$
ne dépendant que des données locales, et
tel que pour toutes données globales et  pour tout $\mu$ le diagramme suivant commute
     \begin{gather}\nonumber
 \xymatrix{
  R\Gamma_{c}(\Cht_{\Delta(\ov v)}^{\leq \mu}/\Xi , R\Psi_{(m_{i})_{i\in I}} \mc S\mc L)
    \ar[d]  \ar[r]^{T_{\delta}} & 
 R\Gamma_{c}(\Cht_{\Delta(\ov v)}^{\leq \mu}/\Xi , R\Psi_{(m_{i}+d_{i})_{i\in I}} \mc S\mc L)  \ar[d]
   \\
 R\Gamma_{c}(\Cht^{\leq \mu}_{\prod_{i\in I} \Frob_{X,\{i\}}^{\deg(v)m_{i}} (\Delta(\ov{\eta_X}))}/\Xi ,  \mc S\mc L)
    \ar[d]^{\mf{sp}^{*}}  &  
   R\Gamma_{c}(\Cht^{\leq \mu}_{\prod_{i\in I} \Frob_{X,\{i\}}^{\deg(v)(m_{i}+d_{i})} (\Delta(\ov{\eta_X}))}/\Xi ,  \mc S\mc L)
    \ar[d]^{\mf{sp}^{*}}
    \\
  R\Gamma_{c}(\Cht^{\leq \mu}_{\prod_{i\in I} \Frob_{X,\{i\}}^{\deg(v)m_{i}} ( \ov{\eta_X^{I}})}/\Xi ,  \mc S\mc L)
    \ar[r]^{\iota_{X}(\delta)^{*}}  &  
    R\Gamma_{c}(\Cht^{\leq \mu}_{\prod_{i\in I} \Frob_{X,\{i\}}^{\deg(v)(m_{i}+d_{i})} ( \ov{\eta_X^{I}})}/\Xi ,  \mc S\mc L)
            } \end{gather}
\end{prop}
\noindent 
Avant de donner la preuve, voici quelques explications sur le diagramme ci-dessus. Les flèches verticales du haut sont la composée
\begin{itemize}
\item des morphismes habituels, pour un trait, de la cohomologie à support compact de la fibre spéciale à coefficients dans les cycles proches vers la cohomologie à support compact de la fibre générique du trait, 
où ici les points génériques des traits sont $\Frob_{X,\{i\}}^{\deg(v)m_{i}} (\Delta(\on{Spec}(\ov K)))$ à gauche et  $\Frob_{X,\{i\}}^{\deg(v)(m_{i}+d_{i})} (\Delta(\on{Spec}(\ov K)))$ à droite, 
\item de l'isomorphisme fourni par   \eqref{plong-FX-FY-K} entre la fibre de la cohomologie à support compact  aux point génériques des traits ci-dessus et la fibre aux points  $\Frob_{X,\{i\}}^{\deg(v)m_{i}} (\Delta(\ov{\eta_X}))$ à gauche et $\Frob_{X,\{i\}}^{\deg(v)(m_{i}+d_{i})} (\Delta(\ov{\eta_X}))$ à droite. 
\end{itemize}
Dans les flèches verticales du bas on désigne encore par $\mf{sp}$ l'image de $\mf{sp}$ par les morphismes de Frobenius partiels $\prod_{i\in I} \Frob_{X,\{i\}}^{\deg(v)m_{i}}$ et  $\prod_{i\in I} \Frob_{X,\{i\}}^{\deg(v)(m_{i}+d_{i})}$. 
 Dans la  ligne du bas  $\iota_{X}$ est le  morphisme  défini dans  \eqref{iota-X}
et  $\iota_{X}(\delta)^{*}$ est l'image inverse par la flèche de spécialisation 
$$\prod_{i\in I} \Frob_{X,\{i\}}^{\deg(v)(m_{i}+d_{i})} ( \ov{\eta_X^{I}})\to \prod_{i\in I} \Frob_{X,\{i\}}^{\deg(v)m_{i}} ( \ov{\eta_X^{I}})$$ associée à $\iota_{X}(\delta)$ et   définie comme  \eqref{fleche-spec-Frob-eta-eta} (en multipliant les $m_{i}$ et $d_{i}$ par $\deg(v)$). 
On note que si 
on compose la   ligne du bas 
 avec les morphismes de Frobenius partiels adéquats
et si on prend la limite inductive sur $\mu$ 
on retrouve  l'action de $\iota_{X}(\delta)\in \on{FWeil} (\eta_X^{I},\ov{\eta_X^{I}})$ 
 sur  $ \varinjlim _{\mu}\restr{\mc H _{ N ^{v}+mv, I, W}^{ \leq\mu,\Lambda}}{\ov{\eta_X^{I}}}  
$ définie après la remarque 8.7 de \cite{coh}  (à la différence près que l'on  considère ici la cohomologie à coefficients dans $\Lambda$ et non dans $E$). 

\begin{rem} Si on prend  la limite inductive sur $\mu$, les deux fl\`eches de gauche deviennent des isomorphismes: celle du haut par le r\'esultat principal de 
\cite{arnaud-cong}, et celle du bas par  \cite{cong-finite}. Nous n'utilisons pas ici ces r\'esultats. 
\end{rem}

\dem 
On applique les théorèmes 2.1 et 8.1 de Orgogozo \cite{orgogozo}
à  $$\mc S\mc L \in D_{c}^{b}(\restr{\Chr_{Y,I,W,(n_{i})_{i\in I},r} ^{mv,(I_{1},...,I_{k})}}{(Y\sm v)^{I}},\Lambda) 
,$$ en prenant comme base $ Y^{I}$. La seule différence par rapport à \cite{orgogozo} est que 
$\Chr_{Y,I,W,(n_{i})_{i\in I},r} ^{mv,(I_{1},...,I_{k})}$ est un champ d'Artin sur $ Y^{I}$ et non pas un schéma. En fait c'est le quotient d'un schéma de type fini par un groupe réductif (agissant trivialement sur $ Y^{I}$ bien sûr). Cela ne pose aucun problème  car 
\begin{itemize}
\item tous les cycles proches sont naturellement dans la catégorie dérivée équivariante, c'est-à-dire associée au champ, 
comme cela est expliqué dans le cas habituel où la base est un trait dans \cite{laumon-moret-bailly} 18.4, 
\item
 les propriétés de constructibilité et de compatibilité au changement de base satisfaites d'après  \cite{orgogozo} se vérifient au niveau des schémas. 
 \end{itemize}
L'article \cite{orgogozo} fournit  une modification $\wt {Y^I}$ de la base $ Y^{I}$, c'est-à-dire un morphisme propre surjectif $\wt {Y^I}\to {Y^I}$ induisant un isomorphisme d'un ouvert partout dense vers un ouvert partout dense, de telle sorte que les conclusions des  théorèmes 2.1 et 8.1 de  \cite{orgogozo} soient satisfaites. On peut supposer de plus que $\wt {Y^I}$ est l'adh\'erence Zariski de son point g\'en\'erique. 
Pour toute flèche de spécialisation $\ov t\to \ov s$ dans $\wt {Y^I}$, en notant $\ov{ s_{0}}$ l'image de $\ov  s$ dans ${Y^I}$, on possède d'après  \cite{orgogozo}
$$\restr{R\Psi(\mc S\mc L)}{(\ov t\to \ov s)} \in D_{c}^{b}(\restr{\Chr_{Y,I,W,(n_{i})_{i\in I},r} ^{mv,(I_{1},...,I_{k})}}{ \ov{ s_{0}}},\Lambda)  . $$

 On note $U$ l'ouvert dense de ${(Y-v)^I}$ (et donc de $\wt {Y^I}$) au-dessus duquel la modification $\wt {Y^I}\to {Y^I}$ est un isomorphisme. 
Evidemment $U$ 
 contient le point générique $\on{Spec} K^{1,...,l}$ (que l'on peut considérer comme un  point asymptote au drapeau de coordonnées qu'on a choisi). On note $s$ le point de $\wt {Y^I}$ au-dessus de 
$s_{0}=\Delta(v)\in {Y^I}$   qui est le ``point exceptionnel du transformé strict'' de ce drapeau, obtenu en appliquant le lemme suivant à la modification $\wt {Y^I}\to {Y^I}$ et au drapeau des coordonnées.  On note que $s$ et $s_{0}$ sont tous les deux égaux à $\on{Spec} \mf k$. 

\begin{lem} \label{drapeau} Soit $\wt S\to S$ une modification, où $S$ est une variété lisse irréductible de dimension $l$ sur $\mf k$ et  $\wt S$ est l'adh\'erence Zariski de son point g\'en\'erique. On se donne
un drapeau complet  \begin{gather}\label{drap-non-modif}
S_{0}\subset S_{1} \subset ... \subset S_{l}=S\end{gather} de sous-variétés lisses fermées irréductibles, avec $S_{j}$ de dimension $j$ pour tout $j$, et $S_{0}$ un point fermé. Alors il existe un unique drapeau 
\begin{gather} \label{drap-modif}\wt S_{0}\subset \wt S_{1} \subset ... \subset \wt S_{l}=\wt S\end{gather} de sous-schémas fermés  de $\wt S$ tel que 
\begin{enumerate}
\item pour tout $j\in \{0,...,l\}$, $\wt S_{j}$ est une modification de $S_{j}$, et  un isomorphisme au-dessus d'un ouvert $S_{j}^{0}$ de $S_{j}$ dont le fermé complémentaire est de codimension $\geq2$, et $\wt S_{j}$ est  l'adh\'erence Zariski de son point g\'en\'erique
\item pour tout $j\in \{0,...,l-1\}$, le point générique de $ S_{j}$ appartient à $ S_{j+1}^{0}$ (puisque le fermé complémentaire est de codimension $\geq2$), et $\wt S_{j}$ est l'adhérence de Zariski dans $\wt  S_{j+1}$ de l'unique  relèvement à $\wt  S_{j+1}$ du point générique de $ S_{j}$. 
\end{enumerate}
\end{lem}
La dernière condition implique en particulier que $\wt S_{0}\to S_{0}$ est un isomorphisme (les deux sont isomorphes \`a $\on{Spec} \mf k$). Le drapeau modifié \eqref{drap-modif} est appelé le ``transformé strict'' du drapeau 
\eqref{drap-non-modif} et on dit que $\wt S_{0}$ est le ``point exceptionnel de ce transformé strict''.

\dem La construction du drapeau \eqref{drap-modif} se fait par récurrence descendante sur $j$.   \cqfd


\noindent{\bf Suite de la démonstration de la \propref{prop-boite-noire}.}
 On note $\ov s$ un point géométrique algébrique au-dessus de $s$ muni d'une 
flèche de spécialisation  \begin{gather}
\label{spec-K1..k-s}\on{Spec} \ov {K^{1,...,l}}\to \ov s\end{gather} 
dans $\wt {Y^I}$, \'etendant le morphisme naturel de $\on{Spec}  K^{1,...,l}$ dans l'hens\'elis\'e de $\wt {Y^I}$ en $s$, qui r\'esulte du \lemref{drapeau} (en effet le morphisme de $K^{1,...,l}$ vers le corps local compl\'et\'e de dimension sup\'erieure associ\'e au drapeau  \eqref{drap-modif}, ou plus pr\'ecis\'ement \`a la valuation $v_\bullet$ que nous d\'efinirons au d\'ebut de la d\'emonstration du lemme \ref{lem-assez},   est un isomorphisme). 
On note que tout automorphisme de $\ov s$ 
respectant \eqref{spec-K1..k-s} est l'identité. 
On voit que  $\on{Gal}(\ov {K^{1,...,l}}/K^{1,...,l})$ agit par automorphismes sur la flèche \eqref{spec-K1..k-s}.  

De plus pour tout $(m_{i})_{i\in I}\in \N^{I}$ on en déduit une flèche de spécialisation dans $\wt {Y^I}$
\begin{gather}\label{K1k-s-spec-Frob}
\prod_{i\in I} \Frob_{Y,\{i\}}^{m_{i}}  (\on{Spec} \ov{K^{1,...,l}})\to
  \Frob_{\mf k} ^{\sum_{i\in I }m_{i}} ( \ov s)=\ov s. 
\end{gather}

Pour tout $(m_{i})_{i\in I}\in \N^{I}$ on définit   une flèche de spécialisation dans $\wt {Y^I}$
\begin{gather}\label{eta-s-spec-Frob}
\prod_{i\in I} \Frob_{Y,\{i\}}^{m_{i}}  (\ov{\eta_{Y}^{I}} )\to
   \ov s 
  \end{gather} comme composée  
\begin{itemize}
\item de l'image par $\prod_{i\in I} \Frob_{Y,\{i\}}^{m_{i}}  $ de la flèche de spécialisation 
\begin{gather}\label{spec-eta-K1..k}\ov{\eta_{Y}^{I}} \to \on{Spec} \ov {K^{1,...,l}}\end{gather}
associée à l'inclusion $ \ov{F_Y^{I}}\subset 
     \ov{K^{1,...,l}} $
  fixée dans   \eqref{diag-plong} (en effet grâce à cette inclusion, $\on{Spec} \ov{F_Y^{I}}$ est l'hensélisé de $\ov{\eta_{Y}^{I}}$ en $\on{Spec}  \ov{K^{1,...,l}} $), 
\item et de la  flèche de spécialisation  \eqref{K1k-s-spec-Frob}. 
\end{itemize}

\begin{lem} \label{lem-assez}
Pour tout $(m_{i})_{i\in I}\in \N^{I}$ assez croissant, on a 
\begin{itemize}
\item $\prod_{i\in I} \Frob_{Y,\{i\}}^{m_{i}} (\Delta(  \eta_{Y} ))$ appartient à $U$, 
\item $s$ est le point exceptionnel du transformé strict dans $\wt {Y^I}$ de l'image de la diagonale par $\prod_{i\in I} \Frob_{Y,\{i\}}^{m_{i}}$ (ce qui a un sens grâce à la condition précédente et au fait que la modification $\wt {Y^I}\to {Y^I}$ est un isomorphisme au-dessus de $U$).  
\end{itemize}
\end{lem} 

On commence par des pr\'eliminaires pour la d\'emonstration du \lemref{lem-assez}. 

On munit  $\Z^l$ de l'ordre lexicographique (inverse), pour lequel $(\mu_i)\geq (\nu_i)$ s'ils sont \'egaux ou s'il existe $i\in\{1,...,l\}$ avec 
$\mu_i>\nu_i$ et $\mu_j=\nu_j$ pour $j\in \{i+1,...,l\}$.

On reprend la situation et les notations du \lemref{drapeau}: 
 $S_{0}\subset S_{1} \subset ... \subset S_{l}=S$ est un drapeau de sous-variétés lisses fermées irréductibles et  
$\wt S_{0}\subset \wt S_{1} \subset ... \subset \wt S_{l}$ un drapeau de sous-schémas fermés  de $\wt S$. 
On rappelle que $S_{0}=\wt S_{0}=\on{Spec} \mf k$. 

On note $ S_0=S_{0}^\lo  \subset S_{1}^\lo \subset ... \subset S_{l}^\lo=S^\lo$ les localis\'es en le point $ S_0$ des $S_i$. 
On note $\wt S_{i}^\lo={\wt S_{i}}\times  _{S_i} S_i^\lo$, si bien que $\wt S_{i}^\lo\to S_i^\lo$ est une modification. 

On note $R$ la  $\mf k$-algèbre locale régulière noethérienne de corps résiduel $\mf k$ et d'id\'eal maximal $\mf m$
telle que $S^\lo=\on{Spec} R$. Soient $t_1,\dots,t_l$  dans $\mf m$ donnant une base sur $\mf k$ de $\mf m/{\mf m}^2$, tels que 
$S_i^\lo=\mb V(t_{i+1},\dots,t_l)$ pour $0\leq i\leq l$. Le drapeau
$$S_\bullet^\lo=\big( \{s_0\}=S_{0}^\lo\subset S_{1}^\lo \subset ... \subset S_{l}^\lo=S^\lo\big)$$
définit une valuation lexicographique
$$( v_1,\dots, v_l)=v_\bullet:  R\setminus \{0\}\to \N^l$$ sur $R$,  de la manière suivante. Procédant par récurrence sur $l$, on suppose la valuation $(v_{1,R/(t_l)},\dots,v_{l-1, R/(t_l)})$ déjà définie sur $R/(t_l)$ (pour $l=1$, $R$ est un anneau de valuation discrète et $v_\bullet$ est bien entendu sa valuation). Pour $f\in R\setminus \{0\}$ on pose alors
$$v_l(f)=v_{t_l}(f)$$
(où $v_{t_l}$ est la valuation $t_l$-adique) et, pour $1\leq i\leq l-1$, 
$$v_i(f)=v_{i, R/(t_l)}\left(\frac{f}{t_l^{v_l(f)}}\right).$$
Cette valuation lexicographique ne dépend que du drapeau $S_\bullet^\lo$, dans la mesure où pour une suite régulière $(t'_1,\dots,t'_l)$ définissant le même drapeau $S_\bullet^\lo$ les images de $t_i$ et $t'_i$ dans $R/(t_{i+1},\dots,t_l)$ sont associées.
On notera encore $$v_\bullet : \on{Frac} R\to \Z^l\cup \{+\infty\}$$
son extension au corps des fractions. Le sous-anneau de $\on{Frac} R$ formé des éléments de valuation $v_\bullet$ positive ou nulle (pour l'ordre lexicographique) sera noté $\mc O_{S_\bullet^\lo}$ et appelé anneau local du drapeau.

Pour $h\in \N^*$, on dit que  $(\mu_{i})_{i\in I}\in (\N^*)^{I}$ est $h$-croissant  si $\mu_1>0$ et pour tout $i=1,...,l-1$, on a
$\mu_{i+1}\geq h  \mu_i$.

\begin{lem}\label{poly-formel} Pour tout $f\in R \setminus \{0\}$ il existe $h$ tel que pour 
 tout $(\mu_{i})_{i\in I}\in (\N^*)^{I}$ $h$-croissant et  tout morphisme   $$c :R\to V$$ 
 (où $V$ est à la fois une $\mf k$-algèbre et un anneau de valuation discrète de corps résiduel $\mf k$, de valuation $v_V$) tel que   
 $$v_V(t_i)=\mu_i, \forall\, 1\leq i\leq l,$$ on ait  
 $$v_V(c(f))=\sum_ {i=1}^l \mu_i v_i(f).$$\end{lem} 

\dem Remplaçant les anneaux locaux $V$ et $R$ par leurs complétés $\widehat V\simeq \mf k[[t]]$ et $\widehat R\simeq \mf k[[t_1,\dots,t_l]]$ on se ramène au cas particulier où $V$ et $R$ sont complets. On note 
$$f=\sum_{(i_1,\dots,i_l)}a_{(i_1,\dots,i_l)}t_1^{i_1}\cdots t_l^{i_l}$$
et on considère le plus petit élément (pour l'ordre lexicographique) $(\alpha_1,\dots,\alpha_l)$ tel que $a_{(\alpha_1,\dots,\alpha_l)}\not=0$ (on note que $\alpha_i=v_i(f)$).
Les autres mon\^omes non nuls intervenant dans $f$  sont tous multiples de $t_1 (\prod_{i=1}^l t_i^{\alpha_i})$ ou 
$t_2 (\prod_{i=2}^l t_i^{\alpha_i})$, $t_3 (\prod_{i=3}^l t_i^{\alpha_i})$, ..., ou $t_l t_l^{\alpha_l}$. 
Par conséquent, tout $h> \alpha_1 +\alpha_2+... +\alpha_{l-1}$ convient. \cqfd

On utilisera le lemme pr\'ec\'edent en combinaison avec le lemme purement combinatoire suivant. 

\begin{lem}\label{lem-combinatoire} Pour tout  $(\alpha_1, ..., \alpha_l)\in \Z^l$ positif (resp. strictement positif) pour l'ordre lexicographique, 
il existe $h$ tel que pour 
 tout $(\mu_{i})_{i\in I}\in (\N^*)^{I}$ $h$-croissant  on ait 
  $$ \sum_ {i=1}^l\mu_i \alpha_i \geq 0 \text{ (resp. }   >0 \text{ ) }.$$
  \end{lem} 
  \cqfd

On note $s$ l'unique point de $\wt S_0^\lo$.

On choisit un ouvert affine $\mathcal U=\on{Spec}  A$ de $\widetilde S^\lo$ contenant le point $s$. La $R$-algèbre $A$ est de type fini et est contenue dans $\on{Frac} R$.  

\begin{lem} \label{A-val} On a $A\subset \mc O_{S_\bullet}$ et un élément de $A$ s'annule au point $s$ si et seulement s'il est de valuation $v_\bullet$ strictement positive.\end{lem}
\dem On procède par récurrence sur $l$. Le cas où $l=1$ est tautologique, car la modification $\widetilde S^\lo\to S^\lo$ est alors un isomorphisme.

On suppose le lemme déjà vérifié pour la modification $\widetilde S_{l-1}^\lo\to S_{l-1}^\lo$ résultant de la condition (1) du \lemref{drapeau} et pour l'ouvert affine $\mathcal  U\cap \widetilde S_{l-1}^\lo$ de $\widetilde S_{l-1}^\lo$.  

Soit $f\in A$. Par la condition (2) du \lemref{drapeau}, la fraction rationnelle $f\in \on{Frac} R$ est définie au point générique de $S_{l-1}^\lo$ et on a donc $v_l(f)\geq 0$. 

Lorsque $v_l(f)>0$ le lemme est déjà démontré pour $f$. 

Lorsque $v_l(f)=0$, la valeur de $f$ en ce point générique (vue comme élément de $\on{Frac} R/(t_l)$) est la restriction $f_{\widetilde S_{l-1}}$ de $f\in A$ au fermé $\wt S_{l-1}^\lo$. 
En particulier, comme le lemme est supposé déjà vérifié pour $\wt S_{l-1}$, on a $$(v_\bullet)_{R/(t_l)} (f_{\widetilde S_{l-1}^\lo})\geq (0,\dots,0)$$ et cette inégalité est stricte si et seulement si $f$ s'annule en $s$, ce qui montre le lemme pour la modification $\wt S^\lo\to S^\lo$ et l'ouvert affine $\mathcal  U=\on{Spec}  A$.\cqfd

La $R$-algèbre $A$ étant de type fini, on en choisit une famille génératrice 
$$\frac{f_1}{g_1}, \dots, \frac{f_N}{g_N}$$
(où, pour tout $1\leq j\leq N$, on a $f_j,g_j\in R$ et $g_j\not =0$). Au-dessus de l'ouvert affine (non-vide) $\mc D_{g_1\cdots g_N}=\on{Spec}  R[(g_1\cdots g_N)^{-1}]$ de $S^\lo$, le morphisme $\mathcal  U\to S^\lo$ est un isomorphisme.

\begin{lem} \label{courbe-drapeau} Il existe $h$ tel que pour 
 tout $(\mu_{i})_{i\in I}\in (\N^*)^{I}$ $h$-croissant et  tout morphisme   $$c :R\to V$$ 
 (où $V$ est à la fois une $\mf k$-algèbre et un anneau de valuation discrète de corps résiduel $\mf k$, de valuation $v_V$) tel que   
 $$v_V(t_i)=\mu_i, \forall\, 1\leq i\leq l,$$ on ait  
  \begin{itemize}
\item a)  le morphisme $\on{Spec}  c$ envoie le point générique de $\on{Spec}  V$ dans l'ouvert $\mc D_{g_1\cdots g_N}$ ;
\item b) ce morphisme se rel\`eve (de  mani\`ere unique) en un morphisme $\on{Spec}  V\to \mathcal  U=\on{Spec}  A$ envoyant le point fermé de $\on{Spec}  V$ sur $s$.
\end{itemize}
\end{lem} 

\dem 
Appliquant le \lemref{poly-formel} à l'élément $g_1\cdots g_N$, 
on obtient déjà un entier $h$ tel que sous l'hypothèse du \lemref{courbe-drapeau} on ait $c(g_1\cdots g_N)\not=0$, ce qui vérifie déjà le point (a) en fournissant un morphisme $R[(g_1\dots g_N)^{-1}]\to \on{Frac} V$.

Dans b) l'unicit\'e du rel\`evement r\'esulte de a) car  la modification est un isomorphisme au-dessus de $\mc D_{g_1\cdots g_N}$. 

Soit $h$ un entier suffisamment grand pour que le \lemref{poly-formel} s'applique à tous les éléments $f_j$ et $g_j\ (1\leq j\leq N)$. 
D'après la conclusion de ce lemme, gr\^ace au \lemref{A-val} et  
quitte à augmenter encore $h$ (gr\^ace au \lemref{lem-combinatoire}), on a alors, pour tout $1\leq j\leq N$,  $v_V\left(\frac{f_j}{g_j}\right)=\sum_ {i=1}^l\mu_i v_i \left(\frac{f_j}{g_j}\right)\geq 0$, ce qui fournit le morphisme $A\to V$ désiré. 

Enfin,  quitte à augmenter encore $h$, on applique une dernière fois le \lemref{poly-formel} (et le \lemref{lem-combinatoire})
à des générateurs (en nombre fini, puisque $R$ et donc aussi $A$ sont noethériens) de l'idéal $I$ de $A$  formé des fonctions s'annulant en $s$ et on obtient $v_V(c(I))>0$. En particulier, l'image réciproque par $c$ de l'idéal maximal de $V$ contient l'idéal (maximal) $I$  de $A$. 
\cqfd

Ceci termine la preuve du \lemref{lem-assez}. En effet $\mu_i=q^{m_i}$ est $h$-croissant pour $h$ assez grand si et seulement si $(m_i)$ est assez croissant au sens de la notation \ref{assez-croissant}. 

\cqfd

\noindent{\bf Fin de la démonstration de la \propref{prop-boite-noire}.}
On choisit $(m_{i})_{i\in I}\in \N^{I}$ tel que 
$(m_{i}+d_{i})_{i\in I}\in \N^{I}$
  et assez croissant pour que les deux propriétés du lemme précédent  soient vérifiées à la fois par $(m_{i})_{i\in I}$ et par $(m_{i}+d_{i})_{i\in I}$. 

Enfin   $\delta\in \on{FWeil}(\ov{K^{1,...,l}}/K^{1,...,l})$ donne lieu à  un diagramme commutatif de flèches de spécialisation dans $\wt {Y^I}$
  $$   \xymatrix{
       \prod_{i\in I} \Frob_{Y,\{i\}}^{m_{i}+d_{i}}  (\ov{\eta_{Y}^{I}})    \ar[d]^{} \ar[rr]^{\iota_{Y}(\delta) } 
       & &  \prod_{i\in I} \Frob_{Y,\{i\}}^{m_{i} }  (\ov{\eta_{Y}^{I}}) \ar[d]_{ } 
     \\
     \prod_{i\in I} \Frob_{Y,\{i\}}^{m_{i}+d_{i}}  (\on{Spec} \ov {K^{1,...,l}})  \ar[dr]^{ } \ar[rr]^{\delta } 
       & &  \prod_{i\in I} \Frob_{Y,\{i\}}^{m_{i} }  (\on{Spec} \ov {K^{1,...,l}}) \ar[dl]_{ } 
     \\
     &   \ov s  &  
    }$$     
où les deux flèches obliques sont \eqref{K1k-s-spec-Frob} pour $(m_{i})_{i\in I}$ et   $(m_{i}+d_{i})_{i\in I}$, et les deux flèches verticales sont   les images par
$ \prod_{i\in I} \Frob_{Y,\{i\}}^{m_{i}+d_{i}}$ et $ \prod_{i\in I} \Frob_{Y,\{i\}}^{m_{i} }$ 
de \eqref{spec-eta-K1..k}. On note que \eqref{eta-s-spec-Frob} est la composée descendante à droite. Le diagramme commutatif ci-dessus associe donc à $\delta$  un morphisme image inverse 
 \begin{gather}\label{morph-delta-fleches}
 \restr{R\Psi(\mc S\mc L)}{\big(\prod_{i\in I} \Frob_{Y,\{i\}}^{m_{i}}  (\ov{\eta_{Y}^{I}}) \to  \ov s \big)}  \to \restr{R\Psi(\mc S\mc L)}{\big(\prod_{i\in I} \Frob_{Y,\{i\}}^{m_{i}+d_{i}}  (\ov{\eta_{Y}^{I}}) \to  \ov s \big)}  
\end{gather} (qui est un isomorphisme). 
%
En vertu de la compatibilité des cycles  proches au changement de base assuré par le théorème 2.1 de \cite{orgogozo}, 
on a \begin{gather}\label{orgo-psimi}\restr{R\Psi(\mc S\mc L)}{\big(\prod_{i\in I} \Frob_{Y,\{i\}}^{m_{i}} (\Delta( \ov{\eta_{Y}})) \to  \ov s \big)}=R\Psi_{(m_{i})_{i\in I}} \mc S\mc L\end{gather} 
et de même avec $(m_{i}+d_{i})_{i\in I}$. Enfin l'image de $\mf{sp}$ par 
$\prod_{i\in I} \Frob_{Y,\{i\}}^{m_{i}} $ fournit, grâce au caract\`ere ULA sur $U$ des faisceaux consid\'er\'es, et \`a un th\'eor\`eme non publi\'e de Gabber  rappel\'e dans le th\'eor\`eme A.9 de \cite{arnaud-cong}, appliqu\'e \`a la suite de fl\`eches de sp\'ecialisation $$\prod_{i\in I} \Frob_{Y,\{i\}}^{m_{i}}  (\ov{\eta_{Y}^{I}}) \to \prod_{i\in I} \Frob_{Y,\{i\}}^{m_{i}}  (\Delta( \ov{\eta_{Y}})) \to  \ov s,$$ 
  un isomorphisme 
\begin{gather}\label{morph-diag-fleches}
\restr{R\Psi(\mc S\mc L)}{\big(\prod_{i\in I} \Frob_{Y,\{i\}}^{m_{i}}  (\Delta( \ov{\eta_{Y}})) \to  \ov s \big)}\isom 
\restr{R\Psi(\mc S\mc L)}{\big(\prod_{i\in I} \Frob_{Y,\{i\}}^{m_{i}}  (\ov{\eta_{Y}^{I}}) \to  \ov s \big)}  
\end{gather}
 et de même avec $(m_{i}+d_{i})_{i\in I}$.  
 En composant (de droite à gauche) 
   \begin{itemize}
 \item l'inverse de \eqref{orgo-psimi}, 
  \item le morphisme  \eqref{morph-diag-fleches}, 
  \item le morphisme  
   \eqref{morph-delta-fleches},
   \item  l'inverse de  \eqref{morph-diag-fleches}   pour $(m_{i}+d_{i})_{i\in I}$ 
   \item et \eqref{orgo-psimi} pour $(m_{i}+d_{i})_{i\in I}$ 
   \end{itemize}
   on obtient 
   $T_{\delta}$ tel que le diagramme commute. Il est clair que  $T_{\delta}$ ne dépend que des données locales au sens hensélien, donc au sens du complété en $v$, d'apr\`es la remarque suivante.  
 \cqfd

\begin{rem}\label{rem-hensel-complete}
Pour les données que l'on considère, 
   la localité au sens hensélien et  celle au sens du complété sont équivalentes. 
\end{rem}

\begin{rem}\label{rem-complement-Tdelta}
En complément on indique deux propriétés des $T_{ \delta}$ qui résultent assez facilement de la démonstration mais ne seront pas utilisées dans la suite: 
\begin{itemize}
\item si on remplace $(m_{i})_{i\in I}$ par $(m_{i}+\alpha_{i})_{i\in I}$ avec $(\alpha_{i})_{i\in I}\in \N^{I}$ croissant, le   
$T_{ \delta}$ associé à $(m_{i}+\alpha_{i})_{i\in I}$ est entrelacé avec celui associé à $(m_{i})_{i\in I}$ par les morphismes de Frobenius partiels 
(plus précisément les correspondances cohomologiques   \eqref{corresp-coho-frob-partiels},  intercalées avec les correspondances cohomologiques \eqref{oubli-cycles-proches} et leurs inverses  pour remettre  les pattes dans l'ordre et pouvoir itérer les morphismes de Frobenius partiels), 
\item $T_{ \delta}$  est 
compatible  à la composition de $\delta$ 
(lorsque  le $I$-uplet $(m_{i})_{i\in I}$ de départ 
est choisi assez positif et croissant 
pour que la composition soit possible). 
\end{itemize}
\end{rem}

\begin{rem} \label{rem-RpsiRpsi} Par le r\'esultat principal de \cite{arnaud-cong}, 
le choix du plongement $ \ov{F_Y^{I}}
   \subset \ov{K^{1,...,l}} $ de \eqref{diag-plong} identifie   $\restr{R\Psi(\mc S\mc L)}{(\ov{\eta_{Y}^{I}} \to\ov  s)}$ à $R\Psi \cdots R\Psi(\mc S\mc L)$ (cycles proches au sens usuel appliqués à chacune des coordonnées dans l'ordre $1,...,l$ de la gauche vers la droite). Cela peut aider à comprendre intuitivement l'action de $\delta$ et la construction de $T_{\delta}$, mais nous n'avons pas besoin ici de ce r\'esultat. 
   \end{rem}

\section{Une variante locale du lemme de Drinfeld}
 
 Le but de ce paragraphe est   de montrer le lemme suivant. 
 
 \begin{lem}\label{lem-surj-Weil-loc}
 Le morphisme \eqref{morph-Weil-loc} 
 $:\on{FWeil}(\ov{K^{1,...,l}}/K^{1,...,l})
   \to   (\on{Weil}(\ov K/K))^{I}$
 (qui est aussi la flèche verticale de droite  du diagramme \eqref{diag-gpes-Weil}) est surjectif. 
 \end{lem}

On pourrait  le montrer à l'aide de polynômes d'Eisenstein. Cependant il est plus naturel de le déduire de la  variante suivante du lemme de Drinfeld pour les corps locaux de dimension supérieure. 

Soit $\mf k$ un corps fini et   
$\mf K$ un corps  contenant $\mf k$. On note $\mr{Frob} _{\mf K/\mf k} :x\mapsto x^{\sharp \mf k}$ le morphisme de Frobenius  de $\mf K$ relatif à $\mf k$. 

\begin{lem}\label{lem-OL-OM-descente}
 On a une équivalence entre  
\begin{itemize}
\item la catégorie des algèbres  finies étales $L$ sur  $\mf K((t))$ avec relèvement $F_{1}$ de $\mr{Frob} _{\mf K/\mf k}\otimes \Id_{t}$ (c'est-à-dire un isomorphisme $\tav L\simeq L$, où $\tav L=L\otimes_{\mf K((t)),\mr{Frob} _{\mf K/\mf k}\otimes \Id_{t}} \mf K((t))$), 
\item la catégorie des algèbres  finies étales $M$ sur $\mf k((t))$ munies d'une action 
de $\mr{Gal}(\mf K^{\mr{sep}}/\mf K)$. 
\end{itemize}
Cette équivalence est caractérisée par un isomorphisme 
$$L \otimes _{\mf K((t))}   \mf K^{\mr{sep}}((t)) \simeq 
 M\otimes _{\mf k((t))}  \mf K^{\mr{sep}}((t)) $$
entrela\c cant  $F_{1}$ et $\Id_{M}\otimes (\mr{Frob} _{\mf K/\mf k}\otimes \Id_{t})$, et compatible aux  actions de $\mr{Gal}(\mf K^{\mr{sep}}/\mf K)$. Autrement dit on a 
\begin{gather}\label{fonct-M-L}
L=\Big( M\otimes _{\mf k((t))} \mf K^{\mr{sep}}((t)) \Big)^{\mr{Gal}(\mf K^{\mr{sep}}/\mf K)}\end{gather}
et 
$$M=\Big(L \otimes _{\mf K((t))}   \mf K^{\mr{sep}}((t))\Big)^{F_{1}}.$$
\end{lem}

On note $\mc O_{L} $ et $\mc O_{M}$ les anneaux d'entiers de $L$ et $M$. 
En complément de  la correspondance ci-dessus  
on a aussi 
\begin{gather}\nonumber\mc O_{L}  \otimes _{\mf K[[t]]}   \mf K^{\mr{sep}}[[t]] \simeq 
 \mc O_{M}\otimes _{\mf k[[t]]}  \mf K^{\mr{sep}}[[t]]   \end{gather}
et pour tout entier $k\in \N$, 
\begin{gather}\label{OLOM}(\mc O_{L}/t^{k}  \mc O_{L} ) \otimes _{\mf K }   \mf K^{\mr{sep}}  \simeq 
 (\mc O_{M}/t^{k}  \mc O_{M} )\otimes _{\mf k }  \mf K^{\mr{sep}} .   \end{gather}
 On va montrer en même temps l'énoncé du lemme et \eqref{OLOM}. 
 
 On commence par rappeler qu'on a une équivalence entre 
 \begin{itemize}
 \item la catégorie des $(E,\theta)$ où $E$ est un $\mf K $-espace  vectoriel  de dimension finie muni d'un isomorphisme $\theta:\tav E\simeq E$, 
 \item la catégorie des $\mf k$-espaces vectoriels $D$ de dimension finie munis d'une action linéaire de $\mr{Gal}(\mf K^{\mr{sep}}/\mf K)$,  
\end{itemize}
et cette dernière équivalence est caractérisée par un isomorphisme 
$$E\otimes _{\mf K}\mf K^{\mr{sep}}\simeq D\otimes _{\mf k}\mf K^{\mr{sep}}$$
compatible au Frobenius et à l'action de $\mr{Gal}(\mf K^{\mr{sep}}/\mf K)$, de sorte que 
$$E=\big( D\otimes _{\mf k}\mf K^{\mr{sep}}\big) ^{\mr{Gal}(\mf K^{\mr{sep}}/\mf K)} \text{ \ \ 
et \ \ }D=\big(E\otimes _{\mf K}\mf K^{\mr{sep}}\big)^{\Frob}.$$

Pour montrer le \lemref{lem-OL-OM-descente} on remarque d'abord que le foncteur, donné par \eqref{fonct-M-L}, de la deuxième catégorie vers la première dans l'énoncé du lemme,  est pleinement fidèle
et qu'on a alors 
$$
\mc O_{L}=\Big( \mc O_{M}\otimes _{\mf k[[t]]} \mf K^{\mr{sep}}[[t]] \Big)^{\mr{Gal}(\mf K^{\mr{sep}}/\mf K)}$$
et, pour tout $k$, 
$$
(\mc O_{L}/t^{k}\mc O_{L}) =\Big( (\mc O_{M}/t^{k}\mc O_{M})\otimes _{\mf k[[t]]} \mf K^{\mr{sep}}[[t]] \Big)^{\mr{Gal}(\mf K^{\mr{sep}}/\mf K)}.$$
Ce foncteur \eqref{fonct-M-L} est essentiellement surjectif car on a le foncteur inverse suivant. Etant donné $L$ avec relèvement de Frobenius partiel dans la catégorie de départ, comme 
 $L$ est une alg\`ebre finie \'etale sur  $\mf K((t))$, $F_1$ envoie $\mc O_L$ (cl\^oture int\'egrale de $\mf K[[t]]$ dans $L$) dans lui-m\^eme et induit un isomorphisme   de $\tav \mc O_L=\mc O_L \otimes_{\mf K[[t]],\mr{Frob} _{\mf K/\mf k}\otimes \Id_{t}} \mf K[[t]]$ dans  $\mc O_L$, et 
 on construit $M$ par la formule $$\mc O_{M}=\varprojlim  \Big((\mc O_{L}/t^{k}  \mc O_{L} ) \otimes _{\mf K }   \mf K^{\mr{sep}}\Big)^{F_{1}}.$$
On montre que 
$L$ est étale sur $\mf K((t))$ si et seulement si $M$ est étale sur $\mf k((t))$.
Ceci termine la preuve  du \lemref{lem-OL-OM-descente}. \cqfd

En complément on remarque que la deuxième catégorie du 
\lemref{lem-OL-OM-descente} est trivialement équivalente
\begin{itemize}
\item à la catégorie des ensembles finis munis d'une action de $\mr{Gal}(\mf K^{\mr{sep}}/\mf K)\times \mr{Gal}\big(\big(\mf k((t))\big)^{\mr{sep}}/\mf k((t))\big)$, 
\item à la catégorie des algèbres  finies étales $M$ sur  $\mf K$ munies d'une action 
de $\mr{Gal}\big(\big(\mf k((t))\big)^{\mr{sep}}/\mf k((t))\big)$. 
\end{itemize}
On a donc une équivalence entre la première catégorie du 
\lemref{lem-OL-OM-descente} et la dernière catégorie ci-dessus. 
Par récurrence sur l'entier $l$ on en déduit  le résultat suivant. 

Pour tout corps fini $\mf k$ on considère les morphismes de Frobenius partiels  
$F_{i}$ sur $\mf k((t_{1}))...((t_{l}))$, définis par $F_{i}(t_{i})=t_{i}^{\sharp \mf l}$ et $F_{i}(t_{j})=t_{j}$ pour $j\neq i$. 

\begin{lem}
Soit $\mf k$ un corps fini. Il y a une équivalence de catégories
entre
\begin{itemize}
\item les algèbres étales sur $\mf k((t_{1}))...((t_{l}))$ munies de relèvement 
des morphismes de Frobenius partiels, commutant entre eux et dont le produit est le morphisme de Frobenius total, 
\item  les ensembles finis munis d'actions continues de 
$\prod_{i=1}^{l}\on{Gal}\big(\big(\mf k((t_{i}))\big)^{\mr{sep}}/\mf k((t_{i})) \big)$. 
\end{itemize}
 \end{lem}

\noindent {\bf Démonstration du \lemref{lem-surj-Weil-loc}.}
On déduit du lemme précédent que le morphisme \eqref{morph-Weil-loc}  est d'image dense (car il induit un isomorphisme des compl\'et\'es profinis). Or il respecte les degrés à valeurs dans $\Z^{I}$, donc  
$$\Ker\big(\on{FWeil}(\ov{K^{1,...,l}}/K^{1,...,l})\to \Z^{I}\big)\to 
\Big( \Ker\big(\on{ Weil}(\ov{K}/K)\to \Z\big)\Big)^{I}$$ est 
un morphisme d'image dense dont l'espace de départ est 
profini, donc compact, donc d'image fermée. Par conséquent ce morphisme est surjectif et  \eqref{morph-Weil-loc}  est surjectif. 
\cqfd

\section{Diagrammes de  correspondances cohomologiques. }
 
  A la suite de  \cite{sga5,varshavsky-fujiwara,brav-var} et comme dans le paragraphe 4.1 de \cite{coh} on appelle correspondance de $X_{1}$ vers $X_{2}$ un  morphisme $(a_{1},a_{2}):A\to X_{1}\times X_{2}$ de champs algébriques localement de type fini sur $\Fq$ tel que $a_{2}$ soit  représentable et de type fini. Alors, pour $\mc F_{1}\in D_{c}^{b}(X_{1}, E)$ et $\mc F_{2}\in D_{c}^{b}(X_{2}, E)$, un morphisme 
  $$u:(a_{2})_{!}(a_{1})^{*}(\mc F_{1})\to \mc F_{2}, \text{ \ ou, de fa\c con équivalente \ } 
 u:a_{1}^{*}(\mc F_{1})\to a_{2}^{!}(\mc F_{2})$$ 
  est appelé une correspondance cohomologique   de 
  $\mc F_{1}$ vers $\mc F_{2}$ (ou de 
    $(X_{1},\mc F_{1})$ vers $(X_{2},\mc F_{2})$) supportée par $A$.
 
 
   Si $(a'_{2},a'_{3}):A'\to X_{2}\times X_{3}$ 
   est une autre correspondance, 
  $\mc F_{3}\in D_{c}^{b}(X_{3}, E)$
   et $u':(a'_{3})_{!}(a'_{2})^{*}(\mc F_{2})\to \mc F_{3}$ est une correspondance cohomologique   de 
     $\mc F_{2}$ vers $\mc F_{3}$ supportée par $A'$, la  composée  $u'\circ u$ est obtenue de la fa\c con suivante. 
     On pose $\wt A=A\times_{X_{2}}A'$. 
     On a un diagramme  commutatif
          $$   \xymatrix{
    &  &   \wt A\ar[dl]^{\wt {a'_{2}}} \ar[dr]_{\wt {a_{2}}} 
       & &
     \\
    &  A\ar[dl]^{a_{1}} \ar[dr]_{a_{2}}&    & A'\ar[dl]^{a'_{2}} \ar[dr]_{a'_{3}} &
     \\
  X_{1}     &     & X_{2} 
  &&  X_{3}
    }$$     
    où le carré est cartésien.      D'où une  correspondance 
     $ 
     (a_{1}\wt {a'_{2}},{a'}_{3}\wt {a_{2}}):\wt A\to X_{1}\times X_{3}$. On définit la correspondance cohomologique   $\wt u=u'\circ u$ de $\mc F_{1}$ vers $\mc F_{3}$ supportée par  $\wt A$ comme la  composée
     \begin{gather*}(a'_{3}\wt {a_{2}})_{!}(a_{1}\wt {a'_{2}})^{*}(\mc F_{1})=  (a'_{3})_{!}  (\wt {a_{2}})_{!} (\wt {a'_{2}} )^{*} (a_{1})^{*} (\mc F_{1}) \\ 
     \simeq (a'_{3})_{!}(a'_{2})^{*}(a_{2})_{!}(a_{1})^{*}(\mc F_{1})\xrightarrow{u} 
     (a'_{3})_{!}(a'_{2})^{*}(\mc F_{2})\xrightarrow{v}  \mc F_{3}, \end{gather*}
     où l'isomorphisme vient du changement de base propre $(\wt a_{2})_{!}(\wt {a'}_{2}) ^{*} 
    \simeq    (a'_{2})^{*}(a_{2})_{!}  $. 
   
   \begin{rappel}\label{rappel-coho}
Une correspondance cohomologique    $((a_{1},a_{2}), u)$ de $(X_1, \mc F_{1})$ vers $(X_2, \mc F_{2})$ comme ci-dessus induit, sous l'hypothèse que $a_{1}$ est propre, un  morphisme sur les  cohomologies à  support compact. En effet soit $Y$  un  champ et 
  $p_{i}:X_{i}\to Y$ (pour  $i=1,2$) des morphismes tels que $p_{1}\circ a_{1}=p_{2}\circ a_{2}$. 
    Alors 
la  correspondance cohomologique  $u$ induit un  morphisme 
  $H(u): p_{1,!}(\mc F_{1}) \to p_{2,!}(\mc F_{2})$ dans  $D_{c}^{b}(Y, E)$ donné  par 
  \begin{gather*}p_{1,!}(\mc F_{1}) \xrightarrow{\on{adj}} p_{1,!} a_{1,*} a_{1}^{*}(\mc F_{1})
  = p_{1,!} a_{1,!} a_{1}^{*}(\mc F_{1}) \\ 
  \xrightarrow{u} p_{1,!} a_{1,!} a_{2}^{!}(\mc F_{2})  =
 p_{2,!} a_{2,!} a_{2}^{!}(\mc F_{2}) \xrightarrow{\on{adj}} p_{2,!}(\mc F_{2})\end{gather*}
  où on a utilisé, dans la deuxième étape, l'hypothèse que  $a_{1}$ est propre. 

 Si $a_{1}$ et $a'_{2}$ sont propres et $p_{i}:X_{i}\to Y$ (pour $i=1,2,3$) sont des morphismes  tels que $p_{1}\circ a_{1}=p_{2}\circ a_{2}$ et $p_{2}\circ a'_{2}=p_{3}\circ a'_{3}$ 
   (et $X_{1},X_{2},X_{3}$ et $Y$ sont des champs de Deligne-Mumford) 
   on a  
   $   H(\wt u)=H(u')\circ H(u): p_{1,!}(\mc F_{1}) \to p_{3,!}(\mc F_{3})$. 
\end{rappel}

\begin{defi}
Un diagramme de correspondances  est un diagramme commutatif 

 $$   \xymatrix{
    &  A \ar[dl]_{a_{1}} \ar[d]_{f} \ar[dr]^{a_{2}} &           \\
  X_{1}  \ar[d]_{f_{1}}  &  B  \ar[dl]^{b_{1}} \ar[dr]_{b_{2}}&   X_{2}\ar[d]_{f_{2}}       \\
  Y_{1}      &     & Y_{2}       }$$     

Il est dit cartésien si le carré de droite $ABY_{2}X_{2}$ est cartésien. 

On suppose maintenant le diagramme cartésien. 
Soient  $\mc G_{1}\in D_{c}^{b}(Y_{1}, E)$,  $\mc G_{2}\in D_{c}^{b}(Y_{2}, E)$  et une correspondance cohomologique  
de $(Y_{1},\mc G_{1})$ vers  $(Y_{2},\mc G_{2})$ donnée par un 
morphisme 
  $v:(b_{2})_{!}(b_{1})^{*}(\mc G_{1})\to \mc G_{2}$

On pose $\mc F_{1}=(f_{1})^{*}(\mc G_{1})\in D_{c}^{b}(X_{1}, E)$ et $\mc F_{2}=(f_{2})^{*}(\mc G_{2})\in D_{c}^{b}(X_{2}, E)$. 

On définit  l'image inverse $u$ de $v$ comme la correspondance cohomologique 
de $(X_{1},\mc F_{1})$ vers  $(X_{2},\mc F_{2})$ donnée par le  morphisme 
  $u:(a_{2})_{!}(a_{1})^{*}(\mc F_{1})\to \mc F_{2}$
  égal à la composée 
   \begin{gather*} (a_{2})_{!}(a_{1})^{*}(\mc F_{1})=(a_{2})_{!}( a_{1})^{*}(f_{1})^{*}(\mc G_{1})=
  (a_{2})_{!} f^{*}(b_{1})^{*}(\mc G_{1}) \\ 
  \simeq (f_{2})^{*}(b_{2})_{!}(b_{1})^{*}(\mc G_{1})
  \xrightarrow{v}   (f_{2})^{*}(\mc G_{2})=\mc F_{2}. 
   \end{gather*}
  où l'isomorphisme au milieu vient de  l'isomorphisme du changement de base propre  $(a_{2})_{!} f^{*}=(f_{2})^{*} (b_{2})_{!}$. 
  \end{defi}

\begin{lem}
La composée de deux diagrammes cartésiens  de  correspondances    
$$   \xymatrix{
    &  A \ar[dl]_{a_{1}} \ar[d]_{f} \ar[dr]^{a_{2}} &          & A' \ar[d]_{f'} \ar[dl]_{a'_{2}} \ar[dr]^{a'_{3}} &
     \\
  X_{1}  \ar[d]_{f_{1}}  &  B  \ar[dl]^{b_{1}} \ar[dr]_{b_{2}}&   X_{2}\ar[d]_{f_{2}}  & B' \ar[dl]^{b'_{2}} \ar[dr]_{b'_{3}} & X_{3}\ar[d]_{f_{3}}
     \\
  Y_{1}      &     & Y_{2}  
  &&  Y_{3}
    }$$     
est un diagramme cartésien de  correspondances. 

De plus, étant donnés 
\begin{itemize}
\item $\mc G_{1}\in D_{c}^{b}(Y_{1}, E)$,  $\mc G_{2}\in D_{c}^{b}(Y_{2}, E)$, $\mc G_{3}\in D_{c}^{b}(Y_{3}, E)$, 
\item     une correspondance cohomologique  
de $(Y_{1},\mc G_{1})$ vers  $(Y_{2},\mc G_{2})$ donnée par un 
morphisme 
  $v:(b_{2})_{!}(b_{1})^{*}(\mc G_{1})\to \mc G_{2}$, 
  \item   une 
correspondance cohomologique  
de $(Y_{2},\mc G_{2})$ vers  $(Y_{3},\mc G_{3})$ donnée par un 
morphisme 
  $v':(b_{3})_{!}(b_{2})^{*}(\mc G_{2})\to \mc G_{3}$, 
  \end{itemize}
  l'image inverse de la composée $v'\circ v$ est égale à la composée des images inverses, comme correspondance cohomologique de 
  $(X_{1},(f_{1})^{*}(\mc G_{1}))$ vers  $(X_{3},(f_{3})^{*}(\mc G_{3}))$.
  \end{lem}
\dem 
On pose $\wt A=A\times_{X_{2}}A'$ et $\wt B=B\times_{Y_{2}}B'$. Alors
$ \wt A=B\times_{Y_{2}}A'=B\times_{Y_{2}} B'\times_{Y_{3}}X_{3}=\wt B\times_{Y_{3}}X_{3}$, de sorte que le diagramme des  correspondances composées est cartésien. Pour montrer que l'image inverse de $v'\circ v$ est égale à la composée des images inverses, la difficulté réside dans la partie centrale et on se ramène facilement au cas où $a_{1}, b_{1}, a_{3}', b_{3}'$ sont l'identité. 
Dans le diagramme commutatif 
\begin{gather} \nonumber  \xymatrix{
 &&     \wt A   \ar[ddll]  \ar[drr]|{\wt {a_{2}}}  \ar[d]_{\wt f} \ar[dll]|{\wt  {a'_{2}}}   &&         \\
     A  \ar[d]_{f} \ar[drr]|<<<<{a_{2}}  &&         \wt B   \ar[drr]|<<<<<{\wt {b_{2}}}   \ar[dll]|<<<<<{\wt {b'_{2}}}      &&  A' \ar[ddll] \ar[d]_{f'} \ar[dll]|<<<<<{a'_{2}}       \\
     B    \ar[drr]|{b_{2}} &&    X_{2}\ar[d]_{f_{2}}   &&  B' \ar[dll]|{b'_{2}}         \\
         &&  Y_{2}  
   &&  
    }\end{gather} 
     les carrés $AX_{2}A'\wt A$  et $BY_{2}B'\wt B$ sont cartésiens, ainsi que les carrés  
     $ABY_{2}X_{2}$ et $\wt A\wt B B'A'$, 
     et les deux isomorphismes de foncteurs
     \begin{gather} \label{isom1} (b'_{2}f')^{*} (b_{2})_{!}= (f')^{*}(b'_{2})^{*}(b_{2})_{!}
     \simeq (f')^{*}(\wt {b_{2}})_{!}(\wt  {b'_{2}})^{*}\simeq 
      (\wt {a_{2}})_{!}  (\wt f)^{*}(\wt  {b'_{2}})^{*}= (\wt {a_{2}})_{!}  (\wt  {b'_{2}}\wt f)^{*}
     \end{gather} et 
     \begin{gather}\label{isom2}(f_{2}a_{2}')^{*} (b_{2})_{!}= (a_{2}')^{*}(f_{2})^{*}(b_{2})_{!}
     \simeq (a_{2}')^{*}  (a_{2})_{!}f^{*}\simeq 
      (\wt {a_{2}})_{!} (\wt {a_{2}'})^{*}f^{*}= (\wt a_{2})_{!}  (f\wt {a_{2}'})^{*}
     \end{gather}
     sont identiques car égaux à l'isomorphisme de changement de base propre dans le carré cartésien $\wt A BY_{2}A'$. On en déduit le résultat car \eqref{isom1} donne l'image inverse de $v'\circ v$  et \eqref{isom2}  donne la composée des images inverses. 
 \cqfd


\section{Fin de la démonstration de la \propref{prop-mfz}}\label{rem-preuve-corresp}

  Comme précédemment on note $K$ un corps local de caractéristique $p$ et de corps résiduel fini 
  $\mf k$.   On fixe un point géométrique algébrique $\ov v$ au-dessus de $v=\on{Spec} \mf k$.  
   
   Soit $I$ un ensemble fini.  
 Soit $(\gamma_{i})_{i\in I}$ un 
  $I$-uplet  d'éléments de   $\on{Weil}(\ov K/K)$, dont l''image dans $ \Z^{I}$ est notée  $(d_{i})_{i\in I}$.  
      
 On pose $\Lambda=  
   \mc O_{E}/\lambda_{E}^{s}\mc O_{E}$. 
   Soit $f$ une fonction sur $\wh G\backslash \wh G ^{I}/\wh G$ définie sur $\Lambda$. 
      
 On va construire, après certains choix, un élément \begin{gather}\label{elem-construit}\mf z_{m,s,I,f,(\gamma_{i})_{i\in I}}\in  
C_{c}(U_{m}\backslash G(K)/U_{m},\Lambda)\end{gather} 
et on montrera à l'aide d'arguments globaux que le résultat ne dépend pas de ces choix. 

   Les choix sont les suivants: 
  \begin{itemize}
  \item  une représentation   de $\wh G ^{I}$ sur un $\Lambda$-module libre de type fini $W$,  
   ainsi que $x\in W$ et $\xi:W\to \Lambda$ invariants par l'action diagonale de $\wh G$ tels que $f$ soit le coefficient de matrice 
   $(g_{i})_{i\in I}\mapsto \s{\xi, (g_{i})_{i\in I}.x}$, 
  \item  un ordre des coordonnées, c'est-à-dire une bijection  $\sigma:I\isom \{1,...,l\}$ et on prend comme partition $(I_{1},...,I_{l})=(\{1\}, ..., \{l\})$
    \item un relèvement  $\delta \in \on{FWeil}(\ov{K^{1,...,l}}/K^{1,...,l})$
  de $(\gamma_{i})_{i\in I}$ 
  (grâce au \lemref{lem-surj-Weil-loc} un tel relèvement existe toujours), 
  \item $(m_{i})_{i\in I}\in \N^{I}$ 
 tel que 
 $(m_{i}+d_{i})_{i\in I}\in \N^{I}$ et 
 assez croissant  pour  que l'énoncé de la \propref{prop-boite-noire}
  soit satisfait.
\end{itemize}

   La \propref{prop-boite-noire}  fournit un morphisme $T_{\delta}:R\Psi_{(m_{i})_{i\in I}} \mc S\mc L\to 
R\Psi_{(m_{i}+d_{i})_{i\in I}} \mc S\mc L$.

    \begin{rem}
 La \remref{rem-complement-Tdelta},  qui renforce l'énoncé de la \propref{prop-boite-noire}, impliquerait que $T_{\delta}$ est indépendant du choix de $(m_{i})_{i\in I}$ et compatible à la composition de $\delta$.  Nous n'en avons pas besoin pour la suite.   
\end{rem}

   La construction   de l'élément \eqref{elem-construit}    consiste  en  une composition de correspondances cohomologiques entre les faisceaux de cycles proches 
  du type   \begin{gather}\label{cycles-proches-composition}R\Psi_{(l_{i})_{i\in I}} \mc S\mc L_{ I,W,(n_{i})_{i\in I},r} ^{nv, \Lambda,(I_{1},...,I_{k})}\in
  D_{c}^{b}\big(\restr{\Chr_{Y,I,W,(n_{i})_{i\in I},r}^{nv,(I_{1},...,I_{k})}}{\Delta(\ov v)}, \Lambda\big). \end{gather}
   Les niveaux   seront de la forme 
$ nv$ avec  
$n$ assez grand (plus grand que $m$ et  assez grand à chaque étape de la composition en fonction de l'étape suivante, si nécessaire). De même $r$ sera choisi  assez grand pour pouvoir être diminué de $1$ à chaque fois que l'on applique un morphisme de Frobenius partiel et il en ira de même avec les $ (n_{i})_{i\in I}$. 

Les points de départ et d'arrivée de cette suite de correspondances cohomologiques sont des champs de chtoucas restreints sans pattes 
munis d'un faisceau lisse $\mc L$ indiquant le niveau. Concrètement ces champs 
 sont   $\bullet/G(\mc O_{nv})$      muni du $\Lambda$-faisceau  lisse associé à la représentation régulière de $G(\mc O_{ mv})$
 (ce qui a un sens car $n\geq m$).

La suite de correspondances cohomologiques que l'on compose pour obtenir 
\eqref{elem-construit} est  la suivante: 
\begin{itemize}
\item le morphisme de création associé à $x\in W$, qui est supporté sur la diagonale et aboutit donc dans un faisceau  \eqref{cycles-proches-composition} avec 
$(l_{i})_{i\in I}=0$, 
\item des morphismes de Frobenius partiels  \eqref{corresp-coho-frob-partiels}     pour passer de 
$(l_{i})_{i\in I}=0$ à $(l_{i})_{i\in I}=(m_{i})_{i\in I}$ 
(en plus ces morphismes doivent être intercalés avec des correspondances cohomologiques \eqref{oubli-cycles-proches} et leurs inverses  pour remettre  les pattes dans l'ordre), 
  \item le morphisme $T_{\delta}$ de la \propref{prop-boite-noire}
  (un morphisme de faisceaux est un cas particulier de correspondance  cohomologique), 
 \item  des morphismes de Frobenius partiels pour passer de 
 $(l_{i})_{i\in I}=(m_{i}+d_{i})_{i\in I}$ à $(l_{i})_{i\in I}=0$,  intercalés avec des  correspondances cohomologiques pour remettre les pattes dans l'ordre, 
 \item un morphisme d'annihilation associé à $\xi$. 
  \end{itemize}
  
  Cette composée n'est pas en général calculable, comme  c'était le cas dans le cas particulier du lemme 6.11 de \cite{coh} (que l'on peut considérer comme un calcul de cette composée en tant que correspondance cohomogique dans un cas très simple). 
  
  
  Comme les espaces de départ et d'arrivée sont des champs de chtoucas restreints sans pattes, la composée totale de ces correspondances cohomologiques est une correspondance cohomologique donnée par 
  \begin{itemize}
  \item la correspondance géométrique constituée par 
   $G(\mc O_{nv}) \backslash Z(\Fq)$ muni des morphismes   $\on{pr_{1}}$ et $\on{pr_{2}}$  vers $\bullet/G(\mc O_{nv})$ et 
  $\bullet/G(\mc O_{K})$, où $Z$ est une réunion finie de strates fermées de la grassmannienne affine en $v$ et $Z(\Fq)\subset G(K)/G(\mc O_{K})$ est muni d'une action de $G(\mc O_{nv})$ pour $n$ assez grand, 
  \item une fonction sur $G(\mc O_{nv})\backslash Z(\Fq)$, à valeurs dans $\Hom(\on{pr_{1}}^{*}(\mc L), \on{pr_{2}}^{*}(\mc L))$  
  où $\mc L$ est le faisceau sur $\bullet/G(\mc O_{nv})$ ou 
  $\bullet/G(\mc O_{K})$  associé à la représentation régulière de $G(\mc O_{mv})$ à coefficients dans $\Lambda$ 
  (en effet on rappelle que $m\leq n$ et on note que $\on{pr_{2}}^{!}(\mc L)=\on{pr_{2}}^{*}(\mc L)$). 
  \end{itemize}
  Cela fournit donc  une 
 fonction à support fini  sur $G(\mc O_{nv})\backslash Z(\Fq)$, à valeurs dans $\Hom(\on{pr_{1}}^{*}(\mc L), \on{pr_{2}}^{*}(\mc L))$, c'est-à-dire, en prolongeant par $0$ de $$ \bullet/U_{m}
 \times_{\bullet/G(\mc O_{nv})} (G(\mc O_{nv}) \backslash Z(\Fq))  \times_{\bullet/G(\mc O_{K}) }  \bullet/U_{m}$$ \`a $U_{m}\backslash G(K)/U_{m}$), 
  une fonction 
    \begin{gather}\label{elem-construit2}
           \mf z_{m,s,I,f,(\gamma_{i})_{i\in I}}\in  
C_{c}(U_{m}\backslash G(K)/U_{m},\mc O_{E}/\lambda_{E}^{s}\mc O_{E}). \end{gather}

  Il est clair que $ \mf z_{m,s,I,f,(\gamma_{i})_{i\in I}}$  ne dépend que des données sur $Y$. Autrement dit  il est local au sens hensélien en $v$, et donc il est local  au sens du complété en $v$, d'apr\`es la remarque \ref{rem-hensel-complete}.

  \begin{lem} Dans le cadre global-local de \eqref{mor-Y-X}
 et \eqref{cond-X-Y-v}, 
l'opérateur d'excursion global $S_{I,f,(\gamma_{i})_{i\in I}}$ agissant sur les formes automorphes à coefficients dans $\mc O_{E}/\lambda_{E}^{s}\mc O_{E}$  est égal à la multiplication à droite  par $  \mf z_{m,s,I,f,(\gamma_{i})_{i\in I}}$. 
\end{lem} 
\dem 
En effet $S_{I,f,(\gamma_{i})_{i\in I}}$ est le morphisme entre groupes de cohomologie à support compact associé grâce au rappel~\ref{rappel-coho} à la composée de  correspondances  entre champs de chtoucas globaux qui sont les images inverses par les morphismes de restriction des correspondances considérées ci-dessus entre champs de chtoucas restreints. On donne ci-dessous les diagrammes cartésiens  grâce auxquels 
les correspondances  entre champs de chtoucas globaux 
sont les images inverses des correspondances  entre champs de chtoucas  restreints considérées ci-dessus. 

La correspondance donnant l'action du morphisme de Frobenius partiel (comme dans 
\eqref{corresp-coho-frob-partiels})   est  

\begin{gather*} 
  \xymatrix{
    \restr{ \Cht_{N^{v}, I,W} ^{(I_{2},...,I_{k},I_{1})}}{\Delta(\ov v)}  \ar[d]_{  \mc R_{ I,W,(n_{i})_{i\in I},r-1} ^{nv,(I_{2},...,I_{k},I_{1})}}       & 
    \restr{  \Cht_{N^{v}, I,W} ^{(I_{1},...,I_{k})}}{\Delta(\ov v)}
     \ar[d]_{  \mc R_{ I,W,(n_{i})_{i\in I},r} ^{nv,(I_{1},...,I_{k})}} \ar[r]^{\Id}
       \ar[l]_{\on{Fr}_{Y,I_{1}}}
  &  \restr{ \Cht_{N^{v}, I,W} ^{(I_{1},...,I_{k})}}{\Delta(\ov v)}  
  \ar[d]_{  \mc R_{ I,W,(n_{i})_{i\in I},r} ^{nv,(I_{1},...,I_{k})}} 
     \\
  \restr{ \Chr_{Y,I,W,(n_{i})_{i\in I},r-1} ^{nv,(I_{2},...,I_{k},I_{1})}}{\Delta(\ov v)}
       &\restr{\Chr_{Y,I,W,(n_{i})_{i\in I},r} ^{nv,(I_{1},...,I_{k})} }{\Delta(\ov v)}
       \ar[l]_{\on{Fr\mc R}_{Y,I_{1}}}  \ar[r]^{\Id}& \restr{ \Chr_{Y,I,W,(n_{i})_{i\in I},r} ^{nv,(I_{1},...,I_{k})}}{\Delta(\ov v)}
    }\end{gather*}

Les correspondances donnant l'oubli de la partition (comme dans \eqref{oubli-cycles-proches}), dans un sens ou dans l'autre, sont 
\begin{gather*} 
  \xymatrix{
   \restr{  \Cht_{N^{v}, I,W} ^{(I)}}{\Delta(\ov v)} \ar[d]_{  \mc R_{ I,W,(n_{i})_{i\in I},r} ^{nv,(I)}}       & 
    \restr{  \Cht_{N^{v}, I,W} ^{(I_{1},...,I_{k})}}{\Delta(\ov v)}
     \ar[d]_{  \mc R_{ I,W,(n_{i})_{i\in I},r} ^{nv,(I_{1},...,I_{k})}} \ar[r]^{\Id}
       \ar[l]
  &  \restr{ \Cht_{N^{v}, I,W} ^{(I_{1},...,I_{k})}}{\Delta(\ov v)}  
  \ar[d]_{  \mc R_{ I,W,(n_{i})_{i\in I},r} ^{nv,(I_{1},...,I_{k})}} 
     \\
  \restr{ \Chr_{Y,I,W,(n_{i})_{i\in I},r} ^{nv,(I)}}{\Delta(\ov v)}
       &\restr{\Chr_{Y,I,W,(n_{i})_{i\in I},r} ^{nv,(I_{1},...,I_{k})} }{\Delta(\ov v)}
       \ar[l] \ar[r]^{\Id}& \restr{ \Chr_{Y,I,W,(n_{i})_{i\in I},r} ^{nv,(I_{1},...,I_{k})}}{\Delta(\ov v)}
    }\end{gather*}
et
\begin{gather*} 
  \xymatrix{
   \restr{ \Cht_{N^{v}, I,W} ^{(I_{1},...,I_{k})}}{\Delta(\ov v)}  
  \ar[d]_{  \mc R_{ I,W,(n_{i})_{i\in I},r} ^{nv,(I_{1},...,I_{k})}} 
     & 
    \restr{  \Cht_{N^{v}, I,W} ^{(I_{1},...,I_{k})}}{\Delta(\ov v)}
     \ar[d]_{  \mc R_{ I,W,(n_{i})_{i\in I},r} ^{nv,(I_{1},...,I_{k})}} \ar[r] 
       \ar[l]_{\Id} 
  &  \restr{  \Cht_{N^{v}, I,W} ^{(I)}}{\Delta(\ov v)} \ar[d]_{  \mc R_{ I,W,(n_{i})_{i\in I},r} ^{nv,(I)}}  
       \\   \restr{ \Chr_{Y,I,W,(n_{i})_{i\in I},r} ^{nv,(I_{1},...,I_{k})}}{\Delta(\ov v)}
       &\restr{\Chr_{Y,I,W,(n_{i})_{i\in I},r} ^{nv,(I_{1},...,I_{k})} }{\Delta(\ov v)}
       \ar[l]_{\Id} \ar[r]&  \restr{ \Chr_{Y,I,W,(n_{i})_{i\in I},r} ^{nv,(I)}}{\Delta(\ov v)}
          }\end{gather*}
De plus $ \restr{  \Cht_{N^{v}, I,W} ^{(I)}}{\Delta(\ov v)}$ et $ \restr{ \Chr_{Y,I,W,(n_{i})_{i\in I},r} ^{nv,(I)}}{\Delta(\ov v)}$ de dépendent pas de $(I)$ (puisque toutes les pattes sont égales au-dessus de  $\Delta(\ov v)$), ce qui exprime les isomorphismes de fusion. 
Enfin les morphismes de création et d'annihilation et les $T_{\delta}$ sont donnés par des morphismes de faisceaux $\restr{\Chr_{Y,I,W,(n_{i})_{i\in I},r} ^{nv,(I_{1},...,I_{k})} }{\Delta(\ov v)}$, et leurs images inverses par $ \mc R_{ I,W,(n_{i})_{i\in I},r} ^{nv,(I_{1},...,I_{k})}$ sont des morphismes de faisceaux sur  
$  \Cht_{N^{v}, I,W} ^{(I_{1},...,I_{k})}$. \cqfd

Comme $N^{v}$ est arbitraire, le \lemref{lem-comme-poincare} montre que cet élement \eqref{elem-construit2} est déterminé de manière unique par les opérateurs d'excursion globaux. 
Donc il ne dépend pas du choix de l'ordre des pattes $\sigma$,  ni du relèvement $\delta$ de $(\gamma_{i})_{i\in I}$, ni du choix de  $(m_{i})_{i\in I}$ dans la construction de $T_{\delta}$. 
Pour la même raison il dépend seulement de la fonction $f:(g_{i})_{i\in I}\mapsto
\s{\xi, (g_{i})_{i\in I} \cdot x}$ et non de $W,x,\xi$. 
C'est pourquoi on l'a noté $\mf z_{m,s,I,f,(\gamma_{i})_{i\in I}}$. 
La construction montre que pour $m,s,I,f$ et $(d_{i})_{i\in I}$ fixés le support de 
$\mf z_{m,s,I,f,(\gamma_{i})_{i\in I}}$est borné  indépendamment de $(\gamma_{i})_{i\in I}$. La compatibilité avec le global montre alors la continuité en $(\gamma_{i})_{i\in I}$. 
Finalement $\mf z_{m,s,I,f,(\gamma_{i})_{i\in I}}$ est exactement l'élément \eqref{elem-construit} recherché, et on a  terminé la démonstration de la \propref{prop-mfz}. 

  
  \begin{rem}\label{rem-support-para5}
  On voit que 
  \begin{itemize}
  \item le support de \eqref{elem-construit} dépend du nombre d'itérations dans la composition de correspondances cohomologiques, donc des entiers $(m_{i})_{i\in I}$ (en effet pour itérer les morphismes de Frobenius partiels, on a besoin des  correspondances donnant l'action des morphismes de Frobenius partiels, ainsi que des correspondances donnant l'oubli de la partition, dans un sens ou dans l'autre, que nous avons écrites dans les trois grands diagrammes ci-dessus),  
    \item les entiers $(m_{i})_{i\in I}$ dépendent de la modification $\wt {Y^I}$ de $Y^{I}$   fournie par \cite{orgogozo}, 
  \item cette modification $\wt {Y^I}$   dépend {\it a priori} de l'entier   $m$ (qui indique le niveau) et de l'entier  $s$ (la puissance de $\ell$ dans les coefficients). 
  \end{itemize}
  
  Il est évident que le support de  \eqref{elem-construit} dépend de $m$, puisqu'à la limite sur $m$ on obtient un élément du centre de Bernstein, donc une distribution invariante par conjugaison.    La seule fa\c con d'expliquer que le support de \eqref{elem-construit}, et donc forcément le nombre d'itérations des morphismes de Frobenius partiels,  tendent vers l'infini quand $m$ tend vers l'infini est que la modification $\wt{Y^{I}}$ de \cite{orgogozo} devienne de plus en plus profonde quand $m$ tend vers l'infini, ce qui n'a rien de surprenant.

  En revanche 
      la modification $\wt {Y^I}$ ne  dépend pas   de l'exposant $s$ de $\ell$ d'apr\`es  
  le corollaire 2.9 de \cite{ito}. Comme ce corollaire concerne le faisceau constant, il faut quelques arguments suppl\'ementaires 
  que nous ne d\'etaillerons pas : des r\'esolutions de Demazure pour \'ecrire les faisceaux de Mirkovic-Vilonen comme facteurs directs d'images directes propres de faisceaux constants, une normalisation  du champ des chtoucas restreints  $\Chr_{Y,I,W,(n_{i})_{i\in I},r}^{nv,(I_{1},...,I_{k})}$
  dans le  rev\^etement \'etale de groupe $G(\mc O_{nv})$ 
  de $\restr{\Chr_{Y,I,W,(n_{i})_{i\in I},r}^{nv,(I_{1},...,I_{k})}}{(Y\sm v)^{I}}$  donn\'e par 
\eqref{mor-GON}
  (apr\`es avoir choisi une pr\'esentation des champs de chtoucas restreints par des sch\'emas pour donner un sens \`a la normalisation). 
  La notion de facteur direct pour des $\Z_\ell$-faisceaux n'est d\'efinie que modulo torsion mais comme des modifications existent pour ces faisceaux de torsion cela ne pose pas de probl\`eme.  
        Le corollaire 2.9 de \cite{ito} implique donc que l'on peut avoir la modification $\wt {Y^I}$ de $Y^{I}$ ind\'ependante de l'entier  $s$. 
        On en d\'eduit que  le support de 
   \eqref{elem-construit} ne dépend pas de $s$. 
      \end{rem}

  \begin{rem}
  On peut se demander pourquoi la modification $\wt {Y^I}$ de $Y^{I}$  n'est pas invariante par les morphismes de Frobenius partiels sur $Y^{I}$ (ce qui la rendrait triviale), puisqu'après tout
    les morphismes de Frobenius partiels sur les champs de chtoucas restreints sont des homéomorphismes locaux 
  à des morphismes lisses près (les morphismes d'oubli du niveau \eqref{oubli-chtR-ChtR}). Il semble que la raison soit la suivante: pour itérer les morphismes de Frobenius partiels on doit les intercaler avec des changements de l'ordre des pattes, c'est-a-dire des morphismes d'oubli des modifications intermédiaires
  \eqref{oubli-cycles-proches} et leurs inverses (vus comme correspondances cohomologiques). 
  Or ces inverses peuvent obliger à  raffiner  la modification de $Y^{I}$.   
      \end{rem}

\section{Facteurs $\gamma$}\label{facteurs-gamma}

En combinant le résultat principal de cet article et des   arguments de \cite{lomeli} et \cite{hl13a}, on construit une bonne théorie des facteurs $\gamma$ locaux dans le cas  des corps de fonctions.  Le point essentiel est que, contrairement aux facteurs $L$ et $\epsilon$ locaux, les facteurs $\gamma$ locaux ne dépendent que du caractère du centre de Bernstein, ou du paramètre de Langlands local à semi-simplification près. 

\subsection{Axiomes pour un  système de facteurs $\gamma$}

On considère une   famille de facteurs  $\gamma(s,\pi \times \chi,\psi) \in \mathbb{C}(T)$, lorsque $K$ est un corps local de caractéristique $p$ contenant $\Fq$, $G$ est un groupe réductif déployé sur $K$, $\rho$ est une représentation finie de dimension $N$ de $\wh G$ sur $ \Qlbar$, $\pi$ est une représentation lisse  irréductible  à coefficients dans  $\Qlbar $ de $G(K)$, $\chi:K^{*}\to \Qlbar^{*}$ est un caractère continu, 
et $\psi:K\to \Qlbar$ un caractère additif non dégénéré. 
Les propriétés que doit satisfaire une telle famille sont les suivantes. 

\begin{itemize}
   \item[(i)] (Naturalité et isomorphisme) $\gamma(s,\rho, \pi, \chi,\psi)$ dépend  de $(K,\pi,\chi,\psi)$ à isomorphisme près. 
      \item[(ii)] (Somme directe) Si $\rho=\rho_{1}\oplus \rho_{2}$, alors    $$\gamma(s,\rho, \pi, \chi,\psi) =\gamma(s,\rho_{1}, \pi,  \chi,\psi) \gamma(s,\rho_{2}, \pi, \chi,\psi). $$
       \item[(iii)]  (Compatibilité avec la théorie du corps de classe)
      Si $G$ est un tore et $\rho$ est un caractère c'est le facteur $\gamma$ abélien de Tate.   
       \item[(iv)]  (Compatibilité à l'induction parabolique) 
      Soit $M$ un sous-groupe de Levi de $M$,  $P$ un parabolique de Levi $M$, $\pi_{M}$ une représentation irréductible de $M$, $\rho_{M}$ la représentation de dimension $N$ de $\wh M$ obtenue en composant $\rho$ avec le morphisme
      $\wh M\to \wh G$,   et $\pi$ un sous-quotient irréductible de l'induite parabolique  $\on{Ind}_{P(K)}^{G(K)}\pi_{M}$. 
      Alors $\gamma(s,\rho, \pi,  \chi,\psi) =\gamma(s,\rho_{M}, \pi_{M},  \chi,\psi)$.

   \item[(v)] (Dépendance en $\psi$). Soit  $a \in K^*$, et soit $\psi^a$ le caractère de $K$ donné par $\psi^a(x) = \psi(ax)$. Alors 
      \begin{equation*}
    \gamma(s,\rho, \pi,  \chi,\psi^a) = 
    \omega_\pi(a)  \chi(a)^N \left| a \right|_F^{N(s - \frac{1}{2})}  \gamma(s,\rho, \pi,  \chi,\psi)
   \end {equation*}
  où   $\omega_\pi:K^{*}\to \Qlbar^{*}$
  est la restriction du caractère central de $\pi$ à $\mb G_{m}\to G$ associé à
  $\wh G\xrightarrow{\rho} GL_{N}\xrightarrow{\det}\mb G_{m}$.

   \item[(vi)] (Stabilité). Si   $\pi_1$ et  $\pi_2$ ont le même caractère central et si 
    $\chi$ est assez ramifié en fonction de   $\pi_1$ et  $\pi_2$, alors 
   \begin{equation*}
     \gamma(s,\rho, \pi_{1},  \chi,\psi)=   \gamma(s,\rho, \pi_{2},  \chi,\psi).
   \end{equation*}
   
   \item[(vii)] (Equation fonctionnelle)
   Soit $X$ une courbe sur $\Fq$, $\pi=\bigotimes_{v\in |X|}\pi_{v}$ une représentation automorphe cuspidale de $G(\mb A)$ non ramifiée   en dehors d'un ensemble fini $S$ de places, alors  \begin{equation*}
      L^S(s,\rho, \pi,  \chi) = \prod_{v \in S}
       \gamma(s,\rho, \pi_v,  \chi_{v},\psi_{v})
        L^S(1-s, \rho, \wt \pi,  \wt \chi). 
   \end{equation*}   
   \end{itemize}

\subsection{Unicité}
Par (iv) on peut se limiter à montrer l'unicité pour les supercuspidales. 
On fixe une courbe $X$ sur un corps $\Fq$ inclus dans $\mf k$ et une place 
$v$ telle que $k(v)=\mf k$. On fixe un isomorphisme $K=F_{v}$. 
Par le même argument de séries de Poincaré que dans la preuve du \lemref{lem-comme-poincare} on trouve une forme automorphe cuspidale $\pi=\bigotimes \pi_{w}$ telle que $\pi_{v}$ soit la représentation supercuspidale de $G(K)$ pour laquelle on veut montrer l'unicité du facteur $\gamma$. On fixe un caractère global  $\chi$ non ramifié en $v$ et très ramifié en toutes les  places autres que $v$ où $\pi$ est ramifié. 
Alors dans l'équation fonctionnelle tous les facteurs $\gamma$ sont connus sauf celui en $v$, qui est donc déterminé de manière unique.  

\subsection{Existence}
On pose $$\gamma(s,\rho, \pi, \chi,\psi)=\gamma(s,\rho\circ \sigma_{ \pi}, \chi,\psi)$$
où $\sigma_{\pi}$ est le paramètre de Langlands local à semi-simplification près associé à $\pi$ par le \thmref{thm-intro}, et le facteur $\gamma$ de droite est le facteur local de Godement-Jacquet habituel pour $GL_{N}$ (en notant $N$ la dimension de $\rho$). Il apparaît ici du côté galoisien mais 
comme  la correspondance locale pour $GL_{N}$ est connue par Laumon-Rapoport-Stuhler \cite{laumon-rapoport-stuhler},   on peut aussi le voir du côté automorphe et il a toutes les propriétés souhaitées (en fait l'existence des facteurs locaux galoisiens pour $GL_{N}$ est due à Artin, Dwork,  Langlands \cite{lan1} et Deligne \cite{del1} et l'équation fonctionnelle est due à Laumon \cite{laumon-eq-fct}). On vérifie alors les propriétés (i)-(vii).  

L'énoncé d'unicité implique que notre construction étend celles des articles de  Henniart-Lomel\'{\i},  Lomel\'{\i},  Gan-Lomel\'{\i}, Ganapathy-Lomel\'{\i}  
\cite{hl11, hl13b, hl13a, gana-lomeli, lomeli, lomeli-AIF, 
lomeli15, gan-lomeli, hl17} mais ces articles montrent d'autres propriétés des facteurs $\gamma$ et aussi des facteurs $L$ et $\epsilon$ associés, que nous ne savons pas obtenir par notre méthode. 

\begin{rem} Luis Lomel\'{\i} nous a fait remarquer que la condition de stabilité (vi) n'est en fait pas nécessaire pour montrer l'unicité d'une théorie de facteurs $\gamma$ vérifiant les autres conditions. L'argument consiste à remplacer l'énoncé naïf de globalisation utilisé ci-dessus (le \lemref{lem-comme-poincare}) par un énoncé plus fin, à savoir le lemme 3.1 de \cite{lomeli15}. 
\end{rem}

\section{Renforcement des résultats avec action sur la cohomologie}

Dans  \cite{cong-coeff-O} (pour le cas déployé, et  \cite{cong-zorro} pour tous les groupes par une approche un peu différente), Cong Xue a défini une action des opérateurs d'excursion $S_{I,f,(\gamma_{i})_{i\in I}}$ sur la cohomologie à support compact et à coefficients finis des champs de chtoucas. 
On montre le renforcement suivant de la \propref{prop-mfz}. 
 On se donne 
\begin{itemize}
\item un ensemble fini $I$, 
\item une fonction $f$ sur le quotient grossier $\wh G \backslash \wh G ^{I}/\wh G$ définie sur $\mc O_{E}$, 
\item un $I$-uplet $(\gamma_{i})_{i\in I}$ d'éléments de  $\on{Weil}(\ov K/K)$, 
\item un entier $m$ assez grand pour que 
$ U_{m} $
soit d'ordre premier à $\ell$, 
\end{itemize}
et, pour tout entier $s$ on note $$\mf z_{m,s,I,f,(\gamma_{i})_{i\in I}}\in  
C_{c}(U_{m}\backslash G(K)/U_{m},\mc O_{E}/\lambda_{E}^{s}\mc O_{E})$$ l'élément construit précédemment (et dont on a déjà montré l'unicité). En passant à la limite sur $s$ on obtient 
$$\mf z_{m,I,f,(\gamma_{i})_{i\in I}}\in  
C_{c}(U_{m}\backslash G(K)/U_{m},\mc O_{E})$$ 
 où $C_{c}$ signifie que la limite à l'infini est $0$ (on s'attend à ce que cet élément soit en fait à support compact). 
 
\begin{prop}\label{prop-mfz-coho} 
Pour $\Fq,X,  \Xi$ comme précédemment, pour toute place $v$ de $X$,  
     en prenant $K=F_{v}$, et pour tout plongement $\ov{ F_X}\subset \ov K$
     (d'où $\on{Weil}(\ov K/K)\subset \on{Gal}(\ov{ F_X}/F_X)$) et pour tout sous-schéma  fini  $N^{v}$ de $X\sm v$,  pour tout ensemble fini $J$ et toute représentation $V$ de $(\wh G)^{J}$ à coefficients dans $ \mc O_{E}$, 
         l'opérateur d'excursion $S_{I,f,(\gamma_{i})_{i\in I}}$   (construit grâce à    \cite{cong-coeff-O} pour le cas déployé, et  \cite{cong-zorro} pour tous les groupes par une approche un peu différente)      agit sur $   \varinjlim _{\mu} \restr{\mc H _{ N, J, V}^{0,\leq\mu, \mc O_{E}}}{\Delta(\ov\eta)}$ par convolution à droite par $\mf z_{m,I,f,(\gamma_{i})_{i\in I}}$. 
                \end{prop}


 \section{Cas des groupes non déployés} \label{cas-non-deploye}  
 
 Tous les énoncés de cet article restent vrais pour $G$ non déployé (après modification convenable de l'énoncé).  
 
 Soit $G$ un groupe réductif sur $K$. 
 Le  $L$-groupe ${}^{L }G$ est alors le   produit   semi-direct
   $\wh G\rtimes \on{Gal}(\wt K/K)$ où  $\wt K$ est l'extension finie galoisienne de $K$  telle  que  $\on{Gal}(\wt K/K)$ soit l'image  de $\on{Gal}(\ov K/K)$ dans le groupe des  automorphismes du diagramme de  Dynkin  de $G$. 
 On rappelle que ce  produit   semi-direct est pris  pour l'action de $\on{Gal}(\wt K/K)$ sur 
 $\wh G$ qui préserve un épinglage, voir \cite{borel-corvallis}. 

On choisit des modèles entiers parahoriques lisses, à fibres géométriquement connexes,  pour $G$ sur $\on{Spec} \mc O_{K}$, $X$ et $Y$, compatibles entre eux. 
Pour toutes les constructions de champs de chtoucas globaux et restreints
 on utilise ces modèles entiers.  
On a besoin de construire des champs de chtoucas même lorsque les pattes rencontrent le  point $v$ où le groupe $G$ n'est pas nécessairement réductif 
(alors que cela n'était pas nécessaire dans \cite{coh}). 
Mais cela n'est pas difficile, on renvoie par exemple à \cite{hartl}. 
D'après le théorème A de \cite{richarz2} le quotient $G(\mc K)/G(\mc O)$ est ind-propre, et plus généralement la grassmannienne affine de Beilinson-Drinfeld est ind-propre sur $Y^{I}$.  Cela permet notamment de justifier la propreté du morphisme d'oubli des modifications intermédiaires dans le \lemref{lem-petit}. 
A part cela on n'a besoin de rien savoir sur la fibre spéciale en $v$, seulement qu'elle existe et supporte les faisceaux de cycles proches. 
Tout le reste marche comme dans le cas déployé. Plus précisément on a les théorèmes suivants.

 Pour énoncer la compatibilité local-global on rappelle le cadre du paragraphe 12 
  de \cite{coh}. 
 Soit  $\Fq$ un corps fini,  $X$ une courbe projective lisse géométriquement irréductible  sur $\Fq$, $F_X$ son corps des fonctions et $\mc A$ ses adèles. Soit $G$ un groupe réductif sur $F_{X}$. 
 On fixe un sous-groupe discret cocompact $\Xi\subset Z_{G}(F_X)\backslash Z_{G}(\mb A)$ (où $Z_{G}$ est le centre de $G$).  
  On choisit un modèle réductif de $G$ sur un ouvert $U$ de $X$ et un modèle parahorique lisse de $G$, à fibres géométriquement connexes,  en tous les points de $X\sm U$. On fixe $N$ un sous-schéma fini de $X$.  On rappelle que  le théorème 12.3   \cite{coh}  fournit 
  une décomposition canonique de  
   $C_{c}(K_{N}\backslash G(\mb A)/K_{N},\Qlbar)$-modules 
 \begin{gather}\label{intro1-dec-canonique-loc}
 C_{c}^{\mr{cusp}}(\Bun_{G,N}(\Fq)/\Xi,\Qlbar)=\bigoplus_{\sigma}
 \mf H_{\sigma},\end{gather}
 où la somme directe dans le membre de droite est indexée par des paramètres de Langlands globaux au sens rappelé dans l'énoncé de ce  théorème. 
 
\begin{thm} \label{thm-intro-loc}
Il existe une application 
\begin{gather}\label{param-local-sigma-pi-loc}\pi\mapsto \sigma_{\pi}\end{gather}
\begin{itemize}
\item []
de l'ensemble des classes d'isomorphismes de représentations lisses admissibles et irréductibles de $G(K)$ définies sur une extension finie 
de $\Ql$ et entières, 
\item [] vers l'ensemble des classes d'isomorphismes de paramètres de Langlands locaux semi-simples, c'est-à-dire 
  des classes de  $\wh G(\Qlbar)$-conjugaison de  morphismes 
       $\sigma:\on{Gal}(\ov K/K)\to {}^{L }G(\Qlbar)$ 
       définis sur une extension finie de  $\Ql$, continus,    semi-simples et faisant commuter le diagramme
       \begin{gather}\label{diag-sigma-loc}
 \xymatrix{
\on{Gal}(\ov K/K) \ar[rr] ^{\sigma}
\ar[dr] 
&& {}^{L} G(\Qlbar) \ar[dl] 
 \\
& \on{Gal}(\wt K/K) }\end{gather}
\end{itemize}
 qui est déterminée de manière unique par les deux propriétés suivantes: 
\begin{itemize}
\item [] a)   $\sigma_{\pi}$ ne dépend que du caractère par lequel le centre de Bernstein agit sur $\pi$, et en dépend ``algébriquement'', 
\item [] b) cette application 
 est compatible avec la paramétrisation globale, au sens suivant. 
 \end{itemize}

  Soient  $\Fq,  X,  N,  G$ et $\Xi$ comme ci-dessus. 
  Alors pour toute représentation $\pi=\bigotimes \pi_{v}$ de $G(\mb A)$ telle que $\pi^{K_{N}}$ apparaisse dans $\mf H_{\sigma}$ (dans \eqref{intro1-dec-canonique-loc} ci-dessus), pour toute place $v$ de $X$ on a égalité entre 
  \begin{itemize}
  \item  le paramètre local $\sigma_{\pi_{v}}$ obtenu en appliquant 
  \eqref{param-local-sigma-pi-loc}  avec   $K$ égal au complété 
 $F_{v}$ de $F_X$ en $v$,   
  \item  le semi-simplifié de la restriction  de $\sigma$ à $\on{Gal}(\ov{F_{v}}/F_{v})$.  
 \end{itemize}
 
 De plus cette application $\pi\mapsto \sigma_{\pi}$ 
 s'étend de fa\c con unique en une application 
 \begin{itemize}
\item  []
de l'ensemble des classes d'isomorphismes de représentations lisses admissibles et irréductibles de $G(K)$ définies sur une extension finie 
de $\Ql$ (pas nécessairement entières)  
\item [] vers l'ensemble des classes d'isomorphismes de paramètres de Weil locaux semi-simples, c'est-à-dire 
  des classes de  $\wh G(\Qlbar)$-conjugaison de  morphismes 
       $\sigma:\on{Weil}(\ov K/K)\to {}^{L} G(\Qlbar)$ 
       définis sur une extension finie de  $\Ql$, continus,    semi-simples et faisant commuter un diagramme analogue à \eqref{diag-sigma-loc}, 
\end{itemize}
  vérifiant a) ci-dessus. 

Cette application est compatible 
  à l'induction parabolique au sens suivant. Si $P$ est un sous-groupe parabolique de $G$ sur $K$, de quotient de Levi $M$, et si $\tau $ est une représentation lisse, admissible et  irréductible de $M(K)$ définie sur une extension finie 
de $\Ql$, et $\pi$ est un sous-quotient de la représentation induite compacte $\on{Ind}_{P(K)}^{G(K)}\tau$ (avec la normalisation unitaire), 
alors $\sigma_{\pi}$ est conjugué à la composée 
$ \on{Weil}(\ov K/K)\xrightarrow{\sigma_{\tau}}  {}^{L} M(\Qlbar) \to  {}^{L} G(\Qlbar)$. 

Enfin elle est compatible aux cas triviaux de  fonctorialité.  
Soit $G'$ un  groupe réductif sur $K$ et $\Upsilon:G\to G'$ un morphisme de groupes sur $K$ dont  l'image est un sous-groupe distingué de $G'$. D'après le paragraphe 2.5 de \cite{borel-corvallis}, on en déduit  ${}^{L}\Upsilon: {}^{L}  G'\to {}^{L}  G$. Alors pour toute représentation lisse admissible irréductible $\pi$ de $G'(K)$, 
le centre de Bernstein de $G(K)$ agit sur $\pi$ par un caractère dont le paramètre local associé par \eqref{param-local-sigma-pi-loc} et a) ci-dessus  est  
${}^{L}\Upsilon\circ \sigma_{\pi}$. 
\end{thm}
\dem
La preuve du théorème est la même que dans les paragraphes~\ref{para-enonce} à~\ref{rem-preuve-corresp}
sauf qu'on avait utilisé \cite{cong-coeff-O}, qui n'est écrit que dans le cas  déployé, et qu'on utilise maintenant \cite{cong-zorro} qui traite le cas de tous les groupes. Pour montrer  dans le cas non déployé    l'assertion (vi)
 du  \thmref{thm-mfz}  (comme à la fin du  paragraphe~\ref{para-enonce})  on a besoin des morphismes   terme constant. Dans \cite{cong}, Cong Xue a défini les morphismes terme constant sur la cohomologie à support compact des champs de chtoucas pour les groupes déployés, et l'extension de leur construction au cas non déployé se fait par la m\^eme m\'ethode. Ils sont compatibles avec les morphismes de création  et d'annihilation. 
\cqfd

On en déduit une théorie des facteurs  $\gamma(s,\pi \times \chi,\psi) \in \mathbb{C}(T)$, où $K$ est un corps local de caractéristique $p$ contenant $\Fq$, $G$ est un groupe réductif   sur $K$, $\rho$ est une représentation finie de dimension $N$ de ${}^{L} G$ sur $ \Qlbar$, $\pi$ est une représentation lisse  irréductible  à coefficients dans  $\Qlbar $ de $G(K)$, $\chi:K^{*}\to \Qlbar^{*}$ est un caractère continu, 
et $\psi:K\to \Qlbar$ un caractère additif non dégénéré. Cette théorie des facteurs $\gamma$ vérifie les propriétés (i) à (vii) ci-dessus (en rempla\c cant 
$\wh M$ et $\wh G$ par ${}^{L}M$ et ${}^{L}G$ dans la propriété (iv)), et elle est caractérisée de manière unique par ces propriétés.

Comme dans le chapitre 14 de \cite{coh}, en utilisant \cite{finkelberg-lysenko} et 
\cite{dennis-sergey} on peut traiter le cas des groupes métaplectiques.

\end{document}